\documentclass
[11pt]{amsart}
\usepackage{amsmath,amsthm, amscd, amssymb, amsfonts}
\newtheorem{Lemma}{Lemma}
\newtheorem{Theorem}{Theorem}
\newtheorem{Prop}{Proposition}

\newtheorem{Cor}{Corollary}

\newcommand{\pf}{\medskip\noindent{\sc Proof: }}

\newcommand{\binr}[3]{\left(\!\! \begin{array}{c} #1 \\ #2 \end{array}\!\!\right)_{#3}}

\begin{document}
\title{On automorphisms of biproducts}
\author{David E. Radford}
\address{University of Illinois at Chicago \\
Department of Mathematics, Statistics and \\ Computer Science (m/c 240) \\
801 South Morgan Street \\
Chicago, IL   60608-7045} \email{radford@uic.edu}
\numberwithin{equation}{section}
\begin{abstract}
{\small  \rm We study certain subgroups of the full group of Hopf algebra automorphisms of a biproduct. In the process interesting subgroups of certain permutation groups come into play.}
\end{abstract}
\maketitle
\date{}

\setcounter{section}{0}
\section*{Introduction}\label{SecIntro}
Biproducts account for many examples of semisimple Hopf algebras over a field $k$. In \cite{RadProjNew} we studied a class of Hopf biproducts $A = B \times H$, where $B = k[\mathcal{G}]$ and $H = k[G]$ are group algebras, and described their coalgebra structures in detail. These biproducts are those of Theorem \ref{ThmFromProjNew}. In important special cases we were able to describe the algebra structure of $A$ in detail as well. A followup to \cite{RadProjNew} is the study Hopf algebra automorphisms of $A$.

Let $A, H$ be any Hopf algebras over $k$ and suppose $A = B \times H$ is a biproduct. The biproduct structure of $A$ is determined by Hopf algebra maps $\pi : A \longrightarrow H$ and $j : H \longrightarrow A$ which satisfy $\pi \circ j = \mathrm{Id}_H$. In this paper we will also think of biproducts as quadruples $(A, H, \pi, j)$, where $A, H$ are Hopf algebras over $k$ and $\pi, j$ are Hopf algebra maps as just described. A notion of automorphism of the biproduct $(A, H, \pi, j)$ is a Hopf algebra automorphism $F$ of $A$ which satisfies $\pi \circ F = \pi$ and $F \circ j = j$. These form a subgroup $\mathrm{Aut}_\mathrm{Hopf}(A, \pi, j)$ of the group of Hopf algebra automorphisms $\mathrm{Aut}_\mathrm{Hopf}(A)$ of $A$ under function composition.

Now $B$ is a Hopf algebra in the Yetter-Drinfel'd category ${}_H^H \mathcal{YD}$. We show that $\mathrm{Aut}_\mathrm{Hopf}(A, \pi, j) \simeq \mathrm{Aut}_{\mathcal{YD}-\mathrm{Hopf}}(B)$, the group of Hopf algebra automorphisms in the category ${}_H^H \mathcal{YD}$ in Theorem \ref{ThmYDAut}. Thus we find the two conditions $\pi \circ F = \pi$ and $F \circ j = j$ too restrictive. We relax the second.

The set $\mathrm{Aut}_\mathrm{Hopf}(A, \pi)$ of Hopf algebra automorphisms $F$ of $A$ satisfying $\pi \circ F = \pi$ is a subgroup of $\mathrm{Aut}_\mathrm{Hopf}(A)$. The study of $\mathrm{Aut}_\mathrm{Hopf}(A, \pi)$ is the focus of this paper.

Let $F \in \mathrm{Aut}_\mathrm{Hopf}(A, \pi)$. Then we define two maps $F_L : B \longrightarrow B$ and $F_R : H \longrightarrow B$ which provide a ``factorization" of $F$. It is always the case that $F_L$ is an algebra automorphism of $B$ and $\nu : \mathrm{Aut}_\mathrm{Hopf}(A, \pi) \longrightarrow \mathrm{Aut}_\mathrm{Alg}(B)$ defined by $\nu (F) = F_L$ is a group homomorphism.

Whether or not $F_L \in \mathrm{Aut}_{\mathcal{YD}_\mathrm{Hopf}}(B)$ is a major theme of this paper. If $\mathrm{Im}(\nu) \subseteq \mathrm{Aut}_{\mathcal{YD}_\mathrm{Hopf}}(B)$ then
\begin{equation}\label{EqIntroSemiDirectProd}
\mathrm{Aut}_\mathrm{Hopf}(A, \pi) \simeq \mathrm{Aut}_{\mathcal{YD}-\mathrm{Hopf}}(B) \times \mathcal{N}(B, H)
\end{equation}
as groups, where the right hand side has a semidirect product structure and $\mathcal{N}(B, H) \simeq \mathrm{Ker}(\nu)$. See Corollary \ref{CorAutSemiDirectProduct}.

Whether or not $F_L \in \mathrm{Aut}_{\mathcal{YD}_\mathrm{Hopf}}(B)$ comes down to whether or not $F_L$ is map of left $H$-modules and a map of coalgebras. Suppose that $B$ is commutative as a $k$-algebra and $H$ is cocommutative. Then $F_L$ is a module map. Hence in this case $F_L \in \mathrm{Aut}_{\mathcal{YD}_\mathrm{Hopf}}(B)$ if and only if $F_L$ is a coalgebra map. Whether or not $F_L$ is a coalgebra map is a fascinating and involved problem to investigate.

We are most interested in the biproducts $A = k[\mathcal{G}] \times k[G]$ of Theorem \ref{ThmFromProjNew} when $\mathcal{G}$ and $G$ are finite groups, $\mathcal{G}$ abelian. Some minor restrictions need to be placed on the field $k$. We show that whether or not the decomposition of $\mathrm{Aut}_\mathrm{Hopf}(A, \pi)$ described in (\ref{EqIntroSemiDirectProd}) holds boils down to the nature of a certain subgroup $\Gamma (\mathbf{G}, \lambda, \sigma)$, where $\mathbf{G}$ is the additive version of $\mathcal{G}$, of the group of permutations $\mathrm{Sym}(\mathbf{G})$ of the set $\mathbf{G}$ under function composition. See Section \ref{SecGammaGLambda}, Theorems \ref{ThmFromProjNew} and \ref{ThmMainOfPaper}. We use the standard notations $\leq$ and $<$ for ``is a subgroup of" and ``is a proper subgroup of" respectively.

The group $\Gamma (\mathbf{G}, \lambda, \sigma)$ is described in terms of certain elements of the character group $\widehat{G}$. This group should be of interest in its own right as well as a generalization studied in Sections \ref{SecAutThetaGGTheta} and \ref{SecAutThetaGGThetaEqual}.

Throughout $k$ is a field and all vector spaces are over $k$, though we use the redundant expression ``over $k$" quite often. The group of units of $k$ is denoted $k^\times$. For vector spaces $U$ and $V$ we drop the subscript $k$ from $\mathrm{Hom}_k(U, V)$, $\mathrm{End}_k(U)$, and $U \otimes_k V$. If $W$ is subspace of $U$ and $f : U \longrightarrow V$ is a map of vector spaces then $f|_W : W \longrightarrow V$ denotes the restriction of $f$ to $W$.

Let $A$ be an algebra. Then $\mathrm{Z}(A)$ denotes the center of $A$ and $1$ or $1_A$ denotes the multiplicative neutral element of $A$. Minor exception, $\mathrm{Id}_U$ denotes the identity map of a vector space $U$.

Let $C$ be a coalgebra. We use a shorthand version of the Heyneman-Sweedler notation for expressing the coproduct in writing $\Delta(c) = c_{(1)} \otimes c_{(2)}$ for $c \in C$. For a coalgebra $C$ and an algebra $D$ over $k$ we let ``$\star$" denote the convolution product of $\mathrm{Hom}(C, D)$. We suggest any of \cite{Abe,DRS,LamRad,Mont,RBook,SweedlerBook} as a Hopf algebra reference. Good general algebra references are \cite{DumFoote,Rotman2002}.

\section{Preliminaries}\label{SecPrelimninaries}
Let $H$ be a Hopf algebra over the field $k$. Then a biproduct of the form $A = B \times H$ is the tensor product $B \otimes H$ as a vector space. The vector space $B$ is an algebra, coalgebra, left $H$-module, and a left $H$-comodule over $k$. There are natural compatibilities among these structures which we will list for the reader's convenience. As an algebra $A = B \# H$ is a smash product and as a coalgebra $A = B \natural H$ is a smash coproduct. The compatibilities are expressed by $B$ is a left $H$-module algebra and coalgebra, by $B$ is a left $H$-comodule algebra and coalgebra, and by connections between the algebra and coalgebra structures of $B$ and the module and comodule structures on $B$. If $C$ is a coalgebra and $(M, \rho)$ is a left $C$-comodule then we write $\rho(m) = m_{(-1)} \otimes m_{(0)} \in C \otimes M$ for all $m \in M$.

A left $H$-module algebra is a left $H$-module $(B, \cdot)$, where $B$ is an algebra over $k$, such that $h{\cdot}1 =\epsilon(h)1$ and $h{\cdot}(bb') = (h_{(1)}{\cdot}b)(h_{(2)}{\cdot}b')$ for all $h \in H$ and $b, b' \in B$. If $(B, \cdot)$ is a left $H$-module algebra the tensor product $B \otimes H$ of vector spaces has an algebra structure, referred to as the smash product, defined by $1_{B \otimes H} = 1_B \otimes 1_H$ and
\begin{equation}\label{EqSmash}
(b \otimes h)(b' \otimes h') = b(h_{(1)}{\cdot}b') \otimes h_{(2)}h'
\end{equation}
for all $b, b' \in B$ and $h, h' \in H$. Typical notation for this algebra is $B \# H$ and tensors $b \otimes h$ are written $b \# h$. Observe $(b \# h)(b' \# h') = bb' \# hh'$ if either $h = 1$ or $b' = 1$. As a consequence the maps $H \longrightarrow B \# H$ and $B \longrightarrow B \# H$ defined by $h \mapsto 1 \# h$ and $b \mapsto b \# 1$ respectively are a one-one algebra maps.

A left $H$-comodule algebra is a left $H$-comodule $(B, \rho)$, where $B$ is an algebra over $k$, such that $\rho(1) = 1_H \otimes 1$ and $\rho(bb') = b_{(-1)}b'_{(-1)} \otimes b_{(0)}b'_{(0)}$ for all $b, b' \in B$. A left $H$-module coalgebra is a left $H$-module $(C, \cdot)$, where $C$ is a coalgebra over $k$, such that $\epsilon (h{\cdot}c) = \epsilon(h)\epsilon(c)$ and $\Delta (h{\cdot}c) = h_{(1)}{\cdot}c_{(1)} \otimes h_{(2)}{\cdot}c_{(2)}$ for all $h \in H$ and $c \in C$.

A left $H$-comodule coalgebra is a left $H$-comodule $(C, \rho)$, where $C$ is a coalgebra over $k$, such that $c_{(-1)}\epsilon(c_{(0)}) = \epsilon(c)1$ and
$$
c_{(1)(-1)}c_{(2)(-1)} \otimes c_{(1)(0)} \otimes c_{(2)(0)} = c_{(-1)} \otimes c_{(0)(1)} \otimes c_{(0)(2)}
$$
for all $c \in C$. If $(C, \rho)$ is a left $H$-comodule coalgebra then the tensor product $C \otimes H$ of vector spaces has a coalgebra structure, referred to as the smash coproduct, defined by $\epsilon_{C \otimes H} = \epsilon_C \otimes \epsilon_H$ and
\begin{equation}\label{EqCoSmash}
\Delta (c \otimes h) = (c_{(1)} \otimes c_{(2)(-1)}h_{(1)}) \otimes (c_{(2)(0)} \otimes h_{(2)})
\end{equation}
for all $c \in C$ and $h \in H$. The usual notation for this coalgebra is $C \natural H$. The tensor $c \otimes h$ is written $c \natural h$.

Let $A = B \times H$ be a biproduct. Generally $B$ with its algebra and coalgebra structure is not a bialgebra over $k$. However:
\begin{equation}\label{EqnEpsilonMult}
\epsilon \in \mathrm{Alg}(B, k), \qquad \Delta(1) = 1
\end{equation}
and
\begin{equation}\label{EqDetltaProd}
\Delta (bb') = b_{(1)}(b_{(2)(-1)}{\cdot}b'_{(1)}) \otimes b_{(2)(0)}b'_{(2)}
\end{equation}
for all $b, b' \in B$. When $B$ is a left $H$-module algebra and coalgebra, is a left $H$-comodule algebra and coalgebra, (\ref{EqnEpsilonMult}) and (\ref{EqDetltaProd}) are satisfied,
\begin{equation}\label{EqYDCondition}
h_{(1)}m_{(-1)} \otimes h_{(2)}{\cdot}m_{(0)} = (h_{(1)}{\cdot}m)_{(-1)}h_{(2)} \otimes (h_{(1)}{\cdot}m)_{(0)}
\end{equation}
holds for all $h \in H$ and $m \in B$, and $\mathrm{Id}_B$ has a convolution inverse, then the vector space $A = B \otimes H$ is a Hopf algebra with the smash product algebra structure and the smash coproduct coalgebra structure. The Hopf algebra $A$ is called the biproduct of $B$ and $H$ and is denoted $A = B \times H$. Tensors $b \otimes h$ are denoted $b \times h$.

Let $A = B \times H$ be the biproduct of $B$ and $H$. Define $\pi : A \longrightarrow H$ by $\pi(b \times h) = \epsilon(b)h$ for $b \in B$, $h \in H$ and $j : H \longrightarrow A$ by $j(h) = 1 \times h$ for $h \in H$. Observe that $\pi$ and $j$ are Hopf algebra maps and $\pi \circ j = \mathrm{Id}_H$. The map $\pi$ is referred to as a Hopf algebra projection of $A$ onto $H$.

Conversely, if $A$ and $H$ are Hopf algebras over $k$ and $\pi : A \longrightarrow H$, $j : H \longrightarrow A$ are Hopf algebra maps which satisfy $\pi \circ j = \mathrm{Id}_H$ then $A \simeq B \times H$ for some $B$. See \cite{RadProj} which covers all the details for this section.

Throughout this paper biproducts will be Hopf algebras over $k$.

\section{Certain Categories Associated with Biproducts}\label{SecBiProdCat}
Let $H$ be a Hopf algebra over the field $k$. We denote by $\mathcal{C}_\mathrm{proj}$ the category whose objects are quadruples $(A, H, \pi, j)$, where $A$ is a Hopf algebra over $k$ and $\pi : A \longrightarrow H$, $j : H \longrightarrow A$ are Hopf algebra maps which satisfy $\pi \circ j = \mathrm{Id}_H$. Morphisms are pairs $(F, f) : (A, H, \pi, j) \longrightarrow (A', H', \pi', j')$, where $F : A \longrightarrow A'$ and $f : H \longrightarrow H'$ are Hopf algebra maps such that the diagrams
\begin{equation}\label{EqSquareCommute}
\mbox{
\begin{picture}(100,100)(0,0)
\put(-4, 0){$H$}
\put(-4, 80){$A$}
\put(98, 0){$H'$}
\put(98, 80){$A'$}
\put(-15,40){$\pi$}
\put(110,40){$\pi'$}
\put(45,-10){$f$}
\put(44,90){$F$}
\put(0, 75){\vector(0,-1){60}}
\put(10, 3){\vector(1, 0){80}}
\put(10, 83){\vector(1, 0){80}}
\put(102, 75){\vector(0, -1){60}}
\end{picture}
\qquad
\raisebox{40pt}{and}
\qquad
\begin{picture}(100,100)(0,0)
\put(-4, 0){$H$}
\put(-4, 80){$A$}
\put(98, 0){$H'$}
\put(98, 80){$A'$}
\put(-15,40){$j$}
\put(110,40){$j'$}
\put(45,-10){$f$}
\put(44,90){$F$}
\put(0, 14){\vector(0,1){60}}
\put(10, 3){\vector(1, 0){80}}
\put(10, 83){\vector(1, 0){80}}
\put(102, 14){\vector(0, 1){60}}
\end{picture}
}
\end{equation}
\medskip

\noindent
commute. The composition of morphisms is given by the composition of components. For an object $(A, H, \pi, j)$ of $\mathcal{C}_\mathrm{proj}$ we let $\mathrm{End}_\mathrm{Hopf}(A, H, \pi, j)$ be the set of all Hopf algebra endomorphisms $F$ of $A$ such that the diagrams
\begin{equation}\label{EqTriangleCommute}
\mbox{
\begin{picture}(100,120)(0,0)
\put(38, 0){$H$}
\put(-2, 78){$A$}
\put(78, 78){$A$}
\put(10, 40){$\pi$}
\put(68,40){$\pi$}
\put(40, 88){$F$}
\put(8, 72){\vector(1,-2){30}}
\put(78, 72){\vector(-1,-2){30}}
\put(10, 80){\vector(1,0){60}}
\end{picture}
\qquad
\raisebox{40pt}{and}
\qquad
\begin{picture}(100,120)(0,0)
\put(38, 0){$H$}
\put(-2, 78){$A$}
\put(78, 78){$A$}
\put(10, 40){$j$}
\put(68,40){$j$}
\put(40, 88){$F$}
\put(38, 12){\vector(-1,2){30}}
\put(48, 12){\vector(1,2){30}}
\put(10, 80){\vector(1,0){60}}
\end{picture} }
\end{equation}
\medskip

\noindent
commute. Thus $\mathrm{End}_\mathrm{Hopf}(A, H, \pi, j)$ is the set of all $F$ such that $(F, \mathrm{Id}_H)$ is an endomorphism of $(A, H, \pi, j)$. Observe that $\mathrm{End}_\mathrm{Hopf}(A, H, \pi, j)$ is a monoid under function composition. The units of this monoid form the set $\mathrm{Aut}_\mathrm{Hopf}(A, H, \pi, j)$ of all Hopf algebra automorphisms $F$ of $A$ such that the diagrams of (\ref{EqTriangleCommute}) commute. It is a subgroup of the group $\mathrm{Aut}_\mathrm{Hopf}(A)$ of all Hopf algebra automorphisms of $A$ group under composition.

As it turns out $\mathrm{Aut}_\mathrm{Hopf}(A, H, \pi, j)$ is a bit too special. We will consider a less restrictive setting in our study of endomorphisms of biproducts by ignoring the map $j$.

Let $\mathcal{C}'_\mathrm{proj}$ be the category whose objects are triples $(A, H, \pi)$ which come from objects $(A, H, \pi, j)$ of $\mathcal{C}_\mathrm{proj}$. Morphisms of the category $\mathcal{C}'_\mathrm{proj}$ are pairs $(F, f) : (A, H, \pi) \longrightarrow (A', H', \pi')$, where $F : A \longrightarrow A'$ and $f : H \longrightarrow H'$ are Hopf algebra maps such that the first diagram of (\ref{EqSquareCommute}) commutes. Again, composition of morphisms is given by composition of components.

For and object $(A, H, \pi)$ of $\mathcal{C}'_\mathrm{proj}$ let $\mathrm{End}_\mathrm{Hopf}(A, H, \pi)$ be the monoid of all Hopf algebra endomorphisms $F$ of $A$ such that the first diagram of (\ref{EqTriangleCommute}) commutes and let $\mathrm{Aut}_\mathrm{Hopf}(A, H, \pi)$  be its set of units. Thus $\mathrm{Aut}_\mathrm{Hopf}(A, H, \pi)$ is the group of Hopf algebra automorphisms $F$ of $A$ such that $\pi \circ F = \pi$ under composition.

In the context of a fixed $H$ we will simplify symbolism by dropping the $H$ in the notation for the two monoids and two groups defined above. Thus we will write $\mathrm{Aut}_\mathrm{Hopf}(A, \pi, j)$ for $\mathrm{Aut}_\mathrm{Hopf}(A, H, \pi, j)$ and $\mathrm{Aut}_\mathrm{Hopf}(A, \pi)$ for $\mathrm{Aut}_\mathrm{Hopf}(A, H, \pi)$.

This paper characterizes elements of $\mathrm{End}_\mathrm{Hopf}(A, \pi)$; in particular elements of $\mathrm{Aut}_\mathrm{Hopf}(A, \pi)$. Some of them arise from the (left-left) Yetter-Drinfel'd category ${}_H^H\mathcal{YD}$ in terms of which $A = B \times H$ is naturally understood.

The objects of ${}_H^H\mathcal{YD}$ are triples $(M, \cdot, \rho)$, where $(M, \cdot)$ is a left $H$-module and $(M, \rho)$ is a left $H$-comodule, such that (\ref{EqYDCondition}) holds
for all $h \in H$ and $m \in M$. Morphisms are functions $f : M \longrightarrow N$ of objects which are left $H$-module and left $H$-comodule maps. Multiplication of morphisms is function composition. If $A = B \times H$ is a biproduct which is a Hopf algebra over $k$ then $B$ is a Hopf algebra of ${}_H^H\mathcal{YD}$. Conversely, if $B$ is a Hopf algebra of ${}_H^H\mathcal{YD}$ the biproduct $A = B \times H$ is defined and is a Hopf algebra over $k$.

The category ${}_H^H\mathcal{YD}$ plays a minor role in this paper. Let $B$ be a Hopf algebra of ${}_H^H\mathcal{YD}$. For Theorem \ref{ThmYDAut} we will need a description of the monoid of Hopf algebra endomorphisms $\mathrm{End}_{\mathcal{YD}-\mathrm{Hopf}}(B)$. Recall $B$ is an algebra and coalgebra over $k$. The elements of $\mathrm{End}_{\mathcal{YD}-\mathrm{Hopf}}(B)$ are those linear endomorphisms of $B$ which are algebra and coalgebra maps and also maps of left $H$-modules and comodules. $\mathrm{Aut}_{\mathcal{YD}-\mathrm{Hopf}}(B)$ is the group of units of this monoid under function composition.

The reader is encouraged to consult basic references \cite{Majid1993Skly,Majid1994Braided,Yetter1990} for Yetter-Drinfel'd categories. See \cite{RadProjNew,Rad-Tow1993} also.
\section{Factorization of Certain Biproduct Endomorphisms}\label{SecAutoFactorize}
Throughout this section $H$ is a Hopf algebra with antipode $S$ and $A = B \times H$ is a biproduct which is a Hopf algebra over $k$. Let $F \in \mathrm{End}_\mathrm{Hopf}(A, \pi)$. The purpose of this section is to show $F$ has a factorization closely related to the factors $B$ and $H$.

Recall from Section \ref{SecPrelimninaries} that $\pi : A \longrightarrow H$ defined by $\pi(b \times h) = \epsilon(b)h$ for $b \in B$, $h \in H$ and $j : H \longrightarrow A$ defined by $j(h) = 1 \times h$ for $h \in H$ are Hopf algebra maps which satisfy $\pi \circ j = \mathrm{Id}_H$. Likewise we define $\Pi : A \longrightarrow B$ and $J : B \longrightarrow A$ by $\Pi(b \times h) = b\epsilon(h)$ for $b \in B$, $h \in H$ and $J(b) = b \times 1$ for $b \in B$. Note that $\Pi \circ J = \mathrm{Id}_B$.  There is a fundamental relationship between these four maps given by
\begin{equation}\label{EqJPiIsjSpi}
J \circ \Pi = \mathrm{Id}_A\star(j \circ S \circ \pi).
\end{equation}

The factorization of $F$ is given in terms of $F_L : B \longrightarrow B$ and $F_R : H \longrightarrow B$ defined by
\begin{equation}\label{EqFLFRDef}
F_L = \Pi \circ F \circ J \;\; \mbox{and} \;\; F_R = \Pi \circ F \circ j.
\end{equation}
\begin{Lemma}\label{LemmaFLFR}
Let $F \in \mathrm{End}_\mathrm{Hopf}(A, \pi)$. Then:
\begin{enumerate}
\item[{\rm (a)}] $F_L(b) \times 1 = F(b \times 1)$ for all $b \in B$.
\item[{\rm (b)}] $F_R(h) \times 1 = F(1 \times h_{(1)})(1 \times S(h_{(2)}))$ for all $h \in H$.
\item[{\rm (c)}] $F(b \times h) = F_L(b)F_R(h_{(1)}) \times h_{(2)}$ for all $b \in B$ and $h \in H$.
\end{enumerate}
\end{Lemma}

\pf
We need to calculate $J  \circ \Pi \circ F$. Let $b \in B$ and $h \in H$. We use (\ref{EqJPiIsjSpi}) to compute
\begin{eqnarray*}
(J \circ \Pi)(F(b \times h)) & = & F((b \times h)_{(1)})((j \circ S \circ \pi)(F((b\times h)_{(2)}))) \\
& = &  F((b \times h)_{(1)})((j \circ S \circ \pi)((b\times h)_{(2)})) \\
& = & F(b_{(1)} \times b_{(2)(-1)}h_{(1)})((j \circ S \circ \pi)(b_{(2)(0)} \times h_{(2)})) \\
& = & F(b_{(1)} \times b_{(2)(-1)}\epsilon(b_{(2)(0)})h_{(1)})(1 \times S(h_{(2)})) \\
& = & F(b_{(1)} \times \epsilon(b_{(2)})h_{(1)})(1 \times S(h_{(2)})) \\
& = & F(b \times h_{(1)})(1 \times S(h_{(2)})).
\end{eqnarray*}
Thus
$$
(J \circ \Pi \circ F)(b \times h) = F(b \times h_{(1)})(1 \times S(h_{(2)}))
$$
for all $b \in B$ and $h \in H$. Parts (a) and (b) follow from these equations. As for part (c), we use parts (a) and (b) to calculate
\begin{eqnarray*}
F(b \times h) & = & F((b \times 1)(1 \times h)) \\
& = & F(b \times 1)F(1 \times h) \\
& = & F(b \times 1)F(1 \times h_{(1)})(1 \times S(h_{(2)}))(1 \times h_{(3)}) \\
& = & (F_L(b) \times 1)(F_R(h_{(1)}) \times 1)(1 \times h_{(2)}) \\
& = & F_L(b)F_R(h_{(1)}) \times h_{(2)}
\end{eqnarray*}
for all $b \in B$ and $h \in H$.
\qed
\medskip

By parts (a) and (b) of the preceding lemma:
\begin{equation}\label{EqnIdASubL}
(\mathrm{Id}_A)_L = \mathrm{Id}_B \;\; \mbox{and} \;\; (\mathrm{Id}_A)_R = \eta \circ \epsilon.
\end{equation}
Since $F_L(1) = 1$ by part (a) of Lemma \ref{LemmaFLFR}, by part (c) of the same:
\begin{equation}\label{EqFOneTimesH}
F(1 \times h) = F_R(h_{(1)}) \times h_{(2)}
\end{equation}
for all $h \in H$. We are now able to compute the factors of a composite.
\begin{Cor}\label{CorFGParts}
Let $F, G \in \mathrm{End}_\mathrm{Hopf}(A, \pi)$. Then:
\begin{enumerate}
\item[{\rm (a)}] $(F \circ G)_L = F_L \circ G_L$.
\item[{\rm (b)}] $(F \circ G)_R = (F_L \circ G_R) \star F_R$.
\end{enumerate}
\end{Cor}

\pf
Let $b \in B$. Then $(F \circ G)_L(b) \times 1 = F_L(G_L(b)) \times 1$ follows by part (a) of Lemma \ref{LemmaFLFR}. Part (a) is established. Let $h \in H$. Using (\ref{EqFOneTimesH}), parts (a) and (c) of Lemma \ref{LemmaFLFR}, and the fact the $F$ is multiplicative, we note that
\begin{eqnarray*}
(F \circ G)_R(h_{(1)}) \times h_{(2)} & = & (F \circ G)(1 \times h) \\
& = & F(G(1 \times h)) \\
& = & F(G_R(h_{(1)}) \times h_{(2)}) \\
& = & F(G_R(h_{(1)}) \times 1)F(1 \times h_{(2)}) \\
& = & (F_L(G_R(h_{(1)})) \times 1)(F_R(h_{(2)}) \times h_{(3)}) \\
& = & F_L(G_R(h_{(1)}))F_R(h_{(2)}) \times h_{(3)}
\end{eqnarray*}
from which $(F \circ G)_R(h) = F_L(G_R(h_{(1)}))F_R(h_{(2)})$ follows. We have shown part (b) holds.
\qed
\medskip

By virtue of Lemma \ref{LemmaFLFR}, to characterize $F$ is a matter of characterizing $F_L$ and $F_R$. Note in particular part (e) of the following describes a commutation relation between $F_L$ and $F_R$.
\begin{Lemma}\label{LemmaFLStructure}
Let $F \in \mathrm{End}_\mathrm{Hopf}(A, \pi)$. Then:
\begin{enumerate}
\item[{\rm (a)}] $F_L : B \longrightarrow B$ is an algebra endomorphism.
\item[{\rm (b)}] $\epsilon \circ F_L = \epsilon$.
\item[{\rm (c)}] $\Delta(F_L(b)) = F_L(b_{(1)})F_R(b_{(2)(-1)}) \otimes F_L(b_{(2)(0)})$ for all $b \in B$.
\item[{\rm (d)}] $\rho(F_L(b)) = b_{(-1)} \otimes F_L(b_{(0)})$ for all $b \in B$.
\item[{\rm (e)}] $F_L(h_{(1)}{\cdot}b)F_R(h_{(2)}) = F_R(h_{(1)})(h_{(2)}{\cdot}F_L(b))$ for all $h \in H$ and $b \in B$.
\end{enumerate}
\end{Lemma}

\pf
We first show part (a). That $F(B \times 1) \subseteq B \times 1$ follows by part (a) of Lemma \ref{LemmaFLFR}. Let $F_r : B \times 1 \longrightarrow B \times 1$ be the map induced by restriction of $F$. The map $J$ thought of as $J : B \longrightarrow B \times 1$ is an algebra isomorphism. Thus $F_L = J^{-1} \circ F_r \circ J$ is an algebra map since it is the composition of such. Part (b) follows since $\epsilon \circ J = \epsilon$ and $F$ is a coalgebra map. Let $b \in B$. To show part (c) we compute the coproduct of $F_L(b) \times 1 = F(b \times 1)$ in two ways. First of all
$$
\Delta(F_L(b) \times 1) = (F_L(b)_{(1)} \times F_L(b)_{(2)(-1)}) \otimes (F_L(b)_{(2)(0)} \times 1)
$$
and secondly, since $F$ is a coalgebra map,
\begin{eqnarray*}
\Delta(F(b \times 1)) & = & F((b \times 1)_{(1)}) \otimes F((b \times 1)_{(2)}) \\
& = & F(b_{(1)} \times b_{(2)(-1)}) \otimes F(b_{(2)(0)} \times 1) \\
& = & (F_L(b_{(1)})F_R(b_{(2)(-1)(1)}) \times b_{(2)(-1)(2)}) \otimes (F_L(b_{(2)(0)}) \times 1);
\end{eqnarray*}
the last equation follows by part (c) of Lemma \ref{LemmaFLFR}. We obtain part (c) by applying $\epsilon$ to the second and fourth factors to the two expressions for the coproduct. Noting that $\epsilon \circ F_R = \epsilon$, which follows by (\ref{EqFOneTimesH}), we use part (b) to obtain part (d) by applying $\epsilon$ to the first and fourth factors of the same coproduct expressions.

It remains to show part (e). For $b \in B$ and $h \in H$ we have
$$
F((1 \times h)(b \times 1)) = F(h_{(1)}{\cdot}b \times h_{(2)}) = F_L(h_{(1)}{\cdot}b)F_R(h_{(2)}) \times h_{(3)}
$$
by part (c) of Lemma \ref{LemmaFLFR}. On the other hand, since $F$ is an algebra map we use (\ref{EqFOneTimesH}) and part (a) of Lemma \ref{LemmaFLFR} to compute
\begin{eqnarray*}
F((1 \times h)(b \times 1)) & = & F(1 \times h)F(b \times 1)) \\
&  = & (F_R(h_{(1)}) \times h_{(2)})(F_L(b) \times 1) \\
& = & F_R(h_{(1)})(h_{(2)}{\cdot}F_L(b)) \times h_{(3)}.
\end{eqnarray*}
Applying $\epsilon$ to the second factor of both expressions for $F((1 \times h)(b \times 1))$ we obtain part (e).
\qed
\medskip

As the reader might suspect, whether or not $F_L$ is a coalgebra map is explained in terms of $F_R$.
\begin{Cor}\label{CorFLCoalgebraMsap}
Let $F \in \mathrm{End}_\mathrm{Hopf}(A, \pi)$. Then $F_L$ is a coalgebra map if and only if $F_R(c_{(-1)}) \otimes c_{(0)} = 1 \otimes c$ for all $c \in \mathrm{Im}(F_L)$.
\end{Cor}

\pf
Suppose $F_R(c_{(-1)}) \otimes c_{(0)} = 1 \otimes c$ for all $c \in \mathrm{Im}(F_L)$. Then by parts (c) and (d) of Lemma \ref{LemmaFLStructure} we have
\begin{eqnarray*}
\Delta(F_L(b)) & = & F_L(b_{(1)})F_R(b_{(2)(-1)}) \otimes F_L(b_{(2)(0)}) \\
& = & F_L(b_{(1)})F_R(F_L(b_{(2)})_{(-1)}) \otimes F_L(b_{(2)})_{(0)} \\
& = & F_L(b_{(1)})1 \otimes F_L(b_{(2)})
\end{eqnarray*}
for all $b \in B$. This calculation and part (b) of Lemma \ref{LemmaFLStructure} imply that $F_L$ is a coalgebra map.

Conversely, suppose that $F_L$ is a coalgebra map. Using parts (a), (c), and (d) of Lemma \ref{LemmaFLStructure} we compute for all $b \in B$ that
\begin{eqnarray*}
F_R(F_L(b)_{(-1)}) \otimes F_L(b)_{(0)} & = & F_L(\epsilon(b_{(1)})1)F_R(F_L(b_{(2)})_{(-1)}) \otimes F_L(b_{(2)})_{(0)} \\
& = & F_L(S(b_{(1)})b_{(2)})F_R(b_{(3)(-1)}) \otimes F_L(b_{(3)(0)}) \\
& = & F_L(S(b_{(1)}))F_L(b_{(2)})F_R(b_{(3)(-1)}) \otimes F_L(b_{(3)(0)}) \\
& = & F_L(S(b_{(1)}))F_L(b_{(2)})_{(1)} \otimes F_L(b_{(2)})_{(2)} \\
& = & F_L(S(b_{(1)}))F_L(b_{(2)}) \otimes F_L(b_{(3)}) \\
& = & F_L(S(b_{(1)})b_{(2)}) \otimes F_L(b_{(3)}) \\
& = & F_L(1) \otimes F_L(b) \\
& = & 1 \otimes F_L(b).
\end{eqnarray*}
We have shown that $F_R(c_{(-1)}) \otimes c_{(0)} = 1 \otimes c$ for all $c \in \mathrm{Im}(F_L)$.
\qed
\medskip

\begin{Lemma}\label{LemmaFRStructure}
Let $F \in \mathrm{End}_\mathrm{Hopf}(A, \pi)$. Then:
\begin{enumerate}
\item[{\rm (a)}] $F_R(hh') = F_R(h_{(1)})(h_{(2)}{\cdot}F_R(h'))$ for all $h, h' \in H$.
\item[{\rm (b)}] $F_R(1) = 1$.
\item[{\rm (c)}] $F_R : H \longrightarrow B$ is a coalgebra map.
\item[{\rm (d)}] $\rho(F_R(h)) = h_{(1)}S(h_{(3)}) \otimes F_R(h_{(2)})$ for all $h \in H$.
\end{enumerate}
\end{Lemma}

\pf
By (\ref{EqFOneTimesH}) we have $F(1 \times h) = F_R(h_{(1)}) \times h_{(2)}$ for all $h \in H$. Since $F$ is an algebra map, $1 \times 1 = F(1 \times 1) = F_R(1) \times 1$ which implies $F_R(1) = 1$. We have established part (b). As for part (a), for $h, h' \in H$ we calculate on one hand
$$
F(1 \times hh') = F_R(h_{(1)}h'_{(1)}) \times h_{(2)}h'_{(2)}
$$
and on the other
\begin{eqnarray*}
F(1 \times hh') & = & F((1 \times h)(1 \times h')) \\
& = & F(1 \times h)F(1 \times h') \\
& = & (F_R(h_{(1)}) \times h_{(2)})(F_R(h'_{(1)}) \times h'_{(2)}) \\
& = & F_R(h_{(1)})(h_{(2)}{\cdot}F_R(h'_{(1)})) \times h_{(3)}h'_{(2)}.
\end{eqnarray*}
Applying $\epsilon$ to the second factor of both expressions for $F(1 \times hh')$ establishes part (a).

Let $h \in H$. To show parts (c) and (d) we compute $\Delta(F(1 \times h))$ in two ways. Since $F$ and $j$ are coalgebra maps
$$
\Delta (F(1 \times h)) = F(1 \times h_{(1)}) \otimes F(1 \times h_{(2)}) = (F_R(h_{(1)}) \times h_{(2)}) \otimes (F_R(h_{(3)}) \times h_{(4)}).
$$
On the other hand
\begin{eqnarray*}
\lefteqn{\Delta(F(1 \times h))} \\
& = & \Delta(F_R(h_{(1)}) \times h_{(2)}) \\
& = & (F_R(h_{(1)})_{(1)} \times F_R(h_{(1)})_{(2)(-1)}h_{(2)}) \otimes (F_R(h_{(1)})_{(2)(0)} \times h_{(3)}).
\end{eqnarray*}
Applying $\epsilon$ to the second and fourth factors of the expressions for $\Delta (F(1 \times h))$ gives $F_R(h_{(1)}) \otimes F_R(h_{(2)}) = \Delta (F_R(h))$ and, since $\epsilon \circ F_R = \epsilon$, to the first and fourth gives $h_{(1)} \otimes F_R(h_{(2)}) = F_R(h_{(1)})_{(-1)}h_{(2)} \otimes F_R(h_{(1)})_{(0)}$. Therefore
$$
\rho(F_R(h)) = F_R(h_{(1)})_{(-1)}h_{(2)}S(h_{3}) \otimes F_R(h_{(1)})_{(0)} = h_{(1)}S(h_{(3)}) \otimes F_R(h_{(2)}).
$$
We have established parts (c) and (d).
\qed
\begin{Cor}\label{CorFLModuleMap}
Let $F \in \mathrm{End}_\mathrm{Hopf}(A, \pi)$. Then:
\begin{enumerate}
\item[{\rm (a)}] $F_L$ is a left $H$-module map if and only if the condition $$F_L(h_{(1)}{\cdot}b)F_R(h_{(2)}) = F_R(h_{(1)})F_L(h_{(2)}{\cdot}b)$$
holds for all $h \in H$ and $b \in B$.
\item[{\rm (b)}]  If $B$ is commutative and $H$ is cocommutative then $F_L$ is a left $H$-module map.
\end{enumerate}
\end{Cor}

\pf
Part (b) follows immediately from part (a). To show part (a), we first note if $F_L$ is a map of left $H$-modules then the condition follows by part (e) of Lemma \ref{LemmaFLStructure}. Suppose the condition holds. Now $F_R$ is a coalgebra map by part (c) of Lemma \ref{LemmaFRStructure}. Using this fact and part (e) of Lemma \ref{LemmaFLStructure}, observe for all $h \in H$ and $b \in B$ that
\begin{eqnarray*}
h{\cdot} F_L(b) & = & \epsilon(F_R(h_{(1)}))(h_{(2)}{\cdot}F_L(b)) \\
& = & S(F_R(h_{(1)}))F_R(h_{(2)})(h_{(3)}{\cdot}F_L(b)) \\
& = & S(F_R(h_{(1)}))F_L(h_{(2)}{\cdot}b)F_R(h_{(3)}) \\
& = & S(F_R(h_{(1)}))F_R(h_{(2)})F_L(h_{(3)}{\cdot}b) \\
& = & \epsilon(F_R(h_{(1)}))F_L(h_{(2)}{\cdot}b) \\
& = & F_L(h{\cdot}b)
\end{eqnarray*}
which shows that $F_L$ is a left $H$-module map.
\qed
\begin{Cor}\label{CorFRAlgebraMap}
Let $F \in \mathrm{End}_\mathrm{Hopf}(A, \pi)$. Then $F_R$ is an algebra map if and only if $h{\cdot}F_R(h') = \epsilon(h)F_R(h')$ for all $h, h' \in H$.
\end{Cor}

\pf
First of all $F_R(1) = 1$ by part (b) of Lemma \ref{LemmaFRStructure}. If the condition $h{\cdot}F_R(h') = \epsilon(h)F_R(h')$ holds for all $h, h' \in H$ then $F_R(hh') = F_R(h)F_R(h')$ for all $h, h' \in H$ by part (a) of Lemma \ref{LemmaFRStructure}. Thus $F_R$ is an algebra map.

Suppose $F_R$ is an algebra map and let $h, h' \in H$. Now $F_R$ is a coalgebra map by part (c) of Lemma \ref{LemmaFRStructure}. Using this fact and part (a) of the same again
\begin{eqnarray*}
\epsilon(h)F_R(h') & = & S(F_R(h_{(1)}))F_R(h_{(2)})F_R(h') \\
& = & S(F_R(h_{(1)}))F_R(h_{(2)}h') \\
& = & S(F_R(h_{(1)}))F_R(h_{(2)})(h_{(3)}{\cdot}F_R(h')) \\
& = & \epsilon(h_{(1)})h_{(2)}{\cdot}F_R(h') \\
& = & h{\cdot}F_R(h').
\end{eqnarray*}
\qed

\begin{Cor}\label{CorhFRhPrime}
Let $F \in \mathrm{End}_\mathrm{Hopf}(A, \pi)$. Then $F_R$ has a convolution inverse $J_R$ defined by $J_R(h) = h_{(1)}{\cdot}F_R(S(h_{(2)}))$ for all $h \in H$.
\end{Cor}

\pf Let $h \in H$. Then by parts (a) and (b) of Lemma \ref{LemmaFRStructure} we have
\begin{eqnarray*}
F_R\star J_R(h) & = & F_R(h_{(1)})(h_{(2)}{\cdot}F_R(S(h_{(3)}))) \\
& = & F_R(h_{(1)}S(h_{(2)})) \\
& = & F_R(\epsilon(h)1) \\
& = & \epsilon(h)1
\end{eqnarray*}
and using the fact that $B$ is a left $H$-module algebra we have
\begin{eqnarray*}
J_R \star F_R(h) & = & (h_{(1)}{\cdot}F_R(S(h_{(2)})))F_R(h_{(3)}) \\
& = & h_{(1)}{\cdot}(F_R(S(h_{(3)}))(S(h_{(2)}){\cdot}F_R(h_{(4)}))) \\
& = & h_{(1)}{\cdot}(F_R(S(h_{(2)})_{(1)})(S(h_{(2)})_{(2)}{\cdot}F_R(h_{(3)}))) \\
& = & h_{(1)}{\cdot}(F_R(S(h_{(2)})h_{(3)})) \\
& = & h_{(1)}{\cdot}(F_R(\epsilon(h)1)) \\
& = & h{\cdot}F_R(1) \\
& = & h{\cdot}1 \\
& = & \epsilon(h)1.
\end{eqnarray*}
\qed
\medskip

We now characterize $\mathrm{End}_\mathrm{Hopf}(A, \pi)$ and $\mathrm{Aut}_\mathrm{Hopf}(A, \pi)$.
\begin{Theorem}\label{ThmEndAut}
Let $A = B \times H$ be a biproduct and $\pi : A \longrightarrow$ be the projection from $A$ onto $H$ and let $\mathcal{F}_{B, H}$ be the set of pairs $(\textit{f}, \textit{g})$, where $\textit{f} : B \longrightarrow B$ and $\textit{g} : H \longrightarrow B$ are maps satisfy the conclusions of Lemma \ref{LemmaFLStructure} and Lemma \ref{LemmaFRStructure} for $F_L$ and $F_R$ respectively. Then:
\begin{enumerate}
\item[{\rm (a)}] The function $\mathcal{F}_{B, H} \longrightarrow  \mathrm{End}_\mathrm{Hopf}(A, \pi)$, described by $(\textit{f}, \textit{g}) \mapsto F$, where $F(b \times h) = \textit{f} \, (b)\textit{g}(h_{(1)}) \times h_{(2)}$ for all $b \in B$ and $h \in H$, is a bijection. Furthermore $F_L = \textit{f}$ and $F_R = \textit{g}$.
\item[{\rm (b)}] Suppose $(\textit{f}, \textit{g}) \in \mathcal{F}_{B, H}$. Then $F \in \mathrm{Aut}_\mathrm{Hopf}(A, \pi)$ if and only if $\textit{f}$ is a bijection.
\end{enumerate}
\end{Theorem}

\pf
Assume the function of part (a) is well-defined. We first observe that $F \mapsto (\Pi \circ F \circ J, \, \Pi \circ F \circ j)$ describes its inverse; see (\ref{EqFLFRDef}). In light of the preceding results, to complete the proof of part (a) we need only show that elements of $\mathcal{F}_{B, H}$ give rise to elements of $\mathrm{End}_\mathrm{Hopf}(A, \pi)$ as indicated.
Let $(\textit{f}, \textit{g}) \in \mathcal{F}_{B, H}$ and let $F$ be defined as in part (a). As the reader might suspect, the proof that $F \in \mathrm{End}_\mathrm{Hopf}(A, \pi)$ is somewhat tedious. We will use Lemmas \ref{LemmaFLFR}, \ref{LemmaFLStructure}, and \ref{LemmaFRStructure} without particular reference initially.

It is easy to see that $\pi \circ F = \pi$. Note that $F(1 \times 1) = \textit{f} \, (1)\textit{g} \, (1) \times 1 = 1 \times 1$ and
\begin{eqnarray*}
\epsilon(F(b \times h)) & = & \epsilon(\textit{f} \, (b)\textit{g}(h_{(1)}) \times h_{(2)}) \\
 & = & \epsilon(\textit{f} \, (b)\textit{g}(h_{(1)}))\epsilon(h_{(2)}) \\
 & = & \epsilon(\textit{f} \, (b)\textit{g}(h)) \\
 & = & \epsilon(\textit{f} \, (b))\epsilon(\textit{g}(h)) \\
 & = & \epsilon(b)\epsilon(h)
\end{eqnarray*}
for $b \in B$ and $h \in H$ which means $\epsilon \circ F = \epsilon$. That $F_L = \textit{f}$ and $F_R = \textit{g}$ is easy to see.

Let $b, b' \in B$ and $h, h' \in H$. Then
\begin{eqnarray*}
\lefteqn{F((b \times h)(b' \times h'))} \\
& = & F(b(h_{(1)}{\cdot}b') \times h_{(2)}h') \\
& = & \textit{f} \, (b(h_{(1)}{\cdot}b'))\textit{g}(h_{(2)}h'_{(1)}) \times h_{(3)}h'_{(2)}  \\
& = & \textit{f} \, (b)\textit{f} \, (h_{(1)}{\cdot}b')\textit{g}(h_{(2)}h'_{(1)}) \times h_{(3)}h'_{(2)}  \\
& = & \textit{f} \, (b)\textit{f} \, (h_{(1)}{\cdot}b')\textit{g}(h_{(2)})(h_{(3)}{\cdot}\textit{g}(h'_{(1)})) \times h_{(4)}h'_{(2)} \\
& = & \textit{f} \, (b)\textit{g}(h_{(1)})(h_{(2)}{\cdot}\textit{f} \, (b'))(h_{(3)}{\cdot}\textit{g}(h'_{(1)})) \times h_{(4)}h'_{(2)} \\
& = & \textit{f} \, (b)\textit{g}(h_{(1)})(h_{(2)}{\cdot}(\textit{f} \, (b')\textit{g}(h'_{(1)}))) \times h_{(3)}h'_{(2)} \\
& = & (\textit{f} \, (b)\textit{g}(h_{(1)}) \times h_{(2)})(\textit{f} \, (b')\textit{g}(h'_{(1)}) \times h'_{(2)}) \\
& = & F(b \times h)F(b' \times h').
\end{eqnarray*}
Therefore $F$ is an algebra map. Using (\ref{EqCoSmash}) and (\ref{EqDetltaProd}) we obtain:
\begin{eqnarray*}
\lefteqn{\Delta(F(b \times h))} \\
& = & \Delta(\textit{f} \, (b)\textit{g}(h_{(1)}) \times h_{(2)}) \\
& = & (\textit{f} \, (b)\textit{g}(h_{(1)}))_{(1)} \times (\textit{f} \, (b)\textit{g}(h_{(1)}))_{(2)(-1)}h_{(2)}) \\
& &  \qquad \qquad \qquad \qquad \qquad \qquad  \otimes ((\textit{f} \, (b)\textit{g}(h_{(1)}))_{(2)(0)} \times h_{(3)}) \\
& = & (\textit{f} \, (b)_{(1)}(\textit{f} \, (b)_{(2)(-1)}{\cdot}\textit{g}(h_{(1)})_{(1)}) \times (\textit{f} \, (b)_{(2)(0)}\textit{g}(h_{(1)})_{(2)})_{(-1)}h_{(2)}) \\
& &    \qquad \qquad \qquad \qquad \qquad \qquad   \otimes ((\textit{f} \, (b)_{(2)(0)}\textit{g}(h_{(1)})_{(2)})_{(0)} \times h_{(3)}).
\end{eqnarray*}
Since $\textit{g}$ is a coalgebra map and $B$ is a left $H$-comodule algebra the last expression
\begin{eqnarray*}
\lefteqn{\phantom{a}} \\
& = & (\textit{f} \, (b)_{(1)}(\textit{f} \, (b)_{(2)(-1)}{\cdot}\textit{g}(h_{(1)})) \times (\textit{f} \, (b)_{(2)(0)}\textit{g}(h_{(2)}))_{(-1)}h_{(3)}) \\
& & \qquad \qquad \qquad \qquad \qquad \qquad \otimes ((\textit{f} \, (b)_{(2)(0)}\textit{g}(h_{(2)}))_{(0)} \times h_{(4)}) \\
& = & (\textit{f} \, (b)_{(1)}(\textit{f} \, (b)_{(2)(-1)}{\cdot}\textit{g}(h_{(1)})) \times \textit{f} \, (b)_{(2)(0)(-1)}\textit{g}(h_{(2)})_{(-1)}h_{(3)})\\
 & & \qquad \qquad \qquad \qquad \qquad \qquad  \otimes (\textit{f} \, (b)_{(2)(0)(0)}\textit{g}(h_{(2)})_{(0)} \times h_{(4)}).
\end{eqnarray*}

Using part (d) of Lemma \ref{LemmaFRStructure} and part (c) of Lemma \ref{LemmaFLStructure} the last expression
\begin{eqnarray*}
\lefteqn{\phantom{a}} \\
& = & (\textit{f} \, (b)_{(1)}(\textit{f} \, (b)_{(2)(-1)}{\cdot}\textit{g}(h_{(1)})) \times \textit{f} \, (b)_{(2)(0)(-1)}h_{(2)}S(h_{(4)})h_{(5)}) \\
& & \qquad \qquad \qquad \qquad \qquad \qquad  \otimes (\textit{f} \, (b)_{(2)(0)(0)}\textit{g}(h_{(3)}) \times h_{(6)}) \\
& = & (\textit{f} \, (b)_{(1)}(\textit{f} \, (b)_{(2)(-1)}{\cdot}\textit{g}(h_{(1)})) \times \textit{f} \, (b)_{(2)(0)(-1)}h_{(2)}) \\
& & \qquad \qquad \qquad \qquad \qquad \qquad  \otimes (\textit{f} \, (b)_{(2)(0)(0)}\textit{g}(h_{(3)}) \times h_{(4)}) \\
& = & (\textit{f} \, (b_{(1)})\textit{g}(b_{(2)(-1)})(\textit{f} \, (b_{(2)(0)})_{(-1)}{\cdot}\textit{g}(h_{(1)})) \times \textit{f} \, (b_{(2)(0)})_{(0)(-1)}h_{(2)}) \\
& & \qquad \qquad \qquad \qquad \qquad \qquad  \otimes (\textit{f} \, (b_{(2)(0)})_{(0)(0)}\textit{g}(h_{(3)}) \times h_{(4)}).
\end{eqnarray*}

Using part (d) of Lemma \ref{LemmaFLStructure}, the coassociative comodule axiom, and part (a) of Lemma \ref{LemmaFRStructure} the last expression
\begin{eqnarray*}
\lefteqn{\phantom{a}} \\
& = & (\textit{f} \, (b_{(1)})\textit{g}(b_{(2)(-1)})(b_{(2)(0)(-1)}{\cdot}\textit{g}(h_{(1)})) \times \textit{f} \, (b_{(2)(0)(0)})_{(-1)}h_{(2)}) \\
& & \qquad \qquad \qquad \qquad \qquad \qquad  \otimes (\textit{f} \, (b_{(2)(0)(0)})_{(0)}\textit{g}(h_{(3)}) \times h_{(4)}) \\
& = & (\textit{f} \, (b_{(1)})\textit{g}(b_{(2)(-1)(1)})(b_{(2)(-1)(2)}{\cdot}\textit{g}(h_{(1)})) \times \textit{f} \, (b_{(2)(0)})_{(-1)}h_{(2)}) \\
& & \qquad \qquad \qquad \qquad \qquad \qquad \otimes (\textit{f} \, (b_{(2)(0)})_{(0)}\textit{g}(h_{(3)}) \times h_{(4)}) \\
& = &  (\textit{f} \, (b_{(1)})\textit{g}(b_{(2)(-1)}h_{(1)}) \times \textit{f} \, (b_{(2)(0)})_{(-1)}h_{(2)}) \\
& & \qquad \qquad \qquad \qquad \qquad \qquad \otimes (\textit{f} \, (b_{(2)(0)})_{(0)}\textit{g}(h_{(3)}) \times h_{(4)}).
\end{eqnarray*}

Using part (d) of Lemma \ref{LemmaFLStructure}, the coassociative comodule axiom again, and (\ref{EqCoSmash}) the last expression
\begin{eqnarray*}
& = & (\textit{f} \, (b_{(1)})\textit{g}(b_{(2)(-1)}h_{(1)}) \times b_{(2)(0)(-1)}h_{(2)}) \\
& & \qquad \qquad \qquad \qquad \qquad \qquad  \otimes (\textit{f} \, (b_{(2)(0)(0)})\textit{g}(h_{(3)}) \times h_{(4)}) \\
& = & (\textit{f} \, (b_{(1)})\textit{g}(b_{(2)(-1)(1)}h_{(1)}) \times b_{(2)(-1)(2)}h_{(2)})\\
& & \qquad \qquad \qquad \qquad \qquad \qquad  \otimes (\textit{f} \, (b_{(2)(0)})\textit{g}(h_{(3)}) \times h_{(4)}) \\
& = & F(b_{(1)} \times b_{(2)(-1)}h_{(1)}) \otimes F(b_{(2)(0)} \times h_{(2)}) \\
& = & F((b \times h)_{(1)}) \otimes F((b \times h)_{(2)}).
\end{eqnarray*}
We have shown that $\Delta \circ F = (F \otimes F) \circ \Delta$. Therefore $F$ is a coalgebra map and consequently is a bialgebra map. Since bialgebra maps of Hopf algebras are Hopf algebra maps, the proof of part (a) is complete.

As for part (b), suppose $F \in \mathrm{Aut}_\mathrm{Hopf}(A, \pi)$. Then $F_L$ and $(F^{-1})_L$ are inverses by (\ref{EqnIdASubL}) and part (a) of Corollary \ref{CorFGParts}. Thus $F_L$ is bijective and $(F_L)^{-1} = (F^{-1})_L$.

Conversely, suppose that $F \in  \mathrm{End}_\mathrm{Hopf}(A, \pi)$ and $F_L$ is bijective. Set $G_L = (F_L)^{-1}$. Now $F_R$ has a convolution inverse $J_R$ by Corollary \ref{CorhFRhPrime}. Set $G_R = G_L \circ J_R = (F_L)^{-1} \circ J_R$ and define $G \in \mathrm{End}(A)$ by $G(b \times h) = G_L(b)G_R(h_{(1)}) \times h_{(2)}$ for all $b \in B$ and $h \in H$. Since $G_L$ is an algebra map we compute
\begin{eqnarray*}
G(F(b \times h)) & = & G(F_L(b)F_R(h_{(1)}) \times h_{(2)}) \\
& = & G_L(F_L(b)F_R(h_{(1)}))G_R(h_{(2)}) \otimes h_{(3)} \\
& = & G_L(F_L(b)F_R(h_{(1)}))G_L(J_R(h_{(2)})) \otimes h_{(3)} \\
& = & bG_L(F_R(h_{(1)})J_R(h_{(2)})) \otimes h_{(3)} \\
& = & bG_L((F_R \star J_R)(h_{(1)})) \times h_{(2)} \\
& = & bG_L(\epsilon(h_{(1)})1) \times h_{(2)} \\
& = & b \times h
\end{eqnarray*}
and
\begin{eqnarray*}
F(G(b \times h)) & = & F(G_L(b)G_R(h_{(1)}) \times h_{(2)}) \\
& = & F_L(G_L(b)G_R(h_{(1)}))F_R(h_{(2)}) \otimes h_{(3)}) \\
& = & F_L(G_L(b)(F_L)^{-1}(J_R(h_{(1)})))F_R(h_{(2)}) \otimes h_{(3)} \\
& = & b(J_R(h_{(1)})F_R(h_{(2)})) \times h_{(3)} \\
& = & b((J_R \star F_R)(h_{(1)})) \times h_{(2)} \\
& = & b(\epsilon(h_{(1)})1) \times h_{(2)} \\
& = & b \times h.
\end{eqnarray*}
We have shown $G \circ F = \mathrm{Id}_A = F \circ G$. Therefore $F$ is bijective and the proof of part (b) is complete.
\qed
\medskip

Apropos of the preceding theorem, let $\mathcal{F}^\bullet_{B, H}$ denote the set of $(\textit{f}, \textit{g}) \in \mathcal{F}_{B, H}$ such that $\textit{f}$ is bijective. Then the correspondence of part (a) induces a bijection $\mathcal{F}^\bullet_{B, H} \longrightarrow  \mathrm{Aut}_\mathrm{Hopf}(A, \pi)$.

When $F(1 \times H) \subseteq 1 \times H$ the structure of $F$ is particularly simple.
\begin{Prop}\label{PropRTrivial}
Let $F \in \mathrm{End}_\mathrm{Hopf}(A, \pi)$. Then the following conditions are equivalent:
\begin{enumerate}
\item[{\rm (a)}] $F(1 \times H) \subseteq 1 \times H$.
\item[{\rm (b)}] $F_R(h) \in k1$ for all $h \in H$.
\item[{\rm (c)}] $F_R(h) = \epsilon(h)1$ for all $h \in H$.
\item[{\rm (d)}] $F_R = \eta \circ \epsilon$.
\item[{\rm (e)}] $F(b \times h) = F_L(b) \times h$ for all $b \in B$ and $h \in H$.
\item[{\rm (f)}] $F(1 \times h) = 1 \times h$ for all $h \in H$.
\end{enumerate}
If $F \in \mathrm{Aut}_\mathrm{Hopf}(A, \pi)$ any of these conditions implies $F_L \in \mathrm{Aut}_{\mathcal{YD}-\mathrm{Hopf}}(B)$.
\end{Prop}

\pf
We first show part (a) implies part (b) implies part (c). Suppose $F(1 \times H) \subseteq 1 \times H$. Then $F_R(H) \subseteq k1$ by (\ref{EqFOneTimesH}). Thus $F_R(h) = \epsilon(h)1$ for all $h \in H$ since $\epsilon \circ F_R = \epsilon$. Part (d) is an equivalent expression of part (c).

Suppose part (d) holds. Then $F(b \times h) = F_L(b) \times h$ for all $b \in B$ and $h \in H$ by virtue of part (c) of Lemma \ref{LemmaFLFR}.

We have shown that part (d) implies part (e). Part (e) implies part (f) as $F_L \,(1) = 1$. Part (f) trivially implies part (a). We have shown the conditions are equivalent. If part (c) holds and $F \in \mathrm{Aut}_\mathrm{Hopf}(A, \pi)$ then $F_L \in \mathrm{Aut}_{\mathcal{YD}-\mathrm{Hopf}}(B)$ follows by Lemma \ref{LemmaFLStructure}.
\qed
\medskip

Let
$$
\mathrm{End}'_\mathrm{Hopf}(A, \pi) = \{F \in \mathrm{End}_\mathrm{Hopf}(A, \pi) \, | \, F(1 \times H) \subseteq 1 \times H\}
$$ and
$$
\mathrm{Aut}'_\mathrm{Hopf}(A, \pi) = \mathrm{Aut}_\mathrm{Hopf}(A, \pi) \cap \mathrm{End}'_\mathrm{Hopf}(A, \pi).
$$
Then $\mathrm{End}'_\mathrm{Hopf}(A, \pi)$ is a submonoid of $\mathrm{End}_\mathrm{Hopf}(A, \pi)$ and $\mathrm{Aut}'_\mathrm{Hopf}(A, \pi)$ is a subgroup of $\mathrm{Aut}_\mathrm{Hopf}(A, \pi)$. By virtue of the preceding proposition
\begin{equation}\label{EqEndPrime}
\mathrm{End}'_\mathrm{Hopf}(A, \pi) = \mathrm{End}_\mathrm{Hopf}(A, \pi, j)
\end{equation}
and
\begin{equation}\label{EqAutPrime}
\mathrm{Aut}'_\mathrm{Hopf}(A, \pi) = \mathrm{Aut}_\mathrm{Hopf}(A, \pi, j).
\end{equation}
Most of the proof of the following theorem is established by the preceding proposition. The remainder of the proof is left to the reader.
\begin{Theorem}\label{ThmYDAut}
Let $A = B \times H$ be a biproduct. There is an isomorphism of monoids $\mathrm{End}_{\mathcal{YD}-\mathrm{Hopf}}(B) \simeq \mathrm{End}_\mathrm{Hopf}(A, \pi, j)$ and an isomorphism of groups $\mathrm{Aut}_{\mathcal{YD}-\mathrm{Hopf}}(B) \simeq \mathrm{Aut}_\mathrm{Hopf}(A, \pi, j)$, given by $f \mapsto F$, where $F(b \times h) = f(b) \times h$ for all $b \in B$ and $h \in H$. \qed
\end{Theorem}
\medskip

There are biproducts found in \cite{RadProjNew} where $\mathrm{Aut}_\mathrm{Hopf}(A, \pi) = \mathrm{Aut}'_\mathrm{Hopf}(A, \pi)$, hence $\mathrm{Aut}_\mathrm{Hopf}(A, \pi) = \mathrm{Aut}_\mathrm{Hopf}(A, \pi, j) \simeq \mathrm{Aut}_{\mathcal{YD}-\mathrm{Hopf}}(B)$ by the preceding theorem.
\begin{Theorem}\label{ThmSmallCenterB}
Let $A = B \times H$ be a biproduct. Suppose that the left $H$-module action on $B$ is trivial and that $k1$ is the only subcoalgebra in the center of $B$. Then $\mathrm{Aut}_\mathrm{Hopf}(A, \pi) = \mathrm{Aut}_\mathrm{Hopf}(A, \pi, j) \simeq \mathrm{Aut}_{\mathcal{YD}-\mathrm{Hopf}}(B)$.
\end{Theorem}

\pf
Let $F \in \mathrm{Aut}_\mathrm{Hopf}(A, \pi)$. In light of Theorem \ref{ThmYDAut} we need only show that $F \in \mathrm{Aut}'_\mathrm{Hopf}(A, \pi)$. Now $F_L(b)F_R(h) = F_R(h)F_L(b)$ for all $h \in H$ and $b \in B$ by part (e) of Lemma \ref{LemmaFLStructure} since the $H$-module action on $B$ is trivial. It follows that $F_L$ is bijective by part (b) of Theorem \ref{ThmEndAut} since $F$ is an automorphism of $A$. Now $F_R$ is coalgebra map by part (c) of Lemma \ref{LemmaFRStructure}. Therefore $F_R(H)$ is a subcoalgebra of $B$ in the center of $B$. By assumption $F_R(H) = k1$. Therefore $F \in \mathrm{Aut}'_\mathrm{Hopf}(A, \pi)$ by Proposition \ref{PropRTrivial}.
\qed
\medskip

The biproduct constructed for the proof of \cite[Theorem 6]{RadProjNew} satisfies the hypothesis of Theorem \ref{ThmSmallCenterB}. Here $A = k[\mathcal{G}] \times k[\mathbb{Z}_p]$, where $\mathcal{G}$ is a finite non-abelian simple group and $p$ is a prime integer. $\mathrm{Aut}_\mathrm{Hopf}(A, \pi) \simeq \mathrm{Aut}_\theta(\mathcal{G})$ which consists of all automorphisms of $\mathcal{G}$ which commute with a certain $\theta \in \mathrm{Aut}_\mathrm{Group}(\mathcal{G})$.
\section{$\mathrm{Aut}_\mathrm{Hopf}(A, \pi)$ as a subgroup of a semidirect product.}\label{SecAutsubSemidirect}
Let $B$ be an algebra and $C$ be a coalgebra over $k$.  The group $\mathcal{G}(B) = \mathrm{Aut}_{\mathrm{Group}}(B)$ acts on the convolution algebra $\mathrm{Hom}(C, B)$ by $\textit{f} \vartriangleright \textit{g} = \textit{f} \circ \textit{g}$ for all $\textit{f} \in \mathcal{G}(B)$ and $\textit{g} \in \mathrm{Hom}(C, B)$. This action satisfies
$$
\textit{f} \vartriangleright (\eta \circ \epsilon) = \eta \circ \epsilon \;\; \mbox{and} \;\; \textit{f} \vartriangleright (\textit{g} \star \textit{g}') = (\textit{f} \vartriangleright \textit{g}) \star (\textit{f} \vartriangleright \textit{g}')
$$
for all $\textit{f} \in \mathcal{G}(B)$ and $\textit{g}, \textit{g}' \in \mathrm{Hom}(C, B)$. Let $\mathcal{U}(C, B)$ be the group of units of the monoid $(\mathrm{Hom}(C, B), \star)$. Then $\mathcal{G} \vartriangleright \mathcal{U}(C, B) \subseteq \mathcal{U}(C, B)$; thus there is a group homomorphism,
$$
\varphi : \mathcal{G}(B) \longrightarrow \mathrm{Aut}_{\mathrm{Group}}(\mathcal{U}(C, B))
$$
given by $\varphi(\textit{f})(\textit{g}) = \textit{f} \vartriangleright \textit{g}$ for all $\textit{f} \in \mathcal{G}(B)$ and $\textit{g} \in \mathcal{U}(C, B)$. The resulting group $\mathcal{U}(C, B) \rtimes_\varphi \mathcal{G}(B)$ has product given by
$$
(\textit{g}, \textit{f})(\textit{g}', \textit{f}') = (\textit{g} \star (\textit{f} \circ \textit{g}'), \textit{f} \circ \textit{f'}).
$$

Note the action of $\mathcal{G}(B)$ on $\mathcal{U}(C, B)$ by group homomorphisms is also one on $\mathcal{U}(C, B)^{op}$. For a group $G$ the group $G^{op}$ is the group whose underlying set is $G$ and product is given by $g{\cdot}g' = g'g$ for all $g, g' \in G$. As a result of Theorem \ref{ThmEndAut} and Corollary \ref{CorFGParts}:
\begin{Theorem}\label{ThmAutSubSemi}
Suppose $A = B \times H$ is a biproduct and $\pi : A \longrightarrow H$ is the projection from $A$ onto $H$. Then there is a one-one group homomorphism $\mathrm{Aut}_\mathrm{Hopf}(A, \pi) \longrightarrow \mathcal{U}(H, B)^{op} \rtimes_\varphi \mathcal{G}(B)$ which is given by $F \mapsto (F_R, F_L)$ for all $F \in \mathrm{Aut}_\mathrm{Hopf}(A, \pi)$.
\qed
\end{Theorem}
\section{A normal subgroup of $\mathrm{Aut}_\mathrm{Hopf}(A ,\pi)$}\label{SecNormalSubgroup}
Let $A = B \times H$ be a biproduct. By part (a) of Lemma \ref{LemmaFLStructure}, Theorem \ref{ThmEndAut}, and part (a) of Corollary \ref{CorFGParts}, there is a group homomorphism
$$
\nu : \mathrm{Aut}_\mathrm{Hopf}(A, \pi) \longrightarrow \mathrm{Aut}_\mathrm{Alg}(B)
$$
defined by $\nu (F) = F_L$ for all $F \in \mathrm{Aut}_\mathrm{Hopf}(A, \pi)$. Let $\mathcal{N}(B, H)$ the set of all $F_R$'s  where $F_L = \mathrm{Id}_B$. Then the set of pairs $(\mathrm{Id}_B, \textit{g})$, where $\textit{g} \in \mathcal{N}(B, H)$, is the subset of $\mathcal{F}_{B, H}$ which corresponds to $\mathrm{Ker}(\nu)$ in part (a) of Theorem \ref{ThmEndAut}. Therefore $F \in \mathrm{Ker}(\nu)$ if and only if there is a $\textit{g} \in \mathcal{N}(B, H)$ such that $F(b \times h) = b\textit{g}(h_{(1)}) \times h_{(2)}$ for all $b \in B$ and $h \in H$.
\begin{Prop}\label{PropNBH}
Let $A = B \times H$ be a biproduct. Then:
\begin{enumerate}
\item[{\rm (a)}] $\mathcal{N}(B, H)$ consists of those maps $\textit{g} : H \longrightarrow B$ which satisfy:
\begin{enumerate}
\item[{\rm (1)}] $\textit{g}(b_{(-1)}) \otimes b_{(0)} = 1 \otimes b$ for all $b \in B$;
\item[{\rm (2)}] $(h_{(1)}{\cdot}b)\textit{g}(h_{(2)}) = \textit{g}(h_{(1)})(h_{(2)}{\cdot}b)$ for all $h \in H$ and $b \in B$;
\item[{\rm (3)}] $\textit{g}(hh') = \textit{g}(h_{(1)})(h_{(2)}{\cdot}\textit{g}(h'))$ for all $h, h' \in H$;
\item[{\rm (4)}] $\textit{g}(1) = 1$;
\item[{\rm (5)}] $\textit{g}$ is a coalgebra map; and
\item[{\rm (6)}] $\rho(\textit{g}(h)) = h_{(1)}S(h_{(3)}) \otimes \textit{g}(h_{(2)})$ for all $h \in H$.
\end{enumerate}
\item[{\rm (b)}] $\mathcal{N}(B, H)$ is a group under the convolution product.
\end{enumerate}
\end{Prop}

\pf
Part (a) follows by Theorem \ref{ThmEndAut}, Lemmas \ref{LemmaFLStructure} and \ref{LemmaFRStructure}, and Corollary \ref{CorFLCoalgebraMsap}. Part (b) follows by Corollary \ref{CorFGParts} and (\ref{EqnIdASubL}). It is an interesting exercise to establish part (b) directly from part (a).
\qed
\medskip

Observe that if $\textit{f} \in \mathrm{Aut}_{\mathcal{YD}-\mathrm{Hopf}}(B)$ and $\textit{g} \in \mathcal{N}(B, H)$ then $\textit{f} \circ \textit{g} \in \mathcal{N}(B, H)$ by Proposition \ref{PropNBH}. The reader is left with the exercise of showing that
$$
\varphi : \mathrm{Aut}_{\mathcal{YD}-\mathrm{Hopf}}(B) \longrightarrow \mathrm{Aut}_{\mathrm{Group}}(\mathcal{N}(B, H))
$$
given by $\varphi(\textit{f})(\textit{g}) = \textit{f} \circ \textit{g}$ for $\textit{f} \in \mathrm{Aut}_{\mathcal{YD}-\mathrm{Hopf}}(A, \pi)$ and $\textit{g} \in \mathcal{N}(B, H)$ is a well-defined group homomorphism. See Section \ref{SecAutsubSemidirect}.
\begin{Theorem}\label{ThmAutSemiProd}
Let $A = B \times H$ be a biproduct. Then:
\begin{enumerate}
\item[{\rm (a)}] There is a one-one group homomorphism
$$
\Phi_{\mathcal{YD}} : \mathcal{N}(B, H)^{op} \rtimes_\varphi \mathrm{Aut}_{\mathcal{YD}-\mathrm{Hopf}}(B) \longrightarrow \mathrm{Aut}_\mathrm{Hopf}(A, \pi)
$$
given by $\Phi_{\mathcal{YD}}(\textit{g}, \textit{f})(b \times h) = f(b)g(h_{(1)}) \times h_{(2)}$.
\item[{\rm (b)}] $\mathrm{Im}(\Phi_{\mathcal{YD}}) = \{F \in  \mathrm{Aut}_\mathrm{Hopf}(A, \pi) \, | \, F_L \in  \mathrm{Aut}_{\mathcal{YD}-\mathrm{Hopf}}(B)\}$.
\end{enumerate}
\end{Theorem}

\pf To show part (a) let $(\textit{g},\textit{f} \, ), (\textit{g}',\textit{f}\,' ) \in \mathcal{N}(B, H) \rtimes_\varphi \mathrm{Aut}_{\mathcal{YD}-\mathrm{Hopf}}(B)$. That $\Phi_{\mathcal{YD}}(\textit{g}, \textit{f} \, ) \in \mathrm{Aut}_\mathrm{Hopf}(A, \pi)$ follows by Theorem \ref{ThmEndAut} and Proposition \ref{PropNBH}. That $\Phi_{\mathcal{YD}}((\textit{g}, \textit{f} \, )(\textit{g}\,', \textit{f}\,' )) = \Phi_{\mathcal{YD}}(\textit{g}, \textit{f} \, ) \circ \Phi_{\mathcal{YD}}(\textit{g}\,', \textit{f}\,')$ follows by Corollary \ref{CorFGParts}. Thus $\Phi_{\mathcal{YD}}$ is a group homomorphism. That $\mathrm{Ker}(\Phi_{\mathcal{YD}})$ is trivial is easy to see. Therefore $\Phi_{\mathcal{YD}}$ is one-one.  We have shown part (a).

If $F = \Phi_{\mathcal{YD}}(\textit{g}, \textit{f} \,) \in \mathrm{Im}(\Phi_{\mathcal{YD}})$ then $F_L = \textit{f} \in \mathrm{Aut}_{\mathcal{YD}-\mathrm{Hopf}}(B)$. Now let $F \in \mathrm{Aut}_\mathrm{Hopf}(A, \pi)$ and suppose $F_L \in \mathrm{Aut}_{\mathcal{YD}-\mathrm{Hopf}}(B)$. To establish part (c) we need only show that $F_R$ satisfies the conditions of part (a) of Proposition \ref{PropNBH}. Since $F_L$ is a coalgebra map part (a)(1) holds for $F_R$ by Corollary \ref{CorFLCoalgebraMsap}. Since $F_L$ is a bijective $H$-module map part (a)(2) holds by virtue of part (e) of Lemma \ref{LemmaFLStructure}. Parts (a)(3) - (a)(6) follow by Lemma \ref{LemmaFRStructure}. We have established part (b).
\qed
\begin{Cor}\label{CorAutSemiDirectProduct}
Let $A = B \times H$ be a biproduct. If $\mathrm{Im}(\nu) \subseteq \mathrm{Aut}_{\mathcal{YD}-\mathrm{Hopf}}(B)$ then the map $\Phi_{\mathcal{YD}} : \mathcal{N}(B, H)^{op} \rtimes_\varphi \mathrm{Aut}_{\mathcal{YD}-\mathrm{Hopf}}(B) \longrightarrow \mathrm{Aut}_\mathrm{Hopf}(A, \pi)$ is an isomorphism. \qed
\end{Cor}
\begin{Cor}\label{CorAutSemiDirectProductCC}
Let $A = B \times H$ be a biproduct where $B$ is commutative and $H$ is cocommutative. If $F_L$ is a coalgebra map for all $F \in \mathrm{Aut}_\mathrm{Hopf}(A, \pi)$ then the map $\Phi_{\mathcal{YD}} : \mathcal{N}(B, H)^{op} \rtimes_\varphi \mathrm{Aut}_{\mathcal{YD}-\mathrm{Hopf}}(B) \longrightarrow \mathrm{Aut}_\mathrm{Hopf}(A, \pi)$ is an isomorphism.
\end{Cor}

\pf
Let $F \in \mathrm{Aut}_\mathrm{Hopf}(A, \pi)$. Then $F_L$ is a map of left $H$-modules by part (b) of Corollary \ref{CorFLModuleMap}. If $F_L$ is a coalgebra map then $F_L \in \mathrm{Aut}_{\mathcal{YD}-\mathrm{Hopf}}(B)$ by Lemma \ref{LemmaFLStructure} and Theorem \ref{ThmEndAut}. Now the corollary follows by Corollary \ref{CorAutSemiDirectProduct}.
\qed
\medskip

In light of Corollary \ref{CorAutSemiDirectProductCC}, we consider biproducts when $H$ is cocommutative. Suppose $H$ is cocommutative. This is the case, for example, when $H$ is a group algebra. We first revisit Proposition \ref{PropNBH}. Let $h \in H$ and $\textit{g} \in \mathcal{N}(B, H)$. Then (a)(6) is $\rho(\textit{g}(h)) = 1 \otimes \textit{g}(h)$. Let $b \in B$. In any event
$$
\textit{g}(h)b = \textit{g}(h_{(1)})(h_{(2)}{\cdot}(S(h_{(3)}){\cdot}b)) = ((h_{(1)}S(h_{(3)})){\cdot}b)\textit{g}(h_{(2)}).
$$
Thus $\textit{g}(h)b = b\textit{g}(h)$ since $H$ is cocommutative; hence (a)(2) is $\mathrm{Im}(\textit{g}) \subseteq \mathrm{Z}(B)$. Therefore $\mathcal{N}(B, H) \subseteq \mathrm{Z}(\mathrm{Hom}(H, B))$, where $\mathrm{Hom}(H, B)$ has the convolution algebra structure. In particular $\mathcal{N}(B, H)$ is an abelian group.

Let $F \in \mathrm{Aut}_\mathrm{Hopf}(A, \pi)$ and $G \in \mathrm{Ker}(\nu)$. Then $G_R \in \mathrm{Z}(\mathrm{Hom}(H, B))$. We show that $F_L \circ G_R \in \mathcal{N}(B, H)$. Now $F \circ G \circ F^{-1} \in \mathrm{Ker}(\nu)$. Observe that $(F^{-1})_L \circ F_R$ is a left inverse of $(F^{-1})_R$ by part (b) of Corollary \ref{CorFGParts} and is therefore a right inverse as well by Corollary \ref{CorhFRhPrime}. Using Corollary \ref{CorFGParts}, (\ref{EqnIdASubL}), and the fact that $F_L$ is an algebra map, we compute
\begin{eqnarray*}
(F \circ G \circ F^{-1})_R &  = & ((F \circ G)_L \circ (F^{-1})_R) \star (F \circ G)_R \\
& = & ((F_L \circ G_L) \circ (F^{-1})_R) \star ((F_L \circ G_R) \star F_R) \\
& = & (F_L \circ (F^{-1})_R) \star ((F_L \circ G_R) \star F_R) \\
& = & F_L \circ ((F^{-1})_R \star (G_R \star ((F_L)^{-1} \circ F_R)))  \\
& = & F_L \circ (G_R \star ((F^{-1})_R \star ((F_L)^{-1} \circ F_R))) \\
& = & F_L \circ (G_R \star (\eta \circ \epsilon)) \\
& = & F_L \circ G_R.
\end{eqnarray*}
Therefore $F_L \circ G_R \in \mathcal{N}(B, H)$.  The reader is left with the exercise of showing that
$$
\varphi : \mathrm{Im}(\nu) \longrightarrow \mathrm{Aut}_{\mathrm{Group}}(\mathcal{N}(B, H))
$$
given by $\varphi(\textit{f})(\textit{g}) = \textit{f} \circ \textit{g}$ for $\textit{f} \in \mathrm{Im}(\nu)$ and $\textit{g} \in \mathcal{N}(B, H)$ is a well-defined group homomorphism. See Section \ref{SecAutsubSemidirect}.

\begin{Theorem}\label{ThmAutSemiProdHCoComm}
Let $A = B \times H$ be a biproduct where $H$ is cocommutative. Then:
\begin{enumerate}
\item[{\rm (a)}] There is a one-one group homomorphism
$$
\Phi_\nu : \mathcal{N}(B, H)^{op} \rtimes_\varphi \mathrm{Im}(\nu) \longrightarrow \mathrm{Sym}(A)
$$
given by $\Phi_\nu(\textit{g}, \textit{f}\,)(b \times h) = \textit{f}\,(b)\textit{g}\,(h_{(1)}) \times h_{(2)}$.
\item[{\rm (b)}] Let $(\textit{g}, \textit{f}\,) \in \mathcal{N}(B, H)^{op} \rtimes_\varphi \mathrm{Im}(\nu)$. Then $\Phi_\nu(\textit{g}, \textit{f}\,) \in \mathrm{Aut}_\mathrm{Hopf}(A, \pi)$ if and only if $\textit{f}\, \in \mathrm{Aut}_{\mathcal{YD}-\mathrm{Hopf}}(B)$.
\end{enumerate}
\end{Theorem}

\pf We first note that $\Phi_\nu(\textit{g}, \textit{f}\,)$ a linear endomorphism of $A$. Observe that $\Phi_\nu(\eta {\circ} \epsilon, \mathrm{Id}_B) = \mathrm{Id}_A$ and $\Phi_\nu((\textit{g}, \textit{f}\,)(\textit{g}\,', \textit{f}\,')) = \Phi_\nu(\textit{g}, \textit{f}\,) {\circ}\Phi_\nu(\textit{g}\,', \textit{f}\,')$. Thus $\Phi_\nu(\textit{g}, \textit{f}\,) \in \mathrm{Sym}(A)$ and $\Phi_\nu$ is a homomorphism.

Suppose $\Phi_\nu(\textit{g}, \textit{f}\,) = \mathrm{Id}_A$. Then $\textit{f}\,(b)\textit{g}\,(h_{(1)}) \times h_{(2)} = b \times h$ for all $b \in B$ and $h \in H$ from which $(\textit{g}, \textit{f}\,) = (\eta {\circ} \epsilon, \mathrm{Id}_B)$ follows. Hence $\Phi_\nu$ is one-one. We have shown part (a).

Suppose $\Phi_\nu(\textit{g}, \textit{f}\,) \in \mathrm{Aut}_\mathrm{Hopf}(A, \pi)$. In light of part (c) of Theorem \ref{ThmAutSemiProd}, to establish part (b) we need only show that $\textit{f}\, \in \mathrm{Aut}_{\mathcal{YD}-\mathrm{Hopf}}(B)$. In light of Lemma \ref{LemmaFLStructure} and Proposition \ref{PropNBH} we need only show that $\textit{f}\,$ is a module map. Now $\textit{f}\,(h_{(1)}{\cdot}b)\textit{g}\,(h_{(2)}) = \textit{g}\,(h_{(1)})(h_{(2)}{\cdot}\textit{f}\,(b))$ for all $h \in H$ and $b \in B$ by part (e) of Lemma \ref{LemmaFLStructure}. Fix $b$ and define $L, R \in \mathrm{Hom}(H, B)$ by $L(h) = \textit{f}\,(h{\cdot}b)$ and $R(h) = h{\cdot}\textit{f}\,(b)$ for all $h \in H$. Then $L \star \textit{g}\, = \textit{g}\, \star R = R \star \textit{g}\,$; the second equation holds since $\textit{g}\, \in \mathrm{Z}(\mathrm{Hom}(H, B))$. Now $\textit{g}\,$ has a convolution inverse by Corollary \ref{CorhFRhPrime}. Therefore $L = R$ which establishes $\textit{f}\,$ is a module map.
\qed
\section{A related biproduct and automorphism group}\label{SecRelated}
Throughout this section $A = B \times H$ is a biproduct. Here we approach the problem of computing $\mathrm{Im}(\nu)$ of Section \ref{SecNormalSubgroup} by considering a related biproduct.

Let $H'$ be a Hopf subalgebra of $H$ which satisfies $\rho(B) \subseteq H' \otimes B$. Then $B$ is also a left $H'$-comodule. Regard $B$ as a left $H'$-module by $H$-module action restricted to $H'$. Then the biproduct $A' = B \times H'$ is defined and $A'$, regarded as a subspace of $A$, is a Hopf subalgebra of $A$. Let $\pi' : A' \longrightarrow H'$ be the projection onto $H'$. Then $\pi' = \pi|_{B \times H'}$.

Now let $F \in \mathrm{Aut}_\mathrm{Hopf}(A, \pi)$. Then $F(B \times C) \subseteq B \times C$ for all subcoalgebras $C$ of $H$ by part (c) of Lemma \ref{LemmaFLFR}; in particular $F(B \times H') \subseteq B \times H'$. Since this inclusion holds for $F^{-1}$ it follows that $F(B \times H') = B \times H'$. Observe that $\pi' \circ (F|_{B \times H'}) = \pi'$ since $\pi \circ F = \pi$. Thus by (\ref{EqnIdASubL}) and part (a) of Corollary \ref{CorFGParts} we have a group homomorphism
$$
\Theta : \mathrm{Aut}_\mathrm{Hopf}(A, \pi) \longrightarrow \mathrm{Aut}_\mathrm{Hopf}(A', \pi')
$$
defined by $\Theta (F) = F|_{B \times H'}$ for all $F \in \mathrm{Aut}_\mathrm{Hopf}(A, \pi)$. By part (a) of Lemma \ref{LemmaFLFR} and (\ref{EqFOneTimesH}) we see
\begin{equation}\label{EqThetaFLFR}
\Theta(F)_L = F_L \;\; \mbox{and} \;\; \Theta(F)_R = F_R|_{H'}
\end{equation}
for all $F \in \mathrm{Aut}_\mathrm{Hopf}(A, \pi)$. Observe that $\mathrm{Ker}(\Theta)$ corresponds to the set of pairs $(\mathrm{Id}_B, \textit{g}) \in \mathcal{F}_{B, H}$ such that $\textit{g}(h) = \epsilon(h)1$ for all $h \in H'$.

Let $\nu' : \mathrm{Aut}_\mathrm{Hopf}(A', \pi') \longrightarrow \mathrm{Aut}_\mathrm{Alg}(A')$ be the group homomorphism defined by $\nu'(F) = F_L$ for all $F \in \mathrm{Aut}_\mathrm{Hopf}(A', \pi')$. Then $\nu' \circ \Theta = \nu$, where $\nu$ is the map mentioned above. Thus $\mathrm{Im}(\nu) \leq \mathrm{Im}(\nu')$. Which elements of $\mathrm{Im}(\nu')$ belong to $\mathrm{Im}(\nu)$ is easily explained in the context of Theorem \ref{ThmEndAut}. Let $\textit{f} \in \mathrm{Im}(\nu')$. Then $\textit{f} \in \mathrm{Im}(\nu)$ if and only if there are pairs $(\textit{f}, \textit{g}') \in \mathcal{F}_{B, H'}$ and $(\textit{f}, \textit{g}) \in \mathcal{F}_{B, H}$ such that $\textit{g}' = \textit{g}|_{H'}$.

In biproducts of interest to us below $H'$ is the group algebra of a cyclic group and the left $H'$-module action on $B$ is trivial.
\section{A generic example from \cite{RadProjNew}}\label{SecGeneric}
Many examples of \cite{RadProjNew} are described as special cases of \cite[Theorem 2]{RadProjNew}. Our Theorem \ref{ThmFromProjNew} is \cite[Theorem 2]{RadProjNew} less one technical detail. First a bit of background.

Let $V$ be a vector space and $T$ be an endomorphism of $V$. For $\lambda \in k$ let $V_\lambda = \{ v \in V \, | \, T(v) = \lambda v\}$. Then $T$ is \textit{diagonalizable} if $V = \sum_{\lambda \in k} V_\lambda$. In any case the sum is direct.

Suppose $T$ is an automorphism and $V$ has a basis $B$ such that $T(B) = B$ and $(T)$-orbits of $B$ are finite. Then $T$ is diagonalizable if and only if $k^\times$ contains a primitive $r^{th}$ root of unity for every orbit length $r$. Suppose that $T$ is diagonalizable and has finite order $N$. Then $N$ is the least common multiple of orbit lengths. Furthermore $k^\times$ has a primitive $N^{th}$ root of unity $\lambda$ and the cyclic group $U = (\lambda)$ is the subgroup of $k^\times$ generated by the set of eigenvalues of $T$.

If $\mathcal{G}$ is a group and $\theta \in \mathrm{Aut}_{\mathrm{Group}}(\mathcal{G})$ we let $\mbox{\boldmath $\theta$} : k[\mathcal{G}] \mapsto k[\mathcal{G}]$ denote the linear extension of $\theta$ to $k[\mathcal{G}]$.
\begin{Theorem}\label{ThmFromProjNew}
Let $G$ and $\mathcal{G}$ be groups, let $\pi : G \longrightarrow \mathrm{Aut}_{\mathrm{Group}}(\mathcal{G})$ be a group homomorphism, and let $\theta \in \mathrm{Aut}_{\mathrm{Group}}(\mathcal{G})$. Suppose that the $(\theta)$-orbits of $\mathcal{G}$ are finite and that $\mbox{\boldmath $\theta$}$ is diagonalizable. Assume further that:
\begin{enumerate}
\item[{\rm (a)}] $\pi(g) \circ \theta = \theta \circ \pi(g)$ for all $g \in G$.
\item[{\rm (b)}] There exists an isomorphism $U \longrightarrow \mathbf{U}$ $(\lambda \mapsto \mbox{\boldmath $\lambda$})$ of the subgroup $U \subseteq k^\times$ generated by the eigenvalues of $\mbox{\boldmath $\theta$}$ and a subgroup $\mathbf{U}$ of $G$ which satisfies $\mathbf{U} \subseteq \mathrm{Ker}(\pi)\cap \mathrm{Z}(G)$.
\end{enumerate}
Let $H = k[G]$ and $B = k[\mathcal{G}]$ Then:
\begin{enumerate}
\item[{\rm (c)}] $B$ is a left $H$-module determined by $g{\cdot}b = \pi(g)(b)$ for all $g \in G$ and $b \in \mathcal{G}$.
\item[{\rm (d)}] $B$ is a left $H$-comodule where $\rho(b) = \mbox{\boldmath $\lambda$} \otimes b$ for $b \in B_\lambda$.
\item[{\rm (e)}] The biproduct $A = B \times H$ is defined for the Hopf algebra $B$ with these module and comodule structures. Furthermore $A$ is a cosemisimple involutory Hopf algebra over $k$.
\end{enumerate}
\qed
\end{Theorem}
\medskip

About the preceding theorem. Observe that $\rho(B) \subseteq k[\mathbf{U}]$. Our discussion of Section \ref{SecRelated} applies to $H' = k[\mathbf{U}]$. Note that $\mathbf{U}$ is a cyclic group and that the left $H'$-module action on $B$ is trivial. A group of great interest to us is $\mathrm{Aut}_\mathrm{Hopf}(A', \pi')$.

To continue our discussion of Theorem \ref{ThmFromProjNew}. We emphasize that $B = k[\mathcal{G}]$ is a Hopf algebra over $k$ whose coproduct also satisfies (\ref{EqDetltaProd}). Note that a map $\textit{f} : B \longrightarrow B$ is a map of left $H$-comodules if and only if $\textit{f} \, (B_\lambda) \subseteq B_\lambda$ for all $\lambda \in k$; that is if and only if $\mbox{\boldmath $\theta$} \circ \textit{f} = \textit{f} \circ \mbox{\boldmath $\theta$}$.

In many cases $\mathrm{Im}(\pi) \subseteq (\theta)$, for example when the left $H$-module action on $B$ is trivial. Suppose $\mathrm{Im}(\pi) \subseteq (\theta)$. If $\textit{f}$ is a map of left $H$-comodules then $\textit{f}$ is a map of left $H$-modules as well. We revisit Theorem \ref{ThmEndAut} in the next corollary. Its proof is left to the reader:
\begin{Cor}\label{CorFBHfromProjNew}
Let $A = B \times H$ be the biproduct described in Theorem \ref{ThmFromProjNew} and further assume that $\mathrm{Im}(\pi) \subseteq (\theta)$. Then $(\textit{f}, \textit{g}) \in \mathcal{F}^\bullet_{B, H}$ if and only if:
\medskip

\noindent
$\textit{f} : B \longrightarrow B$ is an algebra automorphism such that
\begin{enumerate}
\item[{\rm (a)}] $\epsilon \circ \textit{f} = \epsilon$,
\item[{\rm (b)}] $\Delta \textit{f} \, (b) = \textit{f} \, (b_{(1)})\textit{g}(b_{(2)(-1)}) \otimes \textit{f} \, (b_{(2)(0)})$ for all $b \in B$,
\item[{\rm (c)}] $\mbox{\boldmath $\theta$} \circ \textit{f} = \textit{f} \circ \mbox{\boldmath $\theta$}$; and
\end{enumerate}
$\textit{g} : H \longrightarrow B$ is a Hopf algebra map such that
\begin{enumerate}
\item[{\rm (d)}] $\mathrm{Im}(\textit{g}) \subseteq \mathrm{Z}(B)$ and
\item[{\rm (e)}] $\mbox{\boldmath $\theta$} \circ \textit{g} = \textit{g}$.
\end{enumerate}
\qed
\end{Cor}

When $\mathcal{G}$ is abelian (d) is redundant. We consider the case when $\mathcal{G}$ is abelian in the next section. In this situation:
\begin{Cor}\label{CorMCGAbelian}
Let $A = B \times H$ be the biproduct described in Theorem \ref{ThmFromProjNew}. Suppose that $\mathrm{Im}(\pi) \subseteq (\theta)$ and $\mathcal{G}$ is abelian. Then $\mathcal{N}(B, H)$ consists of the Hopf algebra maps $\textit{g} : H \longrightarrow B$ which satisfy $\mbox{\boldmath $\theta$} \circ \textit{g} = \textit{g}$ and $\textit{g}(\mathbf{U}) = \{ 1 \}$.
\qed
\end{Cor}
\section{The special case $A' = k[\mathcal{G}] \times k[\mathbf{U}]$ when $\mathcal{G}$ is finite abelian}\label{SecGammaGAbelian}
Let $\mathcal{G}$ be a finite abelian group of order $n$, suppose $k$ has a primitive $n^{th}$ root of unity $\omega$, and let $A = k[\mathcal{G}] \times H$ be the biproduct of Theorem \ref{ThmFromProjNew}. Then $B = k[\mathcal{G}]$ is commutative and $H$ is cocommutative. The map $\Phi_\mathcal{YD}$ of Corollary \ref{CorAutSemiDirectProductCC} is a group isomorphism if and only $F_L$ is a coalgebra map for all $F \in \mathrm{Aut}_\mathrm{Hopf}(A, \pi)$. We are interested in when $\Phi_\mathcal{YD}$ is an isomorphism. In light of (\ref{EqThetaFLFR}) we are led to consider the special case $A' = k[\mathcal{G}] \times k[\mathbf{U}]$ of the biproduct described in the preceding section. Recall that the $k[\mathbf{U}]$-module action on $B$ is trivial. Corollary \ref{CorFBHfromProjNew} applies to $A'$.

In this section we study $A'$. We describe the group $\mathrm{Aut}_{\mathrm{Hopf}}(A', \pi')$ in terms permutations and characters in Theorem \ref{ThmTwoFFs}. Our discussion is based on the orthogonal basis of idempotents for $B = k[\mathcal{G}]$.

Write $\mathcal{G} = \mathcal{G}_1 \times \cdots \times \mathcal{G}_t$ as a product of cyclic groups, where $\mathcal{G}_j = (g_j)$ and has order $n_j$ for all $1 \leq j \leq t$. Then $n = n_1 \cdots n_t$ and $\omega_j = \omega^{n/n_j}$ is a primitive $n_j^{th}$ root of unity each $1 \leq j \leq t$. Let $\mathbf{G} = \mathbb{Z}_{n_1} \oplus \cdots \oplus \mathbb{Z}_{n_t}$ with its usual \textit{ring} structure. For $\mathbf{m} = (m_1, \ldots, m_t) \in \mathbf{G}$ we define formal powers
$$
g^{(\mathbf{m})} = (g_1^{m_1}, \ldots, g_t^{m_t}) \;\; \mbox{and} \;\; \omega^{(\mathbf{m})} = \omega_1^{m_1} \cdots \omega_t^{m_t}.
$$
The maps $\mathbf{G} \longrightarrow \mathcal{G}$ and $\mathbf{G} \longrightarrow k^\times$ defined by $\mathbf{m} \mapsto g^{(\mathbf{m})}$ and $\mathbf{m} \mapsto \omega^{(\mathbf{m})}$ respectively are group homomorphisms, the first is an isomorphism. Set
\begin{equation}\label{EqEMDef}
e_{\mathbf{m}}= \sum_{\mathbf{r} \in {\mathbf{G}}} (\omega^{(\mathbf{m}\mathbf{r})}/n)g^{(\mathbf{r})}
\end{equation}
for all $\mathbf{m} \in \mathrm{G}$. Then $\{e_{\mathbf{m}}\}_{\mathbf{m} \in \mathbf{G}}$ is a linear basis for $B$ which satisfies
\begin{equation}\label{EqOrthBasis}
\sum_{\mathbf{m} \in \mathbf{G}} e_{\mathbf{m}} = 1 \;\; \mbox{and} \;\; e_{\mathbf{m}} e_{\mathbf{n}} = \delta_{\mathbf{m}, \mathbf{n}}e_{\mathbf{n}}
\end{equation}
for all $\mathbf{m}, \mathbf{n} \in \mathbf{G}$. Since $\{e_{\mathbf{m}}\}_{\mathbf{m} \in \mathbf{G}}$ is a linear basis for $B = k[\mathcal{G}]$ the algebra structure of $B$ is determined by (\ref{EqOrthBasis}). Observe that
\begin{equation}\label{EqGTimesEm}
g^{(\mathbf{m})}e_{\mathbf{n}} = \omega^{(-\mathbf{m}\mathbf{n})}e_{\mathbf{n}}
\end{equation}
for all $\mathbf{m}, \mathbf{n} \in \mathbf{G}$. Let $\mathbf{m} \in \mathbf{G}$. The map $\alpha_\mathbf{m} : \mathbf{G} \longrightarrow k^\times$ defined by $\alpha_\mathbf{m}(\mathbf{n}) = \omega^{(-\mathbf{m}\mathbf{n})}$ for all $\mathbf{n} \in \mathbf{G}$ belongs to the character group $\widehat{\mathbf{G}} = \mathrm{Char}(\mathbf{G}, k^\times)$. Note
\begin{equation}\label{EqngmSumChar}
g^{(\mathbf{m})} = \sum_{\mathbf{n} \in \mathbf{G}} \alpha_\mathbf{m}(\mathbf{n})e_\mathbf{n}
\end{equation}
by (\ref{EqOrthBasis}) and (\ref{EqGTimesEm}). Therefore the group homomorphism $\mathbf{G} \longrightarrow \widehat{\mathbf{G}}$ defined by $\mathbf{m} \mapsto \alpha_\mathbf{m}$ is one-one and consequently is an isomorphism.

The coalgebra structure of $B$ in terms of this basis is given by
\begin{equation}\label{EQnOrthCoalgebra}
\Delta(e_\mathbf{m}) = \sum_{\mathbf{j} \in \mathbf{G}} e_\mathbf{j} \otimes e_{\mathbf{m} - \mathbf{j}} \;\; \mbox{and} \;\; \epsilon(e_\mathbf{m}) = \delta_{\mathbf{m}, \mathbf{0}}
\end{equation}
for all $\mathbf{m} \in \mathbf{G}$, where $\mathbf{0}$ is the neutral element for the additive group structure of $\mathbf{G}$. Note that the unique integral $\Lambda$ for the Hopf algebra $B$ which satisfies $\epsilon(\Lambda) = 1$ is $\Lambda = e_\mathbf{0}$.

For $\tau \in \mathrm{Sym}(\mathbf{G})$ let $F_\tau$ be the linear automorphism of $B$ defined by $F_\tau(e_\mathbf{m}) = e_{\tau(\mathbf{m})}$ for all $\mathbf{m} \in \mathbf{G}$.
\begin{Lemma}\label{LemmaSymBfG}
The correspondence $\tau \mapsto F_\tau$ describes a group isomorphism $\mathrm{Sym}(\mathbf{G}) \simeq \mathrm{Aut}_\mathrm{Alg}(B)$ which restricts to an isomorphism $\mathrm{Aut}_\mathrm{Group}(\mathbf{G}) \simeq \mathrm{Aut}_\mathrm{Hopf}(B)$.
\end{Lemma}

\pf
Let $\tau, \sigma \in \mathrm{Sym}(\mathbf{G})$. Then $F_\sigma \circ F_\tau = F_{\sigma \circ \tau}$. It is easy to see that $F_\tau \in \mathrm{Aut}_\mathrm{Alg}(B)$ and that the map $\tau \mapsto F_\tau$ is one-one. Let $f \in \mathrm{Aut}_\mathrm{Alg}(B)$. We will show that $f = F_\tau$ for some $\tau \in \mathrm{Sym}(\mathbf{G})$.

The $ke_{\mathbf{m}}$'s are the minimal ideals of $k[\mathcal{G}]$ and $f$ is an automorphism. Thus $f$ permutes the $ke_\mathbf{m}$'s. Therefore there is a $\tau \in \mathrm{Sym}(\mathbf{G})$ such that $f(ke_\mathbf{m}) = ke_{\tau(\mathbf{m})}$ for all $\mathbf{m} \in \mathbf{G}$. Let $\mathbf{m} \in \mathbf{G}$. Since $f(e_\mathbf{m})$ is a non-zero idempotent $f(e_\mathbf{m}) = e_\mathbf{\tau(\mathbf{m})}$. We have shown $f = F_\tau$ which establishes the first isomorphism.

As for the second, let $\tau \in \mathrm{Sym}(\mathbf{G})$. It is easy to see that $\epsilon \circ F_\tau = \epsilon$ if and only if $\tau(\mathbf{0}) = \mathbf{0}$. Let $\mathbf{m} \in \mathbf{G}$. Then
$$
\Delta(F_\tau(e_\mathbf{m})) = \Delta (e_{\tau(\mathbf{m})}) = \sum_{\mathbf{r} \in \mathbf{G}} e_\mathbf{r} \otimes e_{\tau(\mathbf{m}) - \mathbf{r}} = \sum_{\mathbf{r} \in \mathbf{G}} e_{\tau(\mathbf{r})} \otimes e_{\tau(\mathbf{m}) - \tau(\mathbf{r})}
$$
and
$$
(F_\tau \otimes F_\tau)(\Delta (e_\mathbf{m})) = \sum_{\mathbf{r} \in \mathbf{G}} F_\tau(e_\mathbf{r}) \otimes F_\tau(e_{\mathbf{m} - \mathbf{r}}) = \sum_{\mathbf{r} \in \mathbf{G}} e_{\tau(\mathbf{r})} \otimes e_{\tau(e_{\mathbf{m} - \mathbf{r}})}.
$$
Therefore $\Delta \circ F_\tau = (F_\tau \otimes F_\tau) \circ \Delta$ if and only if $\tau(\mathbf{m}) - \tau(\mathbf{r}) = \tau(\mathbf{m} - \mathbf{r})$ for all $\mathbf{m}, \mathbf{r} \in \mathbf{G}$. We have shown that $F_\tau$ is a coalgebra map if and only if $\tau(\mathbf{0}) = \mathbf{0}$ and this last condition holds which is equivalent to saying $\tau \in \mathrm{Aut}_\mathrm{Group}(\mathbf{G})$. We have established the second isomorphism.
\qed
\medskip

We turn our attention to $A' = k[\mathcal{G}] \times k[\mathbf{U}]$. Let $F \in \mathrm{Aut}_\mathrm{Hopf}(A', \pi')$ and $(\textit{f}, \textit{g}) \in \mathcal{F}^\bullet_{k[\mathcal{G}], k[\mathbf{U}]}$ correspond to it according to Theorem \ref{ThmEndAut}. First of all
\begin{equation}\label{EqBoldThetaTau}
\mbox{\boldmath $\theta$} = F_\sigma \;\; \mbox{and} \;\; \textit{f} = F_\tau
\end{equation}
for unique $\sigma \in \mathrm{Aut}_\mathrm{Group}(\mathbf{G})$ and $\tau \in \mathrm{Sym}(\mathbf{G})$ by Lemma \ref{LemmaSymBfG}. Since $\epsilon{\circ}\textit{f} = \epsilon$ by part (a) of Corollary \ref{CorFBHfromProjNew}, by the second equation of (\ref{EQnOrthCoalgebra}) it follows that
\begin{equation}\label{EqTauZeroIsZero}
\tau(\mathbf{0}) = \mathbf{0}.
\end{equation}
By part (c) of Corollary \ref{CorFBHfromProjNew} we have $F_\sigma \circ F_\tau = F_\tau \circ F_\sigma$ and therefore
\begin{equation}\label{EqTauThetaCommute}
\sigma \circ \tau = \tau \circ \sigma.
\end{equation}
Further analysis of $\textit{f} = F_\tau$ is done in terms of the $(\mbox{\boldmath $\theta$})$-orbits of $k[\mathcal{G}]$. We continue with the notation of Theorem \ref{ThmFromProjNew} and draw from the discussion preceding it in Section \ref{SecGeneric}.

Let $N$ be the order of $\mbox{\boldmath $\theta$}$, which is also the order of $\theta$. Then $N$ is the least common multiple of the $(\mbox{\boldmath $\theta$})$-orbit lengths of $k[\mathcal{G}]$, which is also the least common multiple of the lengths of the $(\theta)$-orbit of $\mathcal{G}$. The eigenvalues of $\mbox{\boldmath $\theta$}$ generate a (cyclic) subgroup $U$ of $k^\times$ of order $N$. Choose a generator $\lambda$ for $U$ and let $\mbox{\boldmath $\lambda$} \in \mathbf{U}$ be its counterpart. Then $\lambda$ is a primitive $N^{th}$ root of unity, $\mbox{\boldmath $\lambda$}$ has order $N$, and $\mathbf{U} = (\mbox{\boldmath $\lambda$})$.

Let $b \in B$. Denote the $(\mbox{\boldmath $\theta$})$-orbit of $b$ by $\mathcal{O}_{\theta, b}$ and set $r = |\mathcal{O}_{\theta, b}|$. Then $\mathcal{O}_{\theta, b} = \{b, \mbox{\boldmath $\theta$}(b), \ldots, \mbox{\boldmath $\theta$}^{r-1}(b)\}$, where $\mbox{\boldmath $\theta$}^r(b) = b$, and $r$ divides $|\mbox{\boldmath $\theta$}| = |\mathbf{U}| = N$. In particular $\lambda_r = \lambda^{N/r}$ is a primitive $r^{th}$ root of unity.

For $0 \leq i \leq r - 1$ define $b_{\lambda_r^i} = \sum_{\ell = 0}^{r-1}(\lambda_r^{-i\ell}/r)\mbox{\boldmath $\theta$}^\ell(b)$. It is an easy exercise to see that $\mbox{\boldmath $\theta$}(b_{\lambda_r^i}) = \lambda_r^ib_{\lambda_r^i}$ and that $\mbox{\boldmath $\theta$}^i(b) = \sum_{\ell = 0}^{r-1}\lambda_r^{i\ell}b_{\lambda_r^\ell}$ for all $0 \leq i \leq r-1$. In particular the linear span $\mathrm{sp}(\mathcal{O}_{\theta, b})$ of $\mathcal{O}_{\theta, b}$ has a basis of eigenvectors belonging to distinct eigenvalues of $\theta$ and these constitute the $r^{th}$ roots of unity of $k^\times$.

Let $\mathbf{m} \in \mathbf{G}$. We specialize our discussion and notation for $b = e_\mathbf{m}$. We write $\mathcal{O}_{\theta, \mathbf{m}}$ for $\mathcal{O}_{\theta, e_\mathbf{m}}$, $|\mathbf{m}|$ for $r$, $\lambda_\mathbf{m}$ for $\lambda_r$, and $e_{\mathbf{m}, \lambda_\mathbf{m}^i}$ for $b_{\lambda_r^i}$. Let $\mbox{\boldmath $\lambda$}_\mathbf{m} \in \mathbf{U}$ be the counterpart of $\lambda_\mathbf{m} \in k^\times$. Then $\mbox{\boldmath $\lambda$}_\mathbf{m} = \mbox{ \boldmath $\lambda$}^{N/|\mathbf{m}|}$,
\begin{equation}\label{EqnmOrbit}
\mathcal{O}_{\theta ,\mathbf{m}} = \{e_\mathbf{m}, e_{\sigma(\mathbf{m})}, \ldots, e_{\sigma^{|\mathbf{m}|-1}(\mathbf{m})}\}, \;\; \mbox{where $\sigma^{|\mathbf{m}|}(\mathbf{m}) = \mathbf{m}$},
\end{equation}
\begin{equation}\label{EqnEmLambdaM}
e_{\mathbf{m}, \lambda_\mathbf{m}^i} = \sum_{\ell \in \mathbb{Z}_{|\mathbf{m}|}} (\lambda_\mathbf{m}^{-i\ell}/|\mathbf{m}|)e_{\sigma^\ell(\mathbf{m})},
\end{equation}
and
\begin{equation}\label{EqnThetaEmLambdaM}
e_{\sigma^i(\mathbf{m})} = \sum_{\ell \in \mathbb{Z}_{|\mathbf{m}|}}\lambda_\mathbf{m}^{i\ell}e_{\mathbf{m}, \lambda^\ell_\mathbf{m}}
\end{equation}
for all $0 \leq i \leq |\mathbf{m}|-1$. Now $\mathcal{O}_{\theta, b} = \mathcal{O}_{\theta, \theta(b)}$; thus $|\mathbf{m}| = |\sigma(\mathbf{m})|$. Consequently
\begin{equation}\label{EqLambdaSigmaBF}
\lambda_\mathbf{m} = \lambda_{\sigma(\mathbf{m})} \;\; \mbox{and} \;\; \mbox{ \boldmath $\lambda$}_\mathbf{m} = \mbox{ \boldmath $\lambda$}_{\sigma(\mathbf{m})}.
\end{equation}

Let $\mathbf{s} \in \mathbf{G}$ and $c = e_\mathbf{s}$. We find necessary and sufficient conditions for part (b) of Corollary \ref{CorFBHfromProjNew} to hold for $c$; that is
\begin{equation}\label{EqnDeltaFC}
\Delta(\textit{f}\,(c)) = \textit{f}\,(c_{(1)})\textit{g}(c_{(2)(-1)}) \otimes \textit{f}\,(c_{(2)(0)}).
\end{equation}
First of all
\begin{eqnarray*}
c_{(1)} \otimes c_{(2)(-1)} \otimes c_{(2)(0)} & = & \sum_{\mathbf{m} \in \mathbf{G}} e_{\mathbf{s} - \mathbf{m}} \otimes \rho(e_\mathbf{m}) \\
& = & \sum_{\mathbf{m} \in \mathbf{G}} e_{\mathbf{s} - \mathbf{m}} \otimes \rho(\sum_{i \in \mathbb{Z}_{|\mathbf{m}|}} e_{\mathbf{m}, \lambda_\mathbf{m}^i}) \\
& = & \sum_{\mathbf{m} \in \mathbf{G}}\; \sum_{i \in \mathbb{Z}_{|\mathbf{m}|}} e_{\mathbf{s} - \mathbf{m}} \otimes \mbox{ \boldmath $\lambda$}_\mathbf{m}^i \otimes e_{\mathbf{m}, \lambda_\mathbf{m}^i} \\
& = & \sum_{\mathbf{m} \in \mathbf{G}} \; \sum_{i, \ell \in \mathbb{Z}_{|\mathbf{m}|}} e_{\mathbf{s} - \mathbf{m}} \otimes \mbox{\boldmath $\lambda$}_\mathbf{m}^i \otimes (\lambda_\mathbf{m}^{-i\ell}/|\mathbf{m}|) e_{\sigma^\ell(\mathbf{m})}.
\end{eqnarray*}
We interrupt our calculation to organize several technicalities.
\begin{Lemma}\label{LemmaFGDetails}
Let $A' = k[\mathcal{G}] \times k[\mathbf{U}]$ be the biproduct described above and suppose $(\textit{f}, \textit{g}) \in \mathcal{F}^\bullet_{k[\mathcal{G}], k[\mathbf{U}]}$. Then:
\begin{enumerate}
\item[{\rm (a)}] $\textit{f}^{-1} \circ \textit{g} : k[\mathbf{U}] \longrightarrow k[\mathcal{G}]$ is a Hopf algebra map and $\mbox{\boldmath $\theta$} \circ (\textit{f}^{-1} \circ \textit{g}) = \textit{f}^{-1} \circ \textit{g}$.
\item[{\rm (b)}] Let $\textit{g}_\lambda = (\textit{f}^{-1} \circ \textit{g})(\mbox{\boldmath $\lambda$})$. Then $\textit{g}_\lambda \in \mathcal{G}$, $\mbox{\boldmath $\theta$}(g_\lambda) = g_\lambda$, $g_\lambda^N = 1$, and $g_\lambda^{N/|\mathbf{m}|} = (\textit{f}^{-1} \circ \textit{g})(\mbox{\boldmath $\lambda$}_\mathbf{m})$ for all $\mathbf{m} \in \mathbf{G}$.
\item[{\rm (c)}] There exists $\alpha \in \widehat{\mathbf{G}}$ which is determined by $g_\lambda e_\mathbf{s} = \alpha(\mathbf{s})e_\mathbf{s}$ for all $\mathbf{s} \in \mathbf{G}$. Furthermore $\textit{g}_\lambda^\ell e_\mathbf{s} = \alpha(\mathbf{s})^\ell e_\mathbf{s}$ for all $\ell \in \mathbb{Z}$, $\mathbf{s} \in \mathbf{G}$, $\alpha \circ \sigma = \alpha$, and $\alpha(\mathbf{s})^N = 1$ for all $\mathbf{s} \in \mathbf{G}$.
\item[{\rm (d)}] For all $\mathbf{s}, \mathbf{m} \in \mathbf{G}$ there exists an integer $\ell(\mathbf{s}, \mathbf{m})$ which is determined by $0 \leq \ell(\mathbf{s}, \mathbf{m}) < |\mathbf{m}|$ and $(\alpha(\mathbf{s})/\alpha(\mathbf{m}))^{N/|\mathbf{m}|} = \lambda_{\mathbf{m}}^{\ell(\mathbf{s}, \mathbf{m})}$.
\item[{\rm (e)}] $\ell(\sigma(\mathbf{s}), \mathbf{m}) = \ell(\mathbf{s}, \mathbf{m}) = \ell(\mathbf{s}, \sigma(\mathbf{m}))$ for all $\mathbf{s}, \mathbf{m} \in \mathbf{G}$.
\end{enumerate}
\end{Lemma}

\pf
To show part (a) we first note by Corollary \ref{CorFBHfromProjNew} that $\textit{g} : k[\mathbf{U}] \longrightarrow k[\mathcal{G}]$ is a Hopf algebra map which satisfies $\mbox{\boldmath $\theta$} \circ \textit{g} = \textit{g}$ and $\textit{f} : k[\mathcal{G}] \longrightarrow k[\mathcal{G}]$ is an algebra automorphism which satisfies $\mbox{\boldmath $\theta$} \circ \textit{f} = \textit{f} \circ \mbox{\boldmath $\theta$}$ and $\epsilon \circ \textit{f} = \epsilon$. Thus $\textit{f}^{-1} \circ \textit{g}$ is an algebra map which satisfies $\mbox{\boldmath $\theta$} \circ (\textit{f}^{-1} \circ \textit{g}) = \textit{f}^{-1} \circ \textit{g}$ and $\epsilon \circ (\textit{f}^{-1} \circ \textit{g}) = \epsilon$.

Let $h \in k[\mathbf{U}]$. The condition $\mbox{\boldmath $\theta$} \circ \textit{g} = \textit{g}$ implies $\rho(\textit{g}(h)) = 1 \otimes \textit{g}(h)$; thus
$$
\textit{g}(h)_{(1)} \otimes \textit{g}(h)_{(2)(-1)} \otimes \textit{g}(h)_{(2)(0)} = \textit{g}(h)_{(1)} \otimes 1 \otimes \textit{g}(h)_{(2)}
$$
since $\textit{g}$ is a coalgebra map. Since $(\textit{f}, \textit{g})$ corresponds to an automorphism, $(\textit{f}^{-1}, \textit{g}') \in \mathcal{F}^\bullet_{k[\mathcal{G}], k[\mathbf{U}]}$ for some $\textit{g}'$ by (\ref{EqnIdASubL}), part (a) of Corollary \ref{CorFGParts}, and Theorem \ref{ThmEndAut}. Now we apply part (b) of Corollary \ref{CorFBHfromProjNew} to the pair $(\textit{f}^{-1}, \textit{g}')$ to conclude that $\textit{f}^{-1} \circ \textit{g}$ is a coalgebra map. This concludes the proof of part (a).

Note that $\textit{g}_\lambda \in \mathrm{G}(k[\mathcal{G}]) = \mathcal{G}$ since $\mbox{\boldmath $\lambda$} \in \mathbf{U} = \mathrm{G}(k[\mathbf{U}])$ and $\textit{f}^{-1} \circ \textit{g}$ is a coalgebra map. The remainder of part (b) easily follows by part (a). As for part (c), the existence of $\alpha$ follows from the bijection $\mathbf{G} \longrightarrow \widehat{\mathbf{G}}$ described just after (\ref{EqOrthBasis}). Now $\mbox{\boldmath $\theta$}(g_\lambda) = g_\lambda$ by part (b). The calculation
$$
g_\lambda e_{\sigma(\mathbf{s})} = \mbox{\boldmath $\theta$}(g_\lambda e_\mathbf{s}) = \mbox{\boldmath $\theta$}(\alpha(\mathbf{s})e_\mathbf{s}) = \alpha(\mathbf{s})e_{\sigma(\mathbf{s})}
$$
for all $\mathbf{s} \in \mathbf{G}$ shows that $\alpha = \alpha \circ \sigma^{-1}$, or equivalently
\begin{equation}\label{EqAlphaThetaCircTheta}
\alpha \circ \sigma = \alpha.
\end{equation}
Since $g_\lambda^N = 1$ the remainder of the proof of part (c) easily follows.

To show part (d), we recall $\alpha(\mathbf{n})^N = 1$ for all $\mathbf{n} \in \mathbf{G}$ by part (c). Let $\mathbf{s}, \mathbf{m} \in \mathbf{G}$. Then $(\alpha(\mathbf{s})/\alpha(\mathbf{m}))^{N/|\mathbf{m}|}$ is an $|\mathbf{m}|^{th}$  root of unity. Since $\lambda_{\mathbf{m}} = \lambda^{N/|\mathbf{m}|}$ is a primitive $|\mathbf{m}|^{th}$ root of unity, there exists a unique solution to $(\alpha(\mathbf{s})/\alpha(\mathbf{m}))^{N/|\mathbf{m}|} = \lambda_{\mathbf{m}}^\ell$ where $0 \leq \ell < |\mathbf{m}|$. We have shown that part (d) holds. Part (e) follows from part (d), (\ref{EqLambdaSigmaBF}), and the fact that $\alpha \circ \sigma = \alpha$.
\qed
\medskip

Using the preceding lemma we calculate
\begin{eqnarray*}
\lefteqn{\textit{f}\, (c_{(1)})\textit{g}(c_{(2)(-1)}) \otimes \textit{f}\, (c_{(2)(0)})} \\
& = & \sum_{\mathbf{m} \in \mathbf{G}} \; \sum_{i, \ell \in \mathbb{Z}_{|\mathbf{m}|}} \textit{f}\,(e_{\mathbf{s} - \mathbf{m}})\textit{g}(\mbox{\boldmath $\lambda$}_\mathbf{m}^i) \otimes (\lambda_\mathbf{m}^{-i\ell}/|\mathbf{m}|)\textit{f}\,(e_{\sigma^\ell(\mathbf{m})}) \\
& = & \sum_{\mathbf{m} \in \mathbf{G}} \; \sum_{i, \ell \in \mathbb{Z}_{|\mathbf{m}|}} \textit{f}\,(e_{\mathbf{s} - \mathbf{m}}((\textit{f}^{-1} \circ \textit{g})(\mbox{\boldmath $\lambda$}_\mathbf{m}))^i) \otimes (\lambda_\mathbf{m}^{-i\ell}/|\mathbf{m}|)\textit{f}\,(e_{\sigma^\ell(\mathbf{m})}) \\
& = & \sum_{\mathbf{m} \in \mathbf{G}} \; \sum_{i, \ell \in \mathbb{Z}_{|\mathbf{m}|}} \textit{f}\,(e_{\mathbf{s} - \mathbf{m}}\alpha(\mathbf{s} - \mathbf{m})^{(N/|\mathbf{m}|)i}) \otimes (\lambda_\mathbf{m}^{-i\ell}/|\mathbf{m}|)\textit{f}\,(e_{\sigma^\ell(\mathbf{m})}) \\
& = & \sum_{\mathbf{m} \in \mathbf{G}} \; \sum_{i, \ell \in \mathbb{Z}_{|\mathbf{m}|}} \alpha(\mathbf{s} - \mathbf{m})^{(N/|\mathbf{m}|)i})(\lambda_\mathbf{m}^{-i\ell}/|\mathbf{m}|)e_{\tau(\mathbf{s} - \mathbf{m})} \otimes e_{\tau(\sigma^\ell(\mathbf{m}))} \\
& = & \sum_{\mathbf{m} \in \mathbf{G}} \; \sum_{i, \ell \in \mathbb{Z}_{|\mathbf{m}|}} (\lambda_\mathbf{m}^{(\ell(\mathbf{s}, \mathbf{m})-\ell)i}/|\mathbf{m}|)e_{\tau(\mathbf{s} - \mathbf{m})} \otimes e_{\tau(\sigma^\ell(\mathbf{m}))} \\
& = & \sum_{\mathbf{m} \in \mathbf{G}} e_{\tau(\mathbf{s} - \mathbf{m})} \otimes e_{\tau(\sigma^{\ell(\mathbf{s}, \mathbf{m})}(\mathbf{m}))}.
\end{eqnarray*}
On the other hand
$$
\Delta(\textit{f}\,(c)) = \Delta (e_{\tau(\mathbf{s})}) = \sum_{\mathbf{m} \in \mathbf{G}} e_{\tau(\mathbf{s}) - \mathbf{m}} \otimes e_\mathbf{m} = \sum_{\mathbf{m} \in \mathbf{G}} e_{\tau(\mathbf{s}) - \tau(\mathbf{m})} \otimes e_{\tau(\mathbf{m})}.
$$
Therefore (\ref{EqnDeltaFC}) holds if and only if
\begin{equation}\label{EqnPartcHolds}
\sum_{\mathbf{m} \in \mathbf{G}} e_{\tau(\mathbf{s}) - \tau(\mathbf{m})} \otimes e_{\tau(\mathbf{m})} = \sum_{\mathbf{m} \in \mathbf{G}} e_{\tau(\mathbf{s} - \mathbf{m})} \otimes e_{\tau(\sigma^{\ell(\mathbf{s}, \mathbf{m})}(\mathbf{m}))}.
\end{equation}
Note that (\ref{EqnPartcHolds}) holds if and only if it holds when the index of summation runs over a $(\mbox{\boldmath $\theta$})$-orbit of $\mathbf{G}$ for all orbits.

Let $\mathbf{m} \in \mathbf{G}$ and $\mathcal{O} = \mathcal{O}_{\theta, \mathbf{m}}$. Since $\sigma^{\ell(\mathbf{s}, \mathbf{m})}$ permutes $\mathcal{O}$, and $\ell(\mathbf{s}, \mathbf{n}) = \ell(\mathbf{s}, \mathbf{m})$ for all $\mathbf{n} \in \mathcal{O}$ by part (e) of Lemma \ref{LemmaFGDetails}, the equation of (\ref{EqnPartcHolds}) when the index of summation runs over $\mathcal{O}$ can be written
\begin{equation}\label{EqnIndexOrbits}
\sum_{\mathbf{n} \in \mathcal{O}} e_{\tau(\mathbf{s}) - \tau(\sigma^{\ell(\mathbf{s}, \mathbf{m})}(\mathbf{n}))} \otimes e_{\tau(\sigma^{\ell(\mathbf{s}, \mathbf{m})}(\mathbf{n}))}  = \sum_{\mathbf{n} \in \mathcal{O}}e_{\tau(\mathbf{s} - \mathbf{n})} \otimes e_{\tau(\sigma^{\ell(\mathbf{s}, \mathbf{m})}(\mathbf{n}))}.
\end{equation}
Therefore part (b) of Corollary \ref{CorFBHfromProjNew} holds for all $b \in k[\mathcal{G}]$ if and only if
\begin{equation}\label{EqnCor7b}
\tau(\mathbf{s} - \mathbf{m}) = \tau(\mathbf{s}) - \tau(\sigma^{\ell(\mathbf{s}, \mathbf{m})}(\mathbf{m}))
\end{equation}
for all $\mathbf{s}, \mathbf{m} \in \mathbf{G}$.

We are now in a position to characterize the elements of $\mathrm{Aut}_\mathrm{Hopf}(A', \pi')$, or equivalently the elements of $\mathcal{F}^\bullet_{k[\mathcal{G}], k[\mathbf{U}]}$, in terms of certain elements of $\mathrm{Sym}(\mathbf{G})$ and $\widehat{\mathbf{G}}$. For $\mathbf{m} \in \mathbf{G}$ observe that $|\mathbf{m}|$ is also the length of the $(\sigma)$-orbit of $\mathbf{m}$. See (\ref{EqnmOrbit}).

Let $\mathbf{G}$ be any finite abelian group and $\sigma \in \mathrm{Aut}_\mathrm{Group}(\mathbf{G})$. Suppose $k^\times$ contain a primitive $|\sigma|^{th}$ root of unity $\lambda$. Let $\mathcal{F}_{\mathbf{G}, \sigma, \lambda}$ be the set of all pairs $(\tau, \alpha)$, where $\tau \in \mathrm{Sym}(\mathbf{G})$ and $\alpha \in \widehat{\mathbf{G}}$, such that:
\begin{enumerate}
\item[{\rm (F.a)}] $\tau \circ \sigma = \sigma \circ \tau$;
\item[{\rm (F.b)}] $\alpha \circ \sigma = \alpha$;
\item[{\rm (F.c)}] $\alpha^N = 1$, where $N$ is the order of $\sigma$; and
\item[{\rm (F.d)}] $\tau(\mathbf{s} - \mathbf{m}) = \tau(\mathbf{s}) - \tau(\sigma^{\ell_\alpha(\mathbf{s}, \mathbf{m})}(\mathbf{m}))$ for all $\mathbf{s}, \mathbf{m} \in \mathbf{G}$, where $\ell_\alpha(\mathbf{s}, \mathbf{m})$ is the integer determined by the conditions $0 \leq \ell_\alpha(\mathbf{s}, \mathbf{m}) < |\mathbf{m}|$ and  $(\alpha(\mathbf{s})/\alpha(\mathbf{m}))^{N/|\mathbf{m}|} = \lambda_{\mathbf{m}}^{\ell_\alpha(\mathbf{s}, \mathbf{m})}$.
\end{enumerate}

Observe that if $\beta \in \widehat{\mathbf{G}}$ and satisfies $\beta \circ \sigma = \beta$ then $\beta \circ \tau \in \widehat{\mathbf{G}}$ by (F.b) and (F.e). If $(\tau, \alpha) \in \mathcal{F}_{\mathbf{G}, \sigma, \lambda}$ then
\begin{enumerate}
\item[{\rm (F.e)}]
$\ell_\alpha(\sigma(\mathbf{s}), \mathbf{m}) = \ell_\alpha(\mathbf{s}, \mathbf{m}) = \ell_\alpha(\mathbf{s}, \sigma(\mathbf{m}))$
\end{enumerate}
for all $\mathbf{s}, \mathbf{m} \in \mathbf{G}$ for all $\mathbf{s}, \mathbf{m} \in \mathbf{G}$ since $\alpha \circ \sigma = \alpha$ and $|\sigma(\mathbf{m})| = |\mathbf{m}|$ for all $\mathbf{m} \in \mathbf{G}$.
\begin{Theorem}\label{ThmTwoFFs}
Let $A' = k[\mathcal{G}] \times k[\mathbf{U}]$ be the biproduct described above, where $\mathbf{U} = (\mbox{\boldmath $\lambda$})$, $\mbox{\boldmath $\theta$} = F_\sigma$, and $\lambda \in k^\times$ is a primitive $|\sigma|^{th}$ root of unity. Then:
\begin{enumerate}
\item[{\rm (a)}] There is a bijection $\mathcal{F}^\bullet_{k[\mathcal{G}], k[\mathbf{U}]} \longrightarrow \mathcal{F}_{\mathbf{G}, \sigma, \lambda}$ given by $(\textit{f}, \textit{g}) \mapsto (\tau, \alpha)$, where $\textit{f} = F_\tau$ and $\alpha$ is determined by $\textit{g}(\mbox{\boldmath $\lambda$})e_{\tau(\mathbf{s})} = \alpha(\mathbf{s})e_{\tau(\mathbf{s})}$ for all $\mathbf{s} \in \mathbf{G}$.
\item[{\rm (b)}] With the identification of $\mathcal{F}_{\mathbf{G}, \sigma, \lambda}$ and $\mathcal{F}^\bullet_{k[\mathcal{G}], k[\mathbf{U}]}$ via the bijection of part (a), the element $F \in \mathrm{Aut}_\mathrm{Hopf}(A', \pi')$ corresponding to $(\tau, \alpha)$ as in part (a) of Theorem \ref{ThmEndAut} is determined by
    $$
    F(e_\mathbf{s} \times \mbox{\boldmath $\lambda$}^\ell) = \alpha(\mathbf{s})^\ell e_{\tau(\mathbf{s})} \times \mbox{\boldmath $\lambda$}^\ell
    $$
    for all $\mathbf{s} \in \mathbf{G}$ and $0 \leq \ell < |\mathbf{U}|$.
\end{enumerate}
\end{Theorem}

\pf
Let $(\textit{f}, \textit{g}) \in \mathcal{F}^\bullet_{k[\mathcal{G}], k[\mathbf{U}]}$. Then $\textit{f} = F_\tau$ for a unique $\tau \in \mathrm{Sym}(\mathbf{G})$ by (\ref{EqBoldThetaTau}). Let $\textit{g}_\lambda = (\textit{f}^{-1} \circ \textit{g})(\mbox{\boldmath $\lambda$})$. There exists $\alpha \in \widehat{\mathbf{G}}$ such that $g_\lambda e_\mathbf{s} = \alpha (\mathbf{s})e_\mathbf{s}$ for all $\mathbf{s} \in \mathbf{G}$ by part (c) of Lemma \ref{LemmaFGDetails}. Since $\textit{f}$ is an algebra map $\textit{f}(g_\lambda) \textit{f}(e_\mathbf{s}) = \textit{f}(g_\lambda e_\mathbf{s}) = \alpha (\mathbf{s})\textit{f}(e_\mathbf{s})$ for all $\mathbf{s} \in \mathbf{G}$. Thus $\textit{g}(\mbox{\boldmath $\lambda$})e_{\tau(\mathbf{s})} = \alpha(\mathbf{s})e_{\tau(\mathbf{s})}$ for all $\mathbf{s} \in \mathbf{G}$. The preceding equation determines $\alpha$. That $(\tau, \alpha) \in \mathcal{F}_{\mathbf{G}, \sigma, \lambda}$ follows by (\ref{EqTauThetaCommute}), (\ref{EqAlphaThetaCircTheta}), parts (c)--(e) of Lemma \ref{LemmaFGDetails}, and (\ref{EqnCor7b}). The association $(\textit{f}, \textit{g}) \mapsto (\tau, \alpha)$ defines a function $\mathbf{f} : \mathcal{F}^\bullet_{k[\mathcal{G}], k[\mathbf{U}]} \longrightarrow \mathcal{F}_{\mathbf{G}, \sigma, \lambda}$.

Now let $(\tau, \alpha) \in \mathcal{F}_{\mathbf{G}, \sigma, \lambda}$. Then $\textit{f} = F_\tau$ is an algebra automorphism of $k[\mathcal{G}]$ by Lemma \ref{LemmaSymBfG}. Let ${\sf g} \in k[\mathcal{G}]$ be defined by ${\sf g}e_\mathbf{s} = \alpha(\mathbf{s})e_\mathbf{s}$ for all $\mathbf{s} \in \mathbf{G}$. Then ${\sf g} \in \mathcal{G}$; see (\ref{EqngmSumChar}) and the subsequent remark. Since $\textit{f}$ is a algebra map $\textit{f}({\sf g})e_{\tau(\mathbf{s})} = \alpha(\mathbf{s})e_{\tau(\mathbf{s})}$, or equivalently $\textit{f}({\sf g})e_\mathbf{s} = (\alpha \circ \tau^{-1})(\mathbf{s})e_\mathbf{s}$, for all $\mathbf{s} \in \mathbf{G}$. By virtue of (F.a) and the remark preceding the statement of the theorem $\alpha \circ \tau^\ell \in \widehat{\mathbf{G}}$ for all $\ell \geq 0$. Since $\tau$ has finite order $\alpha \circ \tau^{-1} \in \widehat{\mathbf{G}}$. This means $\textit{f}({\sf g}) \in \mathcal{G}$ as well. Now ${\sf g}^N = 1$ by (F.c). Thus $\textit{f}({\sf g})^N = 1$.

Since $\mbox{\boldmath $\lambda$}$ has order $N$ there is a Hopf algebra map $\mathit{g} : k[\mathbf{U}] \longrightarrow k[\mathcal{G}]$ determined by $\textit{g}(\mbox{\boldmath $\lambda$}) = \textit{f}({\sf g})$. To show that $(\textit{g}, \textit{f}) \in \mathcal{F}^\bullet_{k[\mathcal{G}], k[\mathbf{U}]}$ we need only show that the conditions (a) -- (e) of Corollary \ref{CorFBHfromProjNew} are satisfied for $\textit{f}$ and $\textit{g}$. Since $\mathcal{G}$ is abelian (d) is satisfied. Using the fact that $\mbox{\boldmath $\theta$}$ is an algebra map, (F.a), and (F.b), we compute
$$
\mbox{\boldmath $\theta$}(\textit{g}(\mbox{\boldmath $\lambda$}))e_{\tau(\sigma(\mathbf{s}))} = \mbox{\boldmath $\theta$}(\textit{g}(\mbox{\boldmath $\lambda$})e_{\tau(\mathbf{s})}) = \alpha(\mathbf{s})e_{\sigma(\tau(\mathbf{s}))} = \alpha (\sigma(\mathbf{s}))e_{\tau(\sigma(\mathbf{s}))}
$$
for all $\mathbf{s} \in \mathbf{G}$. Therefore $\mbox{\boldmath $\theta$}(\textit{g}(\mbox{\boldmath $\lambda$}))e_{\tau(\mathbf{s})} = \alpha(\mathbf{s})e_{\tau(\mathbf{s})} = \textit{g}(\mbox{\boldmath $\lambda$})e_{\tau(\mathbf{s})}$ for all $\mathbf{s} \in \mathbf{G}$ which implies $\mbox{\boldmath $\theta$}(\textit{g}(\mbox{\boldmath $\lambda$})) = \textit{g}(\mbox{\boldmath $\lambda$})$. Since $\mbox{\boldmath $\lambda$}$ generates $k[\mathbf{U}]$ as an algebra $\mbox{\boldmath $\theta$} \circ \textit{g} = \textit{g}$. Therefore condition (e) is met.

We observe $\tau(\mathbf{0}) = \mathbf{0}$ by (F.d); note that $\ell_\alpha(\mathbf{s}, \mathbf{s}) = 0$ for all $\mathbf{s} \in \mathbf{G}$. The first equation is equivalent to $\epsilon \circ \textit{f} = \epsilon$. Thus condition (a) is fulfilled. Lemma \ref{LemmaSymBfG} and (F.a) account for the fact that condition (c) holds. As for (b), the reader is left to showing it is fulfilled by retracing the calculations following (\ref{EqnDeltaFC}). We have shown $(\textit{g}, \textit{f}) \in \mathcal{F}^\bullet _{k[\mathcal{G}], k[\mathbf{U}]}$. The association $(\tau, \alpha) \mapsto (\textit{f}, \textit{g})$ defines a function $\mathbf{g} : \mathcal{F}_{\mathbf{G}, \sigma, \lambda}  \longrightarrow \mathcal{F}^\bullet_{k[\mathcal{G}], k[\mathbf{U}]}$. It is an easy exercise to show that $\mathbf{f}$ and $\mathbf{g}$ are inverse functions. This completes our proof of part (a).

Let $F$ be as in part (b). Since $F_L = f$ and $F_R = g$ are algebra maps and $k[\mathcal{G}]$ is commutative
\begin{eqnarray*}
 F(e_\mathbf{s} \times \mbox{\boldmath $\lambda$}^\ell) & = & F_L(e_\mathbf{s})F_R(\mbox{\boldmath $\lambda$}^\ell) \times \mbox{\boldmath $\lambda$}^\ell \\
 & = & g(\mbox{\boldmath $\lambda$})^\ell f(e_\mathbf{s}) \times  \mbox{\boldmath $\lambda$}^\ell \\
 & = & \alpha(\mathbf{s})^\ell e_{\tau(\mathbf{s})} \times \mbox{\boldmath $\lambda$}^\ell.
\end{eqnarray*}
Since the $e_\mathbf{s} \times \mbox{\boldmath $\lambda$}^\ell$'s described in part (b) form a basis for $A'$, the proof of part (b) is complete.
\qed
\section{The group $\Gamma (\mathbf{G}, \lambda, \sigma)$}\label{SecGammaGLambda}
In this section we recast the study of the group $\mathrm{Aut}_\mathrm{Hopf}(A', \pi')$, where $A' = k[\mathcal{G}] \times k[\mathbf{U}]$ is the biproduct of Theorem \ref{ThmTwoFFs}, in terms of permutations and characters. We continue with the notation of the previous section.

Let $\mathcal{G}$ be any finite abelian group and $\sigma \in \mathrm{Aut}_{\mathrm{Group}}(\mathbf{G})$. Assume that $k$ contains a primitive $N^{th}$ root of unity $\lambda$, where $N = |\sigma|$ is the order of $\sigma$. Then we define
$$
\Gamma (\mathbf{G}, \lambda, \sigma) = \{\tau \, | \, (\tau, \alpha) \in \mathcal{F}_{\mathbf{G}, \sigma, \lambda} \;\; \mbox{for some $\alpha \in \widehat{\mathbf{G}}$}\}
$$
We leave it to the reader to show that $\Gamma (\mathbf{G}, \lambda, \sigma)$ is a subgroup of $\mathrm{Sym}(\mathbf{G})$. Observe that group $\mathrm{Aut}_\sigma(\mathbf{G})$ of automorphisms of $\mathbf{G}$ which commute with $\sigma$ is a subgroup of $\Gamma (\mathbf{G}, \lambda, \sigma)$. The function
\begin{equation}\label{EqNuGSigma}
\nu_{\mathbf{G}, \sigma} : \mathcal{F}_{\mathbf{G}, \sigma, \lambda} \longrightarrow \Gamma(\mathbf{G}, \lambda, \sigma) \qquad \qquad (\tau, \alpha) \mapsto \tau
\end{equation}
is a bijection. By definition $\nu_{\mathbf{G}, \sigma}$ is onto. To see that it is one-one suppose $(\tau, \alpha), (\tau, \alpha') \in \mathcal{F}_{\mathbf{G}, \sigma, \lambda}$. Since $\tau$ is one-one $\ell_\alpha = \ell_{\alpha'}$ by (F.d). Fix $\mathbf{m} \in \mathbf{G}$. Then by (F.d) again
$$
\left(\alpha(\mathbf{s})/\alpha(\mathbf{m})\right)^{N/|\mathbf{m}|} = \left(\alpha'(\mathbf{s})/\alpha'(\mathbf{m})\right)^{N/|\mathbf{m}|},
$$
or equivalently
$$
\left(\alpha(\mathbf{s})/\alpha'(\mathbf{s})\right)^{N/|\mathbf{m}|} = \left(\alpha(\mathbf{m})/\alpha'(\mathbf{m})\right)^{N/|\mathbf{m}|},
$$
for all $\mathbf{s} \in \mathbf{G}$. Therefore $\beta \in \widehat{\mathbf{G}}$ defined by $\beta(\mathbf{s}) = \left(\alpha(\mathbf{s})/\alpha'(\mathbf{s})\right)^{N/|\mathbf{m}|}$ for all $\mathbf{s} \in \mathbf{G}$ is constant. Hence $\beta = 1$ which means $\left(\alpha(\mathbf{s})/\alpha'(\mathbf{s})\right)^{N/|\mathbf{m}|} = 1$. Since $N$ is the least common multiple of the lengths of the $(\sigma)$-orbits of $\mathbf{G}$ it follows that $\alpha(\mathbf{s})/\alpha'(\mathbf{s}) = 1$. We have shown $\alpha = \alpha'$ and therefore $\nu_{\mathbf{G}, \sigma}$ is one-one.
\begin{Theorem}\label{ThmAPrtimeBoldG}
Let $A = k[\mathcal{G}] \times k[G]$ be the biproduct of Theorem \ref{ThmFromProjNew}. Then $k$ contains a primitive $N^{th}$ root of unity $\lambda$, where $N = |\theta|$. Assume further that $k$ contains a primitive $|\mathcal{G}|^{th}$ root of unity. Let $A' = k[\mathcal{G}] \times k[\mathbf{U}]$ and let $\mathbf{G}$ be the additive version of $\mathcal{G}$. Then:
\begin{enumerate}
\item[{\rm (a)}] There is an isomorphism of groups
$$
\Phi_{A', \mathbf{G}} : \mathrm{Aut}_{\mathrm{Hopf}}(A', \pi') \longrightarrow \Gamma (\mathbf{G}, \lambda, \sigma)
$$
for some $\sigma \in \mathrm{Aut}_{\mathrm{Group}}(\mathbf{G})$.
\item[{\rm (b)}] Let $F \in \mathrm{Aut}_{\mathrm{Hopf}}(A', \pi')$. Then $F_L$ is a coalgebra map if and only if $\Phi_{A', \mathbf{G}}(F) \in \mathrm{Aut}_{\mathrm{Group}}(\mathbf{G})$.
\end{enumerate}
\end{Theorem}

\pf
$\Phi_{A', \mathbf{G}}$ is the composition of bijections
$$
\mathrm{Aut}_{\mathrm{Hopf}}(A', \pi') \longrightarrow \mathcal{F}^\bullet_{k[\mathcal{G}], k[\mathbf{U}]} \longrightarrow \mathcal{F}_{\mathbf{G}, \sigma, \lambda} \longrightarrow \Gamma (\mathbf{G}, \lambda, \sigma),
$$
where the first is the inverse of the one of part (a) of Theorem \ref{ThmEndAut}, the second is the one of part (a) of Theorem \ref{ThmTwoFFs}, and the third is that of (\ref{EqNuGSigma}). Let $F \in \mathrm{Aut}_{\mathrm{Hopf}}(A', \pi')$. The calculation $F \mapsto (F_L, F_R) \mapsto (F_\tau, \alpha) \mapsto \tau$, where $F_L = F_\tau$ shows that $\Phi_{A', \mathbf{G}}(F) = \tau$. Note that $F_L$ is a coalgebra map if and only if $\tau \in \mathrm{Aut}_{\mathrm{Group}}(\mathbf{G})$ by Lemma \ref{LemmaSymBfG}.
\qed
\medskip

We continue with our general discussion of $\Gamma(\mathbf{G}, \lambda, \sigma)$. Let $\tau \in \Gamma (\mathbf{G}, \lambda, \sigma)$ and let $\alpha \in \widehat{\mathbf{G}}$ satisfy $(\tau, \alpha) \in \mathcal{F}_{\mathbf{G}, \sigma, \lambda}$. Since $\alpha^N = 1$ it follows that $\mathrm{Im}(\alpha) = (\omega)$, where $\omega \in k^\times$ is a primitive $r^{th}$ root of unity where $r | N$. For $i \in \mathbb{Z}_r$ let $\mathbf{G}_i = \alpha^{-1}(\omega^i)$. Then $\mathbf{G}_0 = \mathrm{Ker}(\alpha)$ and the $\mathbf{G}_i$'s are the cosets of $\mathbf{G}_0$. Note that $\mathbf{G}_i +  \mathbf{G}_j = \mathbf{G}_{i + j}$ for all $i, j \in \mathbb{Z}_r$, where addition of subscripts takes place in $\mathbb{Z}_r$.  Also
\begin{equation}\label{EqSigmaGiGi}
\sigma(\mathbf{G}_i) = \mathbf{G}_i
\end{equation}
for all $i \in \mathbb{Z}_r$ since $\alpha \circ \sigma = \alpha$. A consequence of (\ref{EqSigmaGiGi}) is that each $\mathbf{G}_i$ is the union of $(\sigma)$-orbits of $\mathbf{G}$. For all $i \in \mathbb{Z}_r$ observe that (F.d) implies
\begin{equation}\label{EqTauSMinusM}
\tau(\mathbf{s} - \mathbf{s}') = \tau(\mathbf{s}) - \tau(\mathbf{s}') \;\; \mbox{for all $\mathbf{s}, \mathbf{s}' \in \mathbf{G}_i$}
\end{equation}
as $\ell_\alpha (\mathbf{s}, \mathbf{s}') = 0$ for such $\mathbf{s}, \mathbf{s}'$. As a result
\begin{equation}\label{EqTauSMinusMG0}
\tau(\mathbf{s} - \mathbf{m}) = \tau(\mathbf{s}) - \tau(\mathbf{m}) \;\; \mbox{and} \;\; \tau(\mathbf{s} + \mathbf{m}) = \tau(\mathbf{s}) + \tau(\mathbf{m})
\end{equation}
for all $\mathbf{s} \in \mathbf{G}$ and  $\mathbf{m} \in \mathbf{G}_0$.

\pf
Let $\mathbf{s} \in \mathbf{G}$. Then $\mathbf{s} \in \mathbf{G}_i$ for some $i \in \mathbb{Z}_r$. Let $\mathbf{m} \in \mathbf{G}_0$. Then $\mathbf{s}, \mathbf{s} - \mathbf{m} \in \mathbf{G}_i$ and thus $\tau(\mathbf{m}) = \tau(\mathbf{s} - (\mathbf{s} - \mathbf{m})) = \tau(\mathbf{s}) - \tau(\mathbf{s} - \mathbf{m})$ by (\ref{EqTauSMinusM}) from which the first equation of (\ref{EqTauSMinusMG0}) follows. As for the second we use the first to calculate $\tau(\mathbf{s}) = \tau((\mathbf{s} + \mathbf{m}) - \mathbf{m}) =  \tau(\mathbf{s} + \mathbf{m}) - \tau(\mathbf{m})$ from which the second equation of (\ref{EqTauSMinusMG0}) follows.
\qed
\medskip

Let $\mathbf{s}, \mathbf{s}' \in \mathbf{G}$ and $\mathbf{m}, \mathbf{m}' \in \mathbf{G}_0$. It follows from (\ref{EqTauSMinusMG0}) and (F.d) that
\begin{eqnarray*}
\tau((\mathbf{s} + \mathbf{m}) - (\mathbf{s}' + \mathbf{m}')) & = & \tau((\mathbf{s} - \mathbf{s}') + (\mathbf{m} - \mathbf{m}')) \\
& = & \tau(\mathbf{s}) - \tau(\sigma^{\ell_\alpha(\mathbf{s}, \mathbf{s}')}(\mathbf{s}'))  + \tau(\mathbf{m}) - \tau(\mathbf{m}') \\
& = & \tau(\mathbf{s} + \mathbf{m}) - \tau(\sigma^{\ell_\alpha(\mathbf{s}, \mathbf{s}')}(\mathbf{s}') + \mathbf{m}')
\end{eqnarray*}
on the one hand and
\begin{eqnarray*}
\tau((\mathbf{s} + \mathbf{m}) - (\mathbf{s}' + \mathbf{m}')) & = & \tau(\mathbf{s} + \mathbf{m}) - \tau(\sigma^{\ell_\alpha(\mathbf{s} + \mathbf{m}, \mathbf{s}' + \mathbf{m}')}(\mathbf{s}' + \mathbf{m}'))
\end{eqnarray*}
on the other. Since $\tau$ is one-one
\begin{equation}\label{EqLSSPrime}
\sigma^{\ell_\alpha(\mathbf{s} + \mathbf{m}, \mathbf{s}' + \mathbf{m}')}(\mathbf{s}' + \mathbf{m}') = \sigma^{\ell_\alpha(\mathbf{s}, \mathbf{s}')}(\mathbf{s}') + \mathbf{m}'.
\end{equation}
\begin{Lemma}\label{LemmaTauAlphaGSigma}
Let $\tau \in \mathcal{F}_{\mathbf{G}, \sigma, \lambda}$, let $\mathbf{s}, \mathbf{s}' \in \mathbf{G}$, let $\mathbf{m}, \mathbf{m}' \in \mathbf{G}_0$, and let $r$ be the order of $\alpha$. Then:
\begin{enumerate}
\item[{\rm (a)}] $\ell_\alpha(\mathbf{s}, \mathbf{m}) = 0$.
\item[{\rm (b)}] $\ell_\alpha(\mathbf{s}, \mathbf{s}') = 0$ when $\mathbf{s}, \mathbf{s}' \in \mathbf{G}_i$ for some $i \in \mathbb{Z}_r$.
\item[{\rm (c)}] $\sigma^{\ell_\alpha(\mathbf{s} + \mathbf{m}, \mathbf{s}')}(\mathbf{s}') = \sigma^{\ell_\alpha(\mathbf{s}, \mathbf{s}' - \mathbf{m}')}(\mathbf{s}' - \mathbf{m}') + \mathbf{m}'$.
\item[{\rm (d)}] $\ell_\alpha(\mathbf{s} + \mathbf{m}, \mathbf{s}') = \ell_\alpha(\mathbf{s}, \mathbf{s}')$.
\item[{\rm (e)}] $\alpha(\mathbf{s})^{N/|\mathbf{m}|} = 1$; thus $|\mathbf{m}|$ divides $N/r$.
\item[{\rm (f)}] If the lengths of $(\sigma)$-orbits of $\mathbf{G}$ are the lengths of the $(\sigma)$-orbits of $\mathbf{G}_0$ then $\alpha = 1$, that is $\tau \in \mathrm{Aut}_\mathrm{Group}(\mathbf{G})$.
\end{enumerate}
\end{Lemma}

\pf
$\tau(\mathbf{s} - \mathbf{s}') = \tau(\mathbf{s}) - \tau(\mathbf{s}')$ if and only if $\ell_\alpha(\mathbf{s}, \mathbf{s}') = 0$ by (F.d) since $\tau$ is one-one. Thus part (a) follows by (\ref{EqTauSMinusMG0}). We have noted part (b) holds. Part (c) is a reformulation of (\ref{EqLSSPrime}). As for part (d) we note that part (c) is  $\sigma^{\ell_\alpha(\mathbf{s} + \mathbf{m}, \mathbf{s}')}(\mathbf{s}') = \sigma^{\ell_\alpha(\mathbf{s}, \mathbf{s}')}(\mathbf{s}')$ when we take $\mathbf{m}' = \mathbf{0}$. Thus part (d) follows by (F.d).

By part (a) and (F.d) we have $\alpha(\mathbf{s})^{N/|\mathbf{m}|} = 1$. Therefore $r$ divides $N/|\mathbf{m}|$ which means $|\mathbf{m}|$ divides $N/r$. We have established part (e).

Assume the hypothesis of part (f). Since $N$ is the least common multiple of the lengths of the $(\sigma)$-orbits of $\mathbf{G}$, by part (e) it follows that $N$ divides $N/r$. Therefore $r = 1$ which means $\alpha = 1$. Consequently $\mathbf{G} = \mathbf{G}_0$ and therefore $\ell_\alpha(\mathbf{s}, \mathbf{s}') = 0$ by part (a). We have noted this equation is equivalent to $\tau(\mathbf{s} - \mathbf{s}') = \tau(\mathbf{s}) - \tau(\mathbf{s}')$. The latter implies $\tau$ is a group homomorphism.
\qed
\medskip

Part (f) of  the preceding lemma suggests examination of the relationship between the lengths of the $(\sigma)$-orbits of $\mathbf{G}$ and those of the $(\sigma)$-orbits of $\mathbf{G}_0$. Let $\mathbf{s} \in \mathbf{G}$. Then $\mathbf{s} \in \mathbf{G}_i$ for some $i \in \mathbb{Z}_r$. Since $\sigma(\mathbf{s}) \in \mathbf{G}_i$ by (\ref{EqSigmaGiGi}), we have $\sigma(\mathbf{s}) = \mathbf{s} + \mathbf{m}$ for some $\mathbf{m} \in \mathbf{G}_0$. Thus
$$
\sigma^\ell(\mathbf{s}) = \mathbf{s} + \mathbf{m} + \cdots + \sigma^{\ell - 1}(\mathbf{m})
$$
for all $\ell \geq 1$. In particular $\mathbf{s} = \sigma^{|\mathbf{s}|}(\mathbf{s}) = \mathbf{s} + \mathbf{m} + \cdots + \sigma^{|\mathbf{s}| - 1}(\mathbf{m})$ which means $\mathbf{m} + \cdots + \sigma^{|\mathbf{s}| - 1}(\mathbf{m}) = 0$. Therefore $\sigma^{|\mathbf{s}|}(\mathbf{m}) = \mathbf{m}$ from which we conclude $|\mathbf{m}|$ divides $|\mathbf{s}|$. We have shown $|\mathbf{s}| = s|\mathbf{m}|$ for some positive integer.

Let $\mathbf{f} = \mathbf{m} + \cdots + \sigma^{|\mathbf{m}| - 1}(\mathbf{m})$. Then $\mathbf{f} \in \mathbf{G}_0$, $\sigma(\mathbf{f}) = \mathbf{f}$, and $\sigma^{|\mathbf{m}|}(\mathbf{s}) = \mathbf{s} + \mathbf{f}$. The latter implies $\sigma^{\ell|\mathbf{m}|}(\mathbf{s}) = \mathbf{s} + \ell\mathbf{f}$ for all non-negative integers $\ell$. Therefore $s\mathbf{f} = 0$ and if $\ell$ is a positive integer such that $\ell\mathbf{f} = 0$ then $s$ divides $\ell$. To summarize:
\begin{Lemma}\label{LemmaOrbitLenghtsMU}
Let $\mathbf{s} \in \mathbf{G}$. Then $\sigma(\mathbf{s}) = \mathbf{s} + \mathbf{m}$ for some $\mathbf{m} \in \mathbf{G}_0 = \mathrm{Ker}(\alpha)$. Set $\mathbf{f} = \mathbf{m} + \cdots + \sigma^{|\mathbf{m}|-1}(\mathbf{m})$. Then:
\begin{enumerate}
\item[{\rm (a)}] $\mathbf{f} \in \mathbf{G}_0$ and $\sigma(\mathbf{f}) = \mathbf{f}$;
\item[{\rm (b)}] $|\mathbf{s}| = s|\mathbf{m}|$, where $s$ is the order of $\mathbf{f}$.
\end{enumerate}
\qed
\end{Lemma}
\begin{Theorem}\label{ThmSigmaNoFix}
Let $\mathbf{G}$ be a finite abelian group, suppose $\sigma \in \mathrm{Aut}_\mathrm{Group}(\mathbf{G})$ and has order $N$, and suppose $k$ has a primitive $N^{th}$ root of unity $\lambda$. Then
$\mathrm{Aut}_\sigma(\mathbf{G}) = \Gamma (\mathbf{G}, \lambda, \sigma)$ if any of the following hold:
\begin{enumerate}
\item[{\rm (a)}] $\sigma$ has a unique fixed point.
\item[{\rm (b)}] The orders of $\sigma$ and $|\mathbf{G}|$ are relatively prime.
\item[{\rm (c)}] $N$ is prime $N^2$ does not divide $|\mathbf{G}|$.
\end{enumerate}
\end{Theorem}

\pf
We have noted that $\mathrm{Aut}_\sigma(\mathbf{G}) \leq \Gamma (\mathbf{G}, \lambda, \sigma)$. Let $\tau \in \Gamma (\mathbf{G}, \lambda, \sigma)$. Then $(\tau, \alpha) \in \mathcal{F}_{\mathbf{G}, \sigma, \lambda}$ for some $\alpha \in \widehat{\mathbf{G}}$.

First of all assume part (a) or (b) holds. By Lemma \ref{LemmaOrbitLenghtsMU} the hypothesis of part (f) of Lemma \ref{LemmaTauAlphaGSigma} holds. Therefore $\tau \in \mathrm{Aut}_\sigma(\mathbf{G})$.

Now assume that part (c) holds and set $N = p$. Suppose $\tau \not \in \mathrm{Aut}_\sigma(\mathbf{G})$. Then $\alpha \neq 1$ by (F.d) which means $\alpha$ has order $p$ by (F.c). Since a $(\sigma)$-orbit of $\mathbf{G}$ has length $1$ or $p$, $\mathbf{G}_0$ consists of fixed points of $\sigma$ by part (f) of Lemma \ref{LemmaTauAlphaGSigma} again.

Let $\mathbf{s} \in \mathbf{G}$. Then $\sigma(\mathbf{s}) = \mathbf{s} + \mathbf{m}$ for some $\mathbf{m} \in \mathbf{G}_0$ by Lemma \ref{LemmaOrbitLenghtsMU}. Since $\sigma(\mathbf{m}) = \mathbf{m}$ it follows that $\mathbf{s} = \sigma^p(\mathbf{s}) = \mathbf{s} + p\mathbf{m}$ and therefore $p\mathbf{m} = \mathbf{0}$. Since $|\mathbf{G}| = p|\mathbf{G}_0|$ it follows $p$ does not divide $|\mathbf{G}_0|$ since $p^2$ does not divide $|\mathbf{G}|$. Therefore $\mathbf{m} = \mathbf{0}$ which means $\sigma(\mathbf{s}) = \mathbf{s}$. We have shown $\sigma = \mathrm{Id}_G$, a contradiction. Therefore $\tau \in \mathrm{Aut}_\sigma(\mathbf{G})$ after all.
\qed
\medskip

Closer examination of the proof of part (a) of the preceding theorem reveals:
\begin{Cor}\label{CorTau0OnklyFixed}
Let $(\tau, \alpha) \in \mathcal{F}_{\mathbf{G}, \sigma, \lambda}$. If $\mathbf{0}$ is the only fixed point of $\sigma$ in $\mathbf{G}_0$ then $\tau \in \mathrm{Aut}_\mathrm{Group}(\mathbf{G})$.
\qed
\end{Cor}

We construct examples where $\mathrm{Aut}_\sigma(\mathbf{G}) < \Gamma (\mathbf{G}, \lambda, \sigma)$. To do this we examine what it means for $(\tau, \alpha) \in \mathcal{F}_{\mathbf{G}, \, \sigma}$, where $\sigma$ has prime order $N = p$ and $\tau \not \in \mathrm{Aut}_\sigma(\mathbf{G})$. At one point we find it convenient to place a restriction on $\tau$. We continue with the preceding notation and results of this section without particular reference for the most part.

Let $\mathbf{G}$ be a finite abelian group, $\sigma \in \mathrm{Aut}_\mathrm{Group}(\mathbf{G})$ and has prime order $p$, and $(\tau, \alpha) \in \mathcal{F}_{\mathbf{G}, \, \sigma}$.

Suppose $\tau \not \in \mathrm{Aut}_\mathrm{Group}(\mathbf{G})$. Then $\alpha \neq 1$. Since $\alpha^p = 1$ it follows that $\alpha$ has order $p$. We take $\omega = \lambda$. Thus $\mathbf{G}_i = \alpha^{-1}(\lambda^i)$ for all $0 \leq i \leq p-1$.

The $(\sigma)$-orbits of $\mathbf{G}$ have lengths $1$ or $p$. Now each $\mathbf{G}_i$ is the union of $(\sigma)$-orbits since $\sigma(\mathbf{G}_i) = \mathbf{G}_i$. Since $\tau \not \in \mathrm{Aut}_\mathrm{Group}(\mathbf{G})$, by part (f) of Lemma \ref{LemmaTauAlphaGSigma} it follows that $\mathbf{G}_0$ consists of fixed points of $\sigma$. The set $\mathbf{F}$ of fixed points of $\sigma$ is a subgroup of $\mathbf{G}$. Since $\mathbf{G}_0 \subseteq \mathbf{F}$ and $[\mathbf{G} : \mathbf{G}_0]$ is prime, either $\mathbf{F} =\mathbf{G}_0$ or $\mathbf{F} = \mathbf{G}$. Thus $\mathbf{F} =\mathbf{G}_0$ since $\sigma \neq \mathrm{Id}_\mathbf{G}$. Since $\tau$ and $\sigma$ commute and $\tau$ is bijective, $\tau(\mathbf{F}) = \mathbf{F}$. We have shown $\tau(\mathbf{G}_0) = \mathbf{G}_0$.  The equation $\mathbf{F} = \mathbf{G}_0$ implies that $\mathbf{G}_i$ is the union of $(\sigma)$-orbits of length $p$, where $1 \leq i \leq p-1$.

Each $\tau(\mathbf{G}_i)$ is a coset of $\tau(\mathbf{G}_0) = \mathbf{G}_0$ by (\ref{EqTauSMinusMG0}). Therefore there is a $\eta \in \mathrm{Sym}(\mathbb{Z}_p)$ such that
\begin{equation}\label{EqEtaTauGI}
\eta(0) = 0 \;\; \mbox{and} \;\; \tau(\mathbf{G}_i) = \mathbf{G}_{\eta(i)}
\end{equation}
for all $0 \leq i \leq p-1$.

Now $\mathbf{G}/\mathbf{G}_0$ is cyclic of order $p$. Fix $\mathbf{s} \in \mathbf{G}_1$. Then $\mathbf{G}_i = i\mathbf{s} + \mathbf{G}_0$ for all $0 \leq i \leq p-1$. Since $\sigma(\mathbf{G}_1) = \mathbf{G}_1$,
\begin{equation}\label{EqSigmaSSPlusM}
\sigma(\mathbf{s}) = \mathbf{s} + \mathbf{m}
\end{equation}
for some $\mathbf{m} \in \mathbf{G}_0$. Since $\sigma \neq \mathrm{Id}_\mathbf{G}$ necessarily $\mathbf{m} \neq \mathbf{0}$. Since $\mathbf{s} = \sigma^p(\mathbf{s}) = \mathbf{s} + p\mathbf{m}$ we have shown
\begin{equation}\label{EqMPM}
\mathbf{m} \neq \mathbf{0} \;\; \mbox{and} \;\; p\mathbf{m} = \mathbf{0}.
\end{equation}

Let $\tau_0 = \tau |_{\mathbf{G}_0}$. Since $\tau(\mathbf{G}_0) = \mathbf{G}_0$, by (\ref{EqTauSMinusMG0}) again
\begin{equation}\label{EqTau0AutG0}
\tau_0 \in \mathrm{Aut}_\mathrm{Group}(\mathbf{G}_0).
\end{equation}
Let $0 \leq i \leq p-1$. Since $\tau(\mathbf{G}_i) = \mathbf{G}_{\eta(i)}$, which follows by (\ref{EqEtaTauGI}), we have
\begin{equation}\label{EqTauISEtaISN}
\tau(i\mathbf{s}) = \eta(i)\mathbf{s} + \mathbf{n}_i
\end{equation}
for some $\mathbf{n}_i \in \mathbf{G}_0$. Using (\ref{EqTauSMinusMG0}) again, we deduce $\tau(\mathbf{0}) = \mathbf{0}$ and thus
\begin{equation}\label{EqEta00}
\mathbf{n}_0 = \mathbf{0},
\end{equation}
and also since $\sigma$ is a homomorphism
\begin{equation}\label{EqTauAndSigmaISX}
\tau(i\mathbf{s} + \mathbf{x}) = \eta(i)\mathbf{s} + \mathbf{n}_i + \tau_0(\mathbf{x}) \;\; \mbox{and} \;\; \sigma(i\mathbf{s} + \mathbf{x}) = i\mathbf{s} + i\mathbf{m} + \mathbf{x}
\end{equation}
for all $0 \leq i \leq p-1$ and $\mathbf{x} \in \mathbf{G}_0$. The last two equations describe $\tau$ and $\sigma$ explicitly.

Since $\tau$ and $\sigma$ commute, (\ref{EqTauAndSigmaISX}) implies $\tau_0(i\mathbf{m}) = \eta(i)\mathbf{m}$ for all $0 \leq i \leq p-1$. Thus by virtue of (\ref{EqTau0AutG0}) we have $\eta(i)\mathbf{m} = \tau_0(i\mathbf{m}) = i\tau_0(\mathbf{m}) = i\eta(1)\mathbf{m}$. Now $\mathbf{m} \neq \mathbf{0}$ and has order $p$ by (\ref{EqSigmaSSPlusM}). We have shown
\begin{equation}\label{EqEtaIModP}
\eta(i) \equiv i\eta(1) \,(\mathrm{mod} \; p) \;\; \mbox{and} \;\; \tau_0(i\mathbf{m}) = i\eta(1)\mathbf{m}
\end{equation}
for $0 \leq i \leq p-1$.

We now turn our attention to (F.d). Let $\mathbf{x}, \mathbf{x}' \in \mathbf{G}$. Then
\begin{equation}\label{EqFdXXPrime}
\tau(\mathbf{x} - \mathbf{x}') = \tau(\mathbf{x}) - \tau(\sigma^{\ell_\alpha(\mathbf{x}, \mathbf{x}')}(\mathbf{x}')).
\end{equation}
We first determine $\ell_\alpha(\mathbf{x}, \mathbf{x}')$. Now $\mathbf{x} \in \mathbf{G}_i$ and $\mathbf{x}' \in \mathbf{G}_{i'}$, where $0 \leq i, i' \leq p-1$. Since $\ell_\alpha(\mathbf{x}, \mathbf{x}') = 0$ when $i' = 0$ we will assume $1 \leq i'$. Hence $|\mathbf{x}'| = p$. Therefore the equation
$$
\displaystyle{\left(\alpha(\mathbf{x})/\alpha(\mathbf{x}')\right)^{p/|\mathbf{x}'|} = \lambda_{|\mathbf{x}'|}^{\ell_\alpha(\mathbf{x}, \mathbf{x}')}}
$$
is $\lambda^{i - i'} = \lambda^{\ell_\alpha(\mathbf{x}, \mathbf{x}')}$. We have therefore
\begin{equation}\label{EqEllXXPrime}
\ell_\alpha(\mathbf{x}, \mathbf{x}') \equiv i - i' \,(\mathrm{mod} \; p) \;\; \mbox{and} \;\; 0 \leq \ell_\alpha(\mathbf{x}, \mathbf{x}') \leq p-1
\end{equation}
for all $\mathbf{x} \in \mathbf{G}_i$, $\mathbf{x}' \in \mathbf{G}_{i'}$ and $0 \leq i, i' \leq p-1$, $1 \leq i'$. Since $\sigma^j(i'\mathbf{s}) = i'\mathbf{s} + ji'\mathbf{m}$ for all $j \in \mathbb{Z}$, when $\mathbf{x} = i\mathbf{s}$ and $\mathbf{x}' = i\mathbf{s}'$, (\ref{EqFdXXPrime}) and (\ref{EqEllXXPrime}) imply
\begin{eqnarray*}
\tau((i - i')\mathbf{s}) & = & \tau(i\mathbf{s} - i'\mathbf{s}) \\
& = & \tau(i\mathbf{s}) - \tau(\sigma^{\ell_\alpha(i\mathbf{s}, i'\mathbf{s}')}(i'\mathbf{s})) \\
& = & \tau(i\mathbf{s}) - \tau(\sigma^{i - i'}(i'\mathbf{s})) \\
& = & \tau(i\mathbf{s}) - \tau(i'\mathbf{s} + (i - i')i'\mathbf{m}).
\end{eqnarray*}
Therefore
\begin{equation}\label{EqTauMain}
\tau((i - i')\mathbf{s}) = (\eta(i) - \eta(i'))\mathbf{s} + (\mathbf{n}_i - \mathbf{n}_{i'}) - (i - i')i'\eta(1)\mathbf{m}
\end{equation}
for all $0 \leq i, i' \leq p-1$ and $1 \leq i'$ by (\ref{EqTauAndSigmaISX}) and (\ref{EqEtaIModP}).

We consider what it means for (\ref{EqTauMain}) to hold. It holds when $i = i'$. There are two other cases.
\medskip

\noindent
Case 1: $1 \leq i' < i \leq p-1$.

In this situation (\ref{EqTauMain}) is
\begin{equation}\label{EqCAse1}
\eta(i - i')\mathbf{s} + \mathbf{n}_{i - i'} = (\eta(i) - \eta(i'))\mathbf{s} + (\mathbf{n}_i - \mathbf{n}_{i'}) - (i - i')i'\eta(1)\mathbf{m}.
\end{equation}
This equation implies $\eta(i - i') - (\eta(i) - \eta(i')) \equiv 0 \,(\mathrm{mod} \; p)$ which is already implied by (\ref{EqEtaIModP}). The restriction we wish to place on $\tau$ is
\begin{enumerate}
\item[{\rm (ER)}] $\eta = \mathrm{Id}_{\mathbb{Z}_p}$,
\end{enumerate}
that is $\tau(\mathbf{G}_i) = \mathbf{G}_i$ for all $0 \leq i \leq p-1$. Observe that (ER) holds when $p = 2$. Also $\tau(\mathbf{G}_1) = \mathbf{G}_1$ implies (ER) by (\ref{EqEtaIModP}).

Under the assumption that (ER) holds (\ref{EqCAse1}) is
\begin{equation}\label{EqWithSAOne}
\mathbf{n}_{i - i'} = (\mathbf{n}_i - \mathbf{n}_{i'}) - (i - i')i'\eta(1)\mathbf{m}
\end{equation}
for all $1 \leq i' < i \leq p-1$. With $i' = 1$ the preceding equation can be written
\begin{equation}\label{EqWithEROne}
\mathbf{n}_i = \mathbf{n}_{i-1} + \mathbf{n}_1 + (i-1)\eta(1)\mathbf{m}
\end{equation}
for all $2 \leq i \leq p-1$. We note here that (\ref{EqWithSAOne}) holds, and hence (\ref{EqWithEROne}) also holds, when merely $p\mathbf{s} = \mathbf{0}$. As a result of (\ref{EqWithEROne}) we have
\begin{equation}\label{EqNI}
\mathbf{n}_i = i\mathbf{n}_1 + \binr{i}{2}{}\eta(1)\mathbf{m}
\end{equation}
for $0 \leq i \leq p-1$. This equation is $\mathbf{n}_i = \mathbf{n}_i$ for $i = 0, 1$. For $2 \leq i$ the formula follows by induction on $i$.

Assume (\ref{EqNI}) holds. Then (\ref{EqWithSAOne}) holds by virtue of the identity
$$
\displaystyle{\binr{i-i'}{2}{} =  \binr{i}{2}{} - \binr{i'}{2}{} - (i - i')i'}
$$
for all $0 \leq i' \leq i$.
\medskip

\noindent
Case 2: $0 \leq i < i' \leq p-1$.

In this case $i - i' = -p + r$, where $1 \leq r \leq p-1$. Since $p\mathbf{s} \in \mathbf{G}_0$, the left hand side of (\ref{EqTauMain}) is
\begin{eqnarray*}
\tau((i - i')\mathbf{s}) & = & \tau(-p\mathbf{s} + r\mathbf{s}) \\
& = & \tau(-p\mathbf{s}) + \tau(r\mathbf{s}) \\
& = & \tau(-p\mathbf{s}) + \tau((p + i - i')\mathbf{s}) \\
& = & \tau(-p\mathbf{s}) + \eta(p + i - i')\mathbf{s} + \mathbf{n}_{p + i - i'}.
\end{eqnarray*}
Therefore (\ref{EqTauMain}) becomes, with the restriction (ER) holding,
\begin{eqnarray}\label{EqCase2}
\lefteqn{\tau(-p\mathbf{s}) + p\mathbf{s} + p\mathbf{n}_1 + \binr{p + i - i'}{2}{}\mathbf{m}} \\
 & \qquad  = &  \left(\binr{i}{2}{} - \binr{i'}{2}{} - (i - i')i'\right)\mathbf{m}. \nonumber
\end{eqnarray}
Since
$$
\binr{p + i - i'}{2}{} = \binr{p}{2}{} + \binr{i}{2}{} - \binr{i'}{2}{} - (i' - p)(i - i')
$$
for all $0 \leq i < i' \leq p-1$, and $p\mathbf{m} = \mathbf{0}$, in Case 2, under the assumption (ER) holds, (\ref{EqTauMain}) is
\begin{equation}\label{EqCase2}
\tau_0(p\mathbf{s}) = p\mathbf{s} + p\mathbf{n}_1 + \binr{p}{2}{}\mathbf{m}.
\end{equation}
When $p = 2$ this equation is $\tau(2\mathbf{s}) = 2\mathbf{s} + 2\mathbf{n}_1 + \mathbf{m}$ and when $p > 2$ this equation is $\tau(p\mathbf{s}) = p\mathbf{s} + p\mathbf{n}_1$.
\begin{Prop}\label{PropMainExample}
Let $p$ be a prime integer and suppose $k$ contains a primitive $p^{th}$ root of unity $\lambda$. Let $\mathbf{G}$ be a finite abelian group with a subgroup $\mathbf{G}_0$ of index $p$. Let $\mathbf{s} \in \mathbf{G}$ satisfy $\mathbf{s} + \mathbf{G}$ generates $\mathbf{G}/\mathbf{G}_0$.
\begin{enumerate}
\item[{\rm (a)}] Suppose $\mathbf{m} \in \mathbf{G}_0$ satisfies $\mathbf{m} \neq \mathbf{0}$ and $p\mathbf{m} = \mathbf{0}$. Then there exists $\sigma \in \mathrm{Aut}_\mathrm{Group}(\mathrm{G})$ determined by $\sigma(\mathbf{s}) = \mathbf{s} + \mathbf{m}$ and $\sigma(\mathbf{x}) = \mathbf{x}$ for all $\mathbf{x} \in \mathbf{G}_0$. In particular $\sigma$ has order $p$.
\item[{\rm (b)}] Let $\mathbf{m}$ be as in part (a). Suppose $\tau_0 \in \mathrm{Aut}_\mathrm{Group}(\mathbf{G}_0)$ satisfies $\tau_0(\mathbf{m}) = \mathbf{m}$ and there is an $\mathbf{n} \in \mathbf{G}_0$ such that $\tau_0(2\mathbf{s}) = 2\mathbf{s} + 2\mathbf{n} + \mathbf{m}$, if $p = 2$, and $\tau_0(p\mathbf{s}) = p\mathbf{s} + p\mathbf{n}$, if $p > 2$. Then $(\tau, \alpha) \in \mathcal{F}_{\mathbf{G}, \,\sigma}$, where $\tau$ is defined by
   $$
   \tau(i\mathbf{s} + \mathbf{x}) = i(\mathbf{s} + \mathbf{n}) + \binr{i}{2}{}\mathbf{m} + \tau_0(\mathbf{x})
   $$
for all $0 \leq i \leq p-1$ and $\mathbf{x} \in \mathbf{G}_0$, and the character $\alpha$ is given by $\alpha(\mathbf{G}_i) = \{ \lambda^i\}$ for all $0 \leq i \leq p-1$. Furthermore $\tau \not \in \mathrm{Aut}_\sigma(\mathbf{G})$.
\item[{\rm (c)}] $\mathrm{Aut}_\sigma(\mathbf{G}) < \Gamma (\mathbf{G}, \lambda, \sigma)$.
\end{enumerate}
\end{Prop}

\pf
The reader is left with the straightforward exercise of completing the proof of the proposition.
\qed
\medskip

The conditions of Proposition \ref{PropMainExample} are met in the following example. Let $\mathbf{G} = \mathbb{Z}_p \oplus \cdots \oplus \mathbb{Z}_p \oplus \mathbb{Z}_{p^2}$ where the number of $\mathbb{Z}_p$ summands is at least one. Let $\mathbf{s}$ be a generator of $\mathbb{Z}_{p^2}$. Then $\mathbf{G}_0 = \mathbb{Z}_p \oplus \cdots \oplus \mathbb{Z}_p \oplus (p\mathbf{s})$ is a subgroup of $\mathbf{G}$ of index $p$ and is a direct sum of copies of $\mathbb{Z}_p$. Observe that $p\mathbf{x} = \mathbf{0}$ for all $\mathbf{x} \in \mathbf{G}_0$.

Suppose $p > 2$. Take $\tau_0 = \mathrm{Id}_{\mathbf{G}_0}$ and choose any $\mathbf{m}, \mathbf{n} \in \mathbf{G}_0$, where $\mathbf{m} \neq \mathbf{0}$. Then $\tau_0(p\mathbf{s}) = p\mathbf{s} = p\mathbf{s} + p\mathbf{n}$.

Now suppose $p = 2$. Choose a $\tau_0 \in \mathrm{Aut}_\mathrm{Group}(\mathbf{G}_0)$ such that the $(\tau_0)$-orbit of $2\mathbf{s}$ has length $2$. The function which permutes the last two summands of $\mathbf{G}_0$ is an example. Let $\mathbf{m} = \tau_0(2\mathbf{s}) - 2\mathbf{s} $. Then $\tau_0(2\mathbf{s}) = 2\mathbf{s} + \mathbf{m}$; hence $\mathbf{m} \in \mathbf{G}_0$, $\mathbf{m} \neq \mathbf{0}$, and $2\mathbf{m} = \mathbf{0}$. For any $\mathbf{n} \in \mathbf{G}_0$ observe that  $\tau_0^2(2\mathbf{s}) = 2\mathbf{s} + 2\mathbf{n} + \mathbf{m}$.

When $p > 2$ the group $\mathbf{G} = \mathbb{Z}_p \oplus \cdots \oplus \mathbb{Z}_p$ which is the direct sum of $n$ copies $\mathbb{Z}_p$, where $n \geq 2$, provides examples also. Here we take $\mathbf{G}_0$ to be the direct sum of the first $n-1$ copies of $\mathbb{Z}_p$ and $\mathbf{s}$ to be a generator of the last copy. Choose any $\mathbf{m}, \mathbf{n} \in \mathbf{G}_0$, where $\mathbf{m} \neq \mathbf{0}$. Let $\tau_0 = \mathrm{Id}_{\mathbf{G}_0}$. Since $p\mathbf{x} = \mathbf{0}$ for all $\mathbf{x} \in \mathbf{G}$ the equation $\tau_0(p\mathbf{s}) = p\mathbf{s} + p\mathbf{n}$ is trivially satisfied. We point out that $p\mathbf{m} = \mathbf{0}$ in this example.
\section{The biproduct of Section \ref{SecGeneric} revisited}\label{SecGenericBiproduct}
Here we apply results of Sections \ref{SecRelated} and \ref{SecGammaGLambda} to the biproduct of Section \ref{SecGeneric}.
\begin{Theorem}\label{ThmMainOfPaper}
Let $A = k[\mathcal{G}] \times H$ be the biproduct of Theorem \ref{ThmFromProjNew}, where $\mathcal{G}$ is a finite abelian group. Suppose $k$ contains primitive $|\mathcal{G}|^{th}$ root of unity. Then the one-one group homomorphism
$$
\Phi_{\mathcal{YD}} : \mathcal{N}(B, H)^{op} \rtimes_\varphi \mathrm{Aut}_{\mathcal{YD}-\mathrm{Hopf}}(B) \longrightarrow \mathrm{Aut}_\mathrm{Hopf}(A, \pi)
$$
of Theorem \ref{ThmAutSemiProd} is an isomorphism if any of the following hold:
\begin{enumerate}
\item[{\rm (a)}] $\mbox{\boldmath $\theta$}$ fixes a unique one-dimensional ideal of $k[\mathcal{G}]$.
\item[{\rm (b)}] The orders of $\mbox{\boldmath $\theta$}$ and $|\mathcal{G}|$ are relatively prime.
\item[{\rm (c)}] If $\mbox{\boldmath $\theta$}$ has order $p$, where $p$ is prime, and $p^2$ does not divide $|\mathcal{G}|$.
\end{enumerate}
\end{Theorem}

\pf
We have noted that $k$ contains a primitive $|\mbox{\boldmath $\theta$}|^{th}$ root of unity $\lambda$. Let $A' = k[\mathcal{G}]\times k[\mathbf{U}]$. At this point we bring into play the group homomorphism $\Phi : \mathrm{Aut}_\mathrm{Hopf}(A, \pi) \longrightarrow \mathrm{Aut}_\mathrm{Hopf}(A', \pi')$ of Section \ref{SecRelated} and the isomorphism $\Phi_{A', \mathbf{G}} : \mathrm{Aut}_{\mathrm{Hopf}}(A', \pi') \longrightarrow \Gamma (\mathbf{G}, \lambda, \sigma)$ of Theorem \ref{ThmAPrtimeBoldG}. The conditions (a)-(c) of Theorem \ref{ThmAutSemiProd} translate to the conditions (a)-(c) of Theorem \ref{ThmSigmaNoFix} respectively. Thus $\Gamma (\mathbf{G}, \lambda, \sigma) = \mathrm{Aut}_\sigma(\mathbf{G})$ by the same.

Let $F \in \mathrm{Aut}_\mathrm{Hopf}(A, \pi)$. Then $F_L = \Phi(F)_L$ is a coalgebra map by part (b) of Theorem \ref{ThmAPrtimeBoldG}. At this point the proof follows by Corollary \ref{CorAutSemiDirectProductCC}.
\qed
\medskip

Generally $\Phi_{\mathcal{YD}}$ is not an isomorphism.
\begin{Prop}\label{PropPhiNot}
Suppose $k$ is an algebraically closed field of characteristic $0$. Let $p$ be a prime integer. Then there exists a biproduct $A = k[\mathcal{G}] \times k[\mathbb{Z}_p]$, where $\mathcal{G}$ is a finite abelian group, such that:
\begin{enumerate}
\item[{\rm (a)}] The one-one group map $\Phi_{\mathcal{YD}}$ of Theorem \ref{ThmAutSemiProd} is not onto. In particular there is an $F \in \mathrm{Aut}_\mathrm{Hopf}(A, \pi)$ such that $F_L$ is not a coalgebra map.
\item[{\rm (b)}] The one-one group map $\Phi_\nu : \mathcal{N}(B, H)^{op} \rtimes_\varphi \mathrm{Im}(\nu) \longrightarrow \mathrm{Sym}(A)$ of Theorem \ref{ThmAutSemiProdHCoComm} satisfies $\mathrm{Im}(\Phi_\nu) \not \subseteq  \mathrm{Aut}_\mathrm{Hopf}(A, \pi)$.
\end{enumerate}
\end{Prop}

\pf
In light of Theorems \ref{ThmAutSemiProd}, \ref{ThmAutSemiProdHCoComm}, and \ref{ThmAPrtimeBoldG} we need only find finite abelian groups $\mathbf{G}$ such that $\mathrm{Aut}_\sigma(\mathbf{G}) < \Gamma(\mathbf{G}, \lambda, \sigma)$, where $\sigma$ has order $p$ and $\lambda \in k^\times$ is a primitive $p^{th}$ root of unity. Such examples are constructed after the proof of Proposition \ref{PropMainExample}.
\qed
\section{Another group associated with $\Gamma (\mathbf{G}, \lambda, \sigma)$}\label{SecAutThetaGGTheta}
Let $G$ be an additive group and $\sigma \in \mathrm{Aut}_\mathrm{Group}(G)$.  Let $\mathrm{Aut}_\sigma(G)$ be the set of all $\tau \in \mathrm{Aut}_\mathrm{Group}(G)$ which commute with $\sigma$. Let $\mathrm{Sym}_\sigma^-(G)$ be the set $\tau \in \mathrm{Sym}(G)$ which satisfy $\tau \circ \sigma = \sigma \circ \tau$ and
\begin{equation}\label{EqTauAOrbit}
\tau(a - \mathcal{O}) = \tau(a) - \tau(\mathcal{O})
\end{equation}
for all $a \in G$ and $(\sigma)$-orbits $\mathcal{O}$ of $G$. Observe that $\tau \in \mathrm{Sym}_\sigma^-(G)$ permutes the $(\sigma)$-orbits of $G$ since $\tau \circ \sigma = \sigma \circ \tau$. These four sets of permutations are groups under composition and their relationships are described by
$$
\mathrm{Aut}_\sigma(G) \leq \mathrm{Aut}_\mathrm{Group}(G) \leq \mathrm{Sym}(G) \; \mbox{and} \; \mathrm{Aut}_\sigma(G) \leq \mathrm{Sym}_\sigma^-(G) \leq \mathrm{Sym}(G).
$$
The group $\Gamma(\mathbf{G}, \lambda, \sigma)$ of Section \ref{SecGammaGLambda} satisfies
$$
\Gamma(\mathbf{G}, \lambda, \sigma) \leq \mathrm{Sym}_\sigma^-(\mathbf{G}).
$$
Since $\mathrm{Sym}_\sigma^-(G)$ is such a natural generalization of $\Gamma (\mathbf{G}, \lambda, \sigma)$ perhaps there will be some interest in it.

In this section we derive some of the elementary properties of the group $\mathrm{Sym}_\sigma^-(G)$ and describe examples such that $\Gamma(\mathbf{G}, \lambda, \sigma) < \mathrm{Sym}_\sigma^-(\mathbf{G})$. Let $a \in G$. We let $|a|$ denote the cardinality of the $(\sigma)$-orbit of $G$ which $a$ generates.

The set $F$ of fixed points of $\sigma$ is a subgroup of $G$. Let $\tau \in \mathrm{Sym}_\sigma^-(G)$. Then $\tau(F) = F$ since $\tau$ and $\sigma$ commute and $\tau$ is bijective. Let $m \in F$ and $\tau \in \mathrm{Sym}_\sigma^-(G)$. Observe that $\tau(a - m) = \tau(a) - \tau(m)$ for all $a \in G$ by (\ref{EqTauAOrbit}) since $|m| = 1$. In particular $\tau(0) = \tau(0 - 0) = \tau(0) - \tau(0) = 0$ and thus $\tau(-m) = \tau(0 - m) = \tau(0) - \tau(m) = -\tau(m)$. As a consequence
\begin{equation}\label{EqTauAPlusM}
\tau(a + m) = \tau(a) + \tau(m)
\end{equation}
for all $a \in G$ and $m \in F$. Since $\tau(F) = F$ the last equation implies that $\tau$ permutes the left cosets of $F$ in $G$ and
\begin{equation}\label{EqTau0AutSigma0F}
\tau_0 \in \mathrm{Aut}_{\sigma_0}(F),
\end{equation}

For the remainder of this section $\mathbf{G}$, $\mathbf{G}_0$, $\mathbf{s}$, and $\mathbf{m}$ are as described in Proposition \ref{PropMainExample} and $G = \mathbf{G}$. Then $\mathbf{G}_i := i\mathbf{s} + \mathbf{G}_0$, $0 \leq i \leq p-1$, lists the cosets of $\mathbf{G}_0$ in $\mathbf{G}$. Recall that $\mathbf{m} \in \mathbf{G}_0$, $\mathbf{m} \neq \mathbf{0}$, and $p\mathbf{m} = \mathbf{0}$. Let $\sigma \in \mathrm{Aut}_\mathrm{Group}(\mathbf{G})$ be the automorphism determined by $\sigma(\mathbf{s}) = \mathbf{s} + \mathbf{m}$ and $\sigma(\mathbf{x}) = \mathbf{x}$ for $\mathbf{x} \in \mathbf{G}_0$. Then $F = \mathbf{G}_0$.

Since $\tau$ permutes the cosets of $F = \mathbf{G}_0$ in $\mathbf{G}$, and $\tau(\mathbf{G}_0) = \mathbf{G}_0$, there is an $\eta \in \mathrm{Sym}(\mathbb{Z}_p)$ such that $\tau(\mathbf{G}_i) = \mathbf{G}_{\eta(i)}$ for all $0 \leq i \leq p-1$ and $\eta(0) = 0$. For $0 \leq i \leq p-1$ we have
\begin{equation}\label{EqTauISEta}
\tau(i\mathbf{s}) = \eta(i)\mathbf{s} + \mathbf{n}_i
\end{equation}
for all $0 \leq i \leq p-1$, where $\mathbf{n}_i \in \mathbf{G}_0$. We have shown $\tau(\mathbf{0}) = \mathbf{0}$. Therefore $\mathbf{n}_0 = \mathbf{0}$.  Thus (\ref{EqTauAPlusM}) and the fact that $\sigma$ is a homomorphism imply that (\ref{EqTauAndSigmaISX}) holds, in particular
$$
\tau (i\mathbf{s} + \mathbf{x}) = \eta(i)\mathbf{s} + \mathbf{n}_i + \tau_0(\mathbf{x})
$$
for all $0 \leq i \leq p-1$ and $\mathbf{x} \in \mathbf{G}_0$. Since $\sigma$ and $\tau$ commute (\ref{EqEtaIModP}) holds.

We now turn our attention to (\ref{EqTauAOrbit}). Let $\mathbf{a}, \mathbf{a}' \in \mathbf{G}$. Then $\mathbf{a} \in \mathbf{G}_i$ and $\mathbf{a}' \in \mathbf{G}_{i}'$ for some $0 \leq i, i' \leq p-1$. Write $\mathbf{a} = i\mathbf{s} + \mathbf{x}$ and $\mathbf{a}' = i'\mathbf{s} + \mathbf{x}'$, where $\mathbf{x}, \mathbf{x}' \in \mathbf{G}_0$. Since $p\mathbf{m} = \mathbf{0}$, the $(\sigma)$-orbit of $\mathbf{G}$ which $\mathbf{a}'$ generates is
$$
\mathcal{O} = \{i'\mathbf{s} + \ell i' \mathbf{m} + \mathbf{x}' \, | \, 0 \leq \ell \leq p-1\}.
$$
Therefore (\ref{EqTauAOrbit}) holds if and only if
\begin{eqnarray*}
\lefteqn{\{\tau((i - i')\mathbf{s}) - \ell i'\eta(1)\mathbf{m} \, | \, 0 \leq \ell \leq p-1\}} \\
& = & \{(\eta(i) - \eta(i'))\mathbf{s} + (\mathbf{n}_i - \mathbf{n}_{i'}) - \ell' i'\eta(1)\mathbf{m} \, | \, 0 \leq \ell' \leq p-1\}
\end{eqnarray*}
which holds if and only if
\begin{equation}\label{EqTauLIIPrime}
\tau((i - i')\mathbf{s}) = (\eta(i) - \eta(i'))\mathbf{s} + (\mathbf{n}_i - \mathbf{n}_{i'}) - \ell_{i, i'}i'\eta(1)\mathbf{m}
\end{equation}
for some $0 \leq  \ell_{i, i'} \leq p-1$. Observe that (\ref{EqTauLIIPrime}) holds when $i = i'$, with $\ell_{i, i} = 0$, and when $i' = 0$. We will therefore assume $1 \leq i'$.

From this point on we assume that (ER) holds; that is $\eta(i) = i$ for all $0 \leq i \leq p-1$. We note as a consequence of (\ref{EqEtaIModP}) that $\tau(\mathbf{m}) = \mathbf{m}$. There are two more cases to consider.
\medskip

\noindent
Case 1: $1 \leq i' < i \leq p-1$.

Since (ER) holds, in this situation (\ref{EqTauLIIPrime}) boils down to
$$
\mathbf{n}_{i - i'} = \mathbf{n}_i - \mathbf{n}_{i'} - \ell_{i, i'}i'\mathbf{m}
$$
in this case. When $i' = 1$ the preceding equation can be written
\begin{equation}\label{EqNNIMinus1Ell}
\mathbf{n}_i = \mathbf{n}_{i - 1} + \mathbf{n}_1 + \ell_{i, 1}\mathbf{m}
\end{equation}
for all $1 \leq i \leq p-1$. Thus
$$
\mathbf{n}_i = i\mathbf{n}_1 + \left(\sum_{u = 2}^i \ell_{u, 1} \right)\mathbf{m}
$$
for all $0 \leq i \leq p-1$. For $1 \leq j < i \leq p-1$ we set
$$
s(i, j) = \sum_{u = j + 1}^i \ell_{u, 1}
$$
and set $s(i, j) = 0$ otherwise. We have shown:
\begin{equation}\label{EqNIINSI1}
\mathbf{n}_i = i\mathbf{n}_1 + s(i, 1)\mathbf{m}
\end{equation}
for all $0 \leq i \leq p-1$. Therefore (\ref{EqTauLIIPrime}) can be rewritten
\begin{equation}\label{EqTauLIIPrimeCase1}
i'\ell_{i, i'}\mathbf{m} = (s(i, i - i') - s(i', 1))\mathbf{m}
\end{equation}
for all $1 \leq i' < i \leq p-1$. Since $i'$ is a unit $(\mathrm{mod} \; p)$ and $\mathbf{m}$ has order $p$ there exists a unique solution $  0 \leq \ell_{i, i'} \leq p-1$ to this equation in terms of the $s(i, j)$'s.
\medskip

\noindent
Case 2: $0 \leq i < i' \leq p-1$.

In this case $i - i' = -p + r$, where $1 \leq r \leq p-1$. Since $p\mathbf{s} \in \mathbf{G}_0$ and (ER) holds, in this situation the left hand side of (\ref{EqTauLIIPrime}) is
$$
\tau((i -i')\mathbf{s}) = \tau(-p\mathbf{s} + r\mathbf{s}) = \tau(-p\mathbf{s}) + r\mathbf{s} + \mathbf{n}_r
$$
and therefore (\ref{EqTauLIIPrime}) boils down to
\begin{equation}\label{EqCase2Sum}
\tau_0(p\mathbf{s}) - p\mathbf{s} - p\mathbf{n}_1 = (s(p + i - i', 1) + s(i', i))\mathbf{m} + \ell_{i, i'}i'\mathbf{m}.
\end{equation}
This equation has a solution $0 \leq \ell_{i, i'} \leq p-1$ if and only if
\begin{equation}\label{EqGroupMain}
\tau_0(p\mathbf{s}) - p\mathbf{s} - p\mathbf{n}_1 \in \mathbb{Z}\mathbf{m}.
\end{equation}
Again, such a solution must be unique.
\begin{Prop}\label{PropGroupMain}
Assume the hypothesis of Proposition \ref{PropMainExample}. Suppose $\mathbf{m} \in \mathbf{G}_0$ satisfies $\mathbf{m} \neq \mathbf{0}$ and $p\mathbf{m} = \mathbf{0}$. Let $\sigma \in \mathrm{Aut}_\mathrm{Group}(\mathbf{G})$ be determined by $\sigma(\mathbf{s}) = \mathbf{s} + \mathbf{m}$ and $\sigma(\mathbf{x}) = \mathbf{x}$ for all $\mathbf{x} \in \mathbf{G}_0$.
\begin{enumerate}
\item[{\rm (a)}] Suppose $\tau_0 \in \mathrm{Aut}_\mathrm{Group}(\mathbf{G}_0)$ satisfies $\tau_0(\mathbf{m}) = \mathbf{m}$ and $\mathbf{n} \in \mathbf{G}_0$ is such that $\tau_0(p\mathbf{s}) - p\mathbf{s} - p\mathbf{n} \in \mathbb{Z}\mathbf{m}$. Let $0 \leq \ell_{2, 1}, \ldots, \ell_{p-1, 1} \leq p-1$. Then there is a $\tau \in \mathrm{Sym}_\sigma^-(\mathbf{G})$ given by
    $$
    \tau(i\mathbf{s} + \mathbf{x}) = i(\mathbf{s} + \mathbf{n}) + s(i, 1)\mathbf{m} + \tau_0(\mathbf{x})
    $$
for all $0 \leq i \leq p-1$ and $\mathbf{x} \in \mathbf{G}_0$, where $s(i, 1) = \sum_{u = 2}^i \ell_{u, 1}$ for all $0 \leq i \leq p-1$.
\end{enumerate}
For the $\tau$ of part (a):
\begin{enumerate}
\item[{\rm (b)}] If $\tau \in \mathrm{Aut}_\sigma(\mathbf{G})$ then $\ell_{2, 1} = \cdots = \ell_{p-1, 1} = 0$.
\item[{\rm (c)}] If $\tau \in \Gamma(\mathbf{G}, \lambda, \sigma) \setminus \mathrm{Aut}_\sigma(\mathbf{G})$ then $\ell_{i, 1} = i-1$ for all $2 \leq i \leq p-1$.
\end{enumerate}
Furthermore
\begin{enumerate}
\item[{\rm (d)}] $\mathrm{Aut}_\sigma(\mathbf{G}) < \Gamma (\mathbf{G}, \lambda, \sigma) < \mathrm{Sym}_\sigma^-(\mathbf{G})$ when $p \geq 3$.
\end{enumerate}
\end{Prop}

\pf
Most of the details for the proof of part (a) can be gleaned from the discussion preceding the statement of the proposition. The reader is left with the short exercise of completing the proof. As for part (b), suppose that $\tau$ is a group homomorphism. Since $\tau(\mathbf{s}) = \mathbf{s} + \mathbf{n}$ it follows that $\tau(i\mathbf{s}) = i\mathbf{s} + i\mathbf{n}$ for all $0 \leq i \leq p-1$. Therefore $s(i, 1)\mathbf{m} = \mathbf{0}$ for all $0 \leq i \leq p-1$ which means $\ell_{i, 1} = 0$ for all $2 \leq i \leq p-1$.

As for part (c), suppose $\tau \in \Gamma(\mathbf{G}, \lambda, \sigma) \setminus \mathrm{Aut}_\sigma(\mathbf{G})$. As $\tau(\mathbf{G}_i) = \mathbf{G}_i$ it follows $\eta = \mathrm{Id}_{\mathbb{Z}_p}$. Note that $\mathbf{G}_0$ is the set of fixed points of $\sigma$. By (\ref{EqNI}) we have $s(i, 1)\mathbf{m} = \binr{i}{2}{}\mathbf{m}$ for all $2 \leq i \leq p-1$. Therefore $\ell_{i, 1} = i-1$ for all $2 \leq i \leq p-1$. We have established part (c). Part (d) follows from parts (b) and (c).
\qed
\medskip

The examples we constructed which realize Proposition \ref{PropMainExample} of course realize Proposition \ref{PropGroupMain}.
\section{Yet another example of $\mathrm{Sym}^-_\sigma(G)$}\label{SecAutThetaGGThetaMore}
Let $G$ be an additive group and $\sigma \in \mathrm{Aut}_\mathrm{Group}(G)$. Let $\mathrm{Sym}^+_\sigma(G)$ be the set of all $\tau \in \mathrm{Sym}(G)$ which satisfy $\tau \circ \sigma = \sigma \circ \tau$ and
\begin{equation}\label{EqTauAOrbit}
\tau(a + \mathcal{O}) = \tau(a) + \tau(\mathcal{O})
\end{equation}
for all $a \in G$ and $(\sigma)$-orbits $\mathcal{O}$ of $G$. Again, $\tau \in \mathrm{Sym}_\sigma^+(G)$ permutes the $(\sigma)$-orbits of $G$ since $\tau \circ \sigma = \sigma \circ \tau$. Thus $\mathrm{Sym}^+_\sigma(G)$ is a group under function composition and $\mathrm{Aut}_\sigma(G) \leq \mathrm{Sym}^+_\sigma(G) \leq \mathrm{Sym}(G)$. The groups $\mathrm{Sym}^-_\sigma(G)$ and $\mathrm{Sym}^+_\sigma(G)$ are the same by virtue of the following lemma.
\begin{Lemma}\label{LemmaSymSym}
Let $G$ be an additive group, $\sigma \in \mathrm{Aut}_\mathrm{Group}(G)$, and $\tau \in \mathrm{Sym}(G)$. Then the following are equivalent:
\begin{enumerate}
\item[{\rm (a)}] $\tau(a - \mathcal{O}) = \tau(a) - \tau(\mathcal{O})$ for all $a \in G$ and $(\sigma)$-orbits $\mathcal{O}$ of $G$.
\item[{\rm (b)}] $\tau(a + \mathcal{O}) = \tau(a) + \tau(\mathcal{O})$ for all $a \in G$ and $(\sigma)$-orbits $\mathcal{O}$ of $G$.
\end{enumerate}
\end{Lemma}

\pf If either part (a) or (b) holds then $\tau(0) = 0$. To see this take $a = 0$ and $\mathcal{O} = \{0\}$, where $0$ is the neutral element of $G$. Let $\mathcal{O}$ be a $(\sigma)$-orbit of $G$. Then $-\mathcal{O}$ is an orbit also since $\sigma$ is a homomorphism.

Suppose part (a) holds. Then $\tau(-\mathcal{O}) = -\tau(\mathcal{O})$ as $\tau(-\mathcal{O}) = \tau(0 - \mathcal{O}) = \tau(0) - \tau(\mathcal{O}) = -\tau(\mathcal{O})$. Hence $\tau(a + \mathcal{O}) =  \tau(a - (-\mathcal{O})) = \tau(a) - \tau(-\mathcal{O}) =  \tau(a) - (-\tau(\mathcal{O})) =  \tau(a) + \tau(\mathcal{O})$
for all $a \in G$. We have shown part (a) implies part (b).

Suppose part (b) holds. Let $a \in \mathcal{O}$. Then $0 = \tau(0) \in \tau(a - \mathcal{O}) = \tau(a) + \tau(-\mathcal{O})$ which means $-\tau(a) \in \tau(-\mathcal{O})$. Therefore $-\tau(\mathcal{O}) \subseteq \tau(-\mathcal{O})$. As a result $\tau(-\mathcal{O}) = -(-\tau(-\mathcal{O})) \subseteq -\tau(-(-\mathcal{O})) = -\tau(\mathcal{O})$. We have shown $\tau(-\mathcal{O}) = -\tau(\mathcal{O})$ and consequently part (a) follows from part (b).
\qed
\medskip

Since $\mathrm{Sym}^-_\sigma(G)$ and $\mathrm{Sym}^+_\sigma(G)$ are the same group we use the notation $\mathrm{Sym}_\sigma(G)$ to represent it. To highlight the relationships:
$$
\mathrm{Sym}_\sigma(G) := \mathrm{Sym}^-_\sigma(G) = \mathrm{Sym}^+_\sigma(G).
$$

In the last section we constructed examples with $\mathrm{Aut}_\sigma(G) < \mathrm{Sym}^-_\sigma(G)$. In this section we construct examples where $\mathrm{Aut}_\sigma(G) < \mathrm{Sym}_\sigma(G)$ of a very different sort; these have three $(\sigma)$-orbit lengths instead of two. In the next section we will give sufficient conditions for $\mathrm{Aut}_\sigma(G) = \mathrm{Sym}_\sigma(G)$.

Let $X$ be a set and $A, B \subseteq X$. Recall that $A\setminus B = \{a \in A \, | \, a \not \in B\}$. When $B = \{ b\}$ we write $A \setminus b$ for $A \setminus \{ b\}$. Now suppose $f : X \longrightarrow Y$ is a map of sets. Then $f(A) \setminus f(B) \subseteq f(A \setminus B)$ and
\begin{equation}\label{EqfXYSetMinus}
f(A \setminus B) = f(A) \setminus f(B) \;\; \mbox{when $f$ is one-one.}
\end{equation}
We will use this fact several times in the construction of our examples.

Let $R$ be a finite commutative local ring with maximal ideal $\mathcal{M} = Rm$, where $m^2 = 0 \neq m$. The set of units of $R$ is $R^\times = R \setminus \mathcal{M}$. Since $R = R^\times \cup \mathcal{M}$ it follows that $\mathcal{M} = R^\times m \cup \{ 0 \}$. Let $a \in R$. Then $L_a : R \longrightarrow R$ defined by $L_a(x) = a + x$ for all $x \in R$ belongs to the group $\mathrm{Sym}(R)$ of permutations of $R$ under function composition. Thus $L_a(R^\times) = L_a(R \setminus \mathcal{M}) = L_a(R) \setminus L_a(\mathcal{M})$, or
\begin{equation}\label{EqAPlusTimes}
a + R^\times = R \setminus (a + \mathcal{M}),
\end{equation}
for all $a \in R$ by (\ref{EqfXYSetMinus}). As a consequence
\begin{equation}\label{EqAMTimes}
a + R^\times = R^\times
\end{equation}
for all $a \in \mathcal{M}$.

Let $\mathcal{O}_0 = \{ 0\}$, $\mathcal{O}_1 = R^\times$, and $\mathcal{O}_2 = R^\times m$. Then $R = \mathcal{O}_0 \cup \mathcal{O}_1 \cup \mathcal{O}_2$ and $\mathcal{O}_i \cap \mathcal{ O}_j = \emptyset$ unless $i = j$. Now suppose $\tau \in \mathrm{Sym}(G)$ satisfies $\tau(\mathcal{O}_i) = \mathcal{O}_i$ for $0 \leq i \leq 2$ and the restriction $\tau |_{\mathcal{O}_1} = \mathrm{Id}_{\mathcal{O}_1}$. For such a $\tau$ we have:
\begin{Lemma}\label{LemmaTauAOI}
$\tau (a + \mathcal{O}_i) = \tau(a) + \tau(\mathcal{O}_i)$ for all $ a \in R$ and $0 \leq i \leq 2$.
\end{Lemma}

\pf
Let $a \in R$. Then $a \in \mathcal{O}_j$ for some $0 \leq j \leq 2$. Since $\mathcal{O}_0 = \{ 0 \}$, $\tau (0) = 0$ and $\tau (a + \mathcal{O}_i) = \tau(a) + \tau(\mathcal{O}_i)$ when $j = 0$ or $i = 0$. Thus we may assume $1 \leq i, j \leq 2$.
\medskip

\noindent
\underline{Case 1}: $i = 1$.

We must show $\tau (a + R^\times) = \tau (a) + \tau(R^\times)$, that is $\tau (a + R^\times) = \tau (a) + R^\times$, for $a \in R^\times$ or $a \in R^\times m$. Since $\tau$ is one-one $\tau(a + R^\times) = R \setminus \tau (a + \mathcal{M})$ by (\ref{EqfXYSetMinus}) and (\ref{EqAPlusTimes}).

Suppose $a \in R^\times$. By assumption $\tau |_{\mathcal{O}_1} = \mathrm{Id}_{\mathcal{O}_1}$. Therefore $\tau(a) = a$; also $\tau (a + \mathcal{M}) = a + \mathcal{M}$ by (\ref{EqAMTimes}). As a consequence
$$
\tau(a + R^\times) = R \setminus \tau (a + \mathcal{M}) = R \setminus (a + \mathcal{M}) = a + R^\times = \tau(a) + \tau(R^\times)
$$
since $\tau(R^\times) = R^\times$.

Suppose $a \in R^\times m$. Then $a, \tau (a) \in \mathcal{M}$ which means $a + R^\times = R^\times = \tau(a) + R^\times$ by (\ref{EqAMTimes}). Hence $\tau(a + R^\times) = \tau(a) + \tau(R^\times)$.
\medskip

\noindent
\underline{Case 2}: $i = 2$.

We must show that $\tau(a + R^\times m) = \tau(a) + \tau(R^\times m)$, that is $\tau(a + R^\times m) = \tau(a) + R^\times m$, for $a \in R^\times$ or $a \in R^\times m$. First of all assume that $a \in R^\times$. Since $a + R^\times m \subseteq R^\times$ by (\ref{EqAMTimes}), and $\tau |_{\mathcal{O}_1} = \mathrm{Id}_{\mathcal{O}_1}$ by assumption, $\tau(a + R^\times m) = a + R^\times m = \tau(a) + R^\times m$.

Now suppose $a \in R^\times m$. Since $R^\times m = \mathcal{M} \setminus 0$ and $a \in \mathcal{M}$, by (\ref{EqfXYSetMinus}) we have $a + R^\times m = L_a(\mathcal{M}) \setminus L_a(0) = \mathcal{M} \setminus a$. Therefore $a + R^\times m = \mathcal{M} \setminus a$. As a result $\tau(a) + R^\times m = \mathcal{M} \setminus \tau(a) = \tau(a + R^\times m)$; the first equation follows from the preceding equation since $\tau(a) \in R^\times m$ and the second follows from the same by (\ref{EqfXYSetMinus}). Our proof of the lemma is complete.
\qed
\medskip

We next show that
\begin{equation}\label{EqOrbits}
|\mathcal{O}_1| > |\mathcal{O}_2| \geq |\mathcal{O}_0|,
\end{equation}
which boils down to $|\mathcal{O}_1| > |\mathcal{O}_2|$, or $|R^\times| > |R^\times m|$. Now $|R^\times| \geq |R^\times m|$. Suppose $|R^\times| = |R^\times m|$ and set $r = |R^\times m|$. Then $|\mathcal{M}| = r + 1$ and $|R| = |R^\times| + |\mathcal{M}| = 2r + 1$. Since $\mathcal{M}$ is a proper additive subgroup of $R$, $\ell|\mathcal{M}| = |R|$ for some $\ell \geq 2$. But then $\ell(r + 1) = 2r + 1$ which is not possible. We have established (\ref{EqOrbits}).

Now assume $R^\times$ is cyclic with generator $u$. Then $\sigma \in \mathrm{Aut}_{\mathrm{Group}}(R^+)$, where $R^+$ denotes the underlying additive group of $R$ and $\sigma$ is defined by $\sigma(a) = ua$ for all $a \in R$. For $a \in R$ let $\mathcal{O}_{\sigma, a}$ denote the $(\sigma)$-orbit of $G$ generated by $a$. The $(\sigma)$-orbits of $R$ are $\mathcal{O}_{\sigma, 1} = R^\times = \mathcal{O}_1$, $\mathcal{O}_{\sigma, m} = R^\times m = \mathcal{O}_2$ and $\mathcal{O}_{\sigma, 0} = \{ 0\} = \mathcal{ O}_0$. Let $\tau \in \mathrm{Sym}_\sigma (R^+)$. Since $\sigma  \circ  \tau = \tau  \circ  \sigma$ the $(\sigma)$-orbits of $R$ are permuted by $\tau$. As $\tau (0) = 0$, in light of (\ref{EqOrbits}) we see $\tau(\mathcal{O}_i) = \mathcal{O}_i$ for $0 \leq i \leq 2$. Also $\tau (ua) = u\tau(a)$ for all $a \in R$ since $\sigma$ and $\tau$ commute; thus $\tau(u^\ell a) = u^\ell \tau (a)$ for all $\ell \in \mathbb{Z}$ and $a \in R$.
\begin{equation}\label{EqTauOne}
\mbox{If $\tau \in \mathrm{Aut}_\sigma(R^+)$ and $\tau(1) = 1$ then $\tau = \mathrm{Id}_R$.}
\end{equation}

\pf
Assume the hypothesis of (\ref{EqTauOne}). Since $\tau(u^\ell) = u^\ell \tau(1) = u^\ell$ for all $\ell \in \mathbb{Z}$ and $R^\times$ is generated by $u$, $\tau (a) = a$ for all $a \in R^\times$. Let $a \in \mathcal{M}$. Then $1 + a \in R^\times$ by (\ref{EqAMTimes}). Therefore $1 + a = \tau (1 + a) = \tau(1) + \tau(a) = 1 + \tau(a)$ from which $\tau(a) = a$ follows. Since $R = R^\times \cup \mathcal{M}$ we have established (\ref{EqTauOne}).
\qed
\begin{equation}\label{EqAutThetaTheta}
\mathrm{Aut}_\sigma (R^+) = (\sigma).
\end{equation}

\pf
First of all $(\sigma) \subseteq \mathrm{Aut}_\sigma (R^+)$ since $\sigma \in \mathrm{Aut}_\sigma (R^+)$. Conversely, let $\tau \in \mathrm{Aut}_\sigma (R^+)$. Then $\tau (1) \in \tau(\mathcal{O}_1) = R^\times$ which means $\tau(1) = u^s$ for some $s \in \mathbb{Z}$ since $u$ generates $R^\times$. Since $\sigma^s(1) = u^s$ also, $\sigma^{-s} \circ \tau (1) = 1$, and therefore $\sigma^{-s} \circ \tau = \mathrm{Id}_R$ by (\ref{EqTauOne}). We have shown that $\tau = \sigma^s \in (\sigma)$. Consequently $\mathrm{Aut}_\sigma (R^+) \subseteq (\sigma)$.
\qed
\medskip

Let $r = |R^\times m|$. Then $|\mathcal{M}| = r + 1$. For $s \in \mathbb{Z}$ define $\tau_s \in \mathrm{Sym}(R)$ by $\tau_s(a) = a$ for $a \in R^\times \cup \{ 0 \}$ and $\tau_s(a) = \sigma^s(a)$ for $a \in R^\times m$. Then $\tau_s  \circ  \sigma = \sigma  \circ  \tau_s$, $\tau_s (\mathcal{O}_i) = \mathcal{O}_i$ for $0 \leq i \leq 2$, and $\tau_s|_{\mathcal{O}_1} = \mathrm{Id}_{\mathcal{O}_1}$. Using Lemma \ref{LemmaTauAOI} we conclude that $\tau_s \in \mathrm{Sym}_\sigma(R^+)$. Since $\tau_0 = \mathrm{Id}_R$ and $\tau_s \circ \tau_{s'} = \tau_{s + s'}$ for $s, s' \in \mathbb{Z}$ it follows that $G_2 = \{\tau_s \, | \, s \in \mathbb{Z}\}$ is a subgroup of $\mathrm{Sym}_\sigma(R^+)$. Noting that $r = |R^\times m| = |\mathcal{O}_{\sigma, m}|$ it is an easy exercise to show that the map $\mathbb{Z}_r \longrightarrow G_2$ given by $s \mapsto \tau_s$ is an isomorphism of groups and that $G_2 = \{\tau \in \mathrm{Sym}_\sigma(R^+) \, | \, \tau|_{R^\times} = \mathrm{Id}_{R^\times} \}$.
\begin{Prop}\label{PropMainExamples}
Let $R$ be a finite commutative local ring with maximal ideal $\mathcal{M} = Rm$, where $m^2 = 0 \neq m$. Suppose that $R^\times$ is cyclic with generator $u$ and let $\sigma \in \mathrm{Aut}_\mathrm{Group}(R^+)$ be defined by $\sigma(a) = ua$ for all $a \in R$. Then
$$
\mathrm{Sym}_\sigma (R^+) \simeq \mathrm{Aut}_\sigma (R^+) \times G_2 \simeq \mathbb{Z}_s \times \mathbb{Z}_r,
$$
where $s = |R^\times|$ and $r = |\mathcal{M}| - 1$.
\end{Prop}

\pf
Let $\tau \in \mathrm{Aut}_\sigma (R^+)$ and $\tau' \in G_2$. Then $\tau \circ \tau' \in \mathrm{Sym}_\sigma(R^+)$ so $f : \mathrm{Aut}_\sigma (R^+) \times G_2 \longrightarrow \mathrm{Sym}_\sigma (R^+)$ given by $f((\tau, \tau')) = \tau \circ \tau'$ is well-defined. Now $\tau = \sigma^s$ for some $s \in \mathbb{Z}$ by (\ref{EqAutThetaTheta}). Therefore $\tau \circ \tau' = \tau' \circ \tau$ which implies that $f$ is a group homomorphism. We will show that $f$ is an isomorphism.

Suppose $\tau \circ \tau' = \mathrm{Id}_R$. Then $1 = \tau(\tau'(1)) = \tau(1)$ which implies $\tau = \mathrm{Id}_R$ by (\ref{EqTauOne}). Therefore $\tau' = \mathrm{Id}_R$ as well. We have shown that $f$ is one-one.

To show that $f$ is onto let $\tau'' \in \mathrm{Sym}_\sigma(R^+)$. Then $\tau''(u^\ell) = u^\ell \tau''(1)$ for all $\ell \in \mathbb{Z}$. Now $\tau''(1) = u^s$ for some $s \in \mathbb{Z}$. Since $\sigma^s(u^\ell) = u^{s + \ell} = \tau''(1)u^\ell = \tau''(u^\ell)$ for all $\ell \in \mathbb{Z}$ we have that $(\sigma^{-s} \circ \tau'')|_{R^\times} = \mathrm{Id}_{R^\times}$. Therefore $\sigma^{-s} \circ \tau'' \in G_2$ and $f((\sigma^s, \sigma^{-s} \circ \tau'')) = \tau''$. We have shown that $f$ is onto. Thus $f$ is an isomorphism. The second isomorphism of the theorem follows from the first and the facts that $|(\sigma)| = |R^\times|$, (\ref{EqAutThetaTheta}), and $\mathbb{Z}_r \simeq G_2$; the latter was noted above.
\qed
\begin{Cor}\label{CorExamples}
Let $p$ be a positive prime and $G = \mathbb{Z}_{p^2}$ or $G = \mathbb{Z}_p \oplus \mathbb{Z}_p$. Then there exists a $\sigma \in \mathrm{Aut}_{\mathrm{Group}}(G)$ such that $\mathrm{Aut}_\sigma(G) \simeq \mathbb{Z}_{p(p-1)}$ and $\mathrm{Sym}_\sigma(G) \simeq \mathbb{Z}_{p(p-1)} \times \mathbb{Z}_{(p-1)}$. Hence $\mathrm{Aut}_\sigma(G) < \mathrm{Sym}_\sigma(G)$ when $p > 2$.
\end{Cor}

\pf
The ring $\mathbb{Z}_{p^2}$ has $p^2$ elements, is local, and has maximal ideal $\mathcal{M}$ with $p$ elements generated by $p$. The units of $\mathbb{Z}_{2^2}$ are easily seen to form a cyclic group. For odd primes and all $n \geq 1$ the units of $\mathbb{Z}_{p^n}$ form a cyclic group \cite[Theorem 26]{LED}. The hypothesis of Proposition \ref{PropMainExamples} is satisfied with $R = \mathbb{Z}_{p^2}$.

Let $F = \mathbb{Z}_p$ be the field with $p$ elements and consider the polynomial algebra $F[x]$. The quotient of $R = F[x]/(x^2)$ has $p^2$ elements, is a local ring, and identifying cosets with representatives has basis $\{ 1, x\}$. Its maximal ideal $\mathcal{M} = Fx$ has $p$ elements and $x^2 = 0 \neq x$. In particular $|R^\times| = p(p-1)$; indeed $R^\times = \{a + bx \, | \, a, b \in F, a \neq 0 \}$.  Since $F$ is a finite field $F^\times$ is cyclic. Let $v$ be a generator of $F^\times$. Then $v$ has order $p-1$. Observe that $1 + x$ is a unit of $R$ order $p$. Therefore $u = v(1 + x)$ is a unit of $R$ order $(p-1)p$ and consequently generates $R^\times$. Thus $R$ satisfies the hypothesis of Proposition \ref{PropMainExamples}. As a vector space $R = F \oplus F$ which means $R = \mathbb{Z}_p \oplus \mathbb{Z}_p$ as an abelian group.
\qed
\section{Sufficient conditions for $\mathrm{Aut}_\sigma (G) = \mathrm{Sym}_\sigma(G)$}\label{SecAutThetaGGThetaEqual}
Note that $\mathrm{Aut}_\sigma (G) = \mathrm{Sym}_\sigma(G)$ when $\sigma = \mathrm{Id}_G$. We consider involutions.
\begin{Theorem}\label{TheoremInvolution}
Suppose that $G$ is a finite abelian group of odd order and $\sigma \in \mathrm{Aut}_{\mathrm{Group}}(G)$ has order two. Then $\mathrm{Aut}_\sigma (G) = \mathrm{Sym}_\sigma(G)$.
\end{Theorem}

\pf
Let $\tau \in \mathrm{Sym}_\sigma(G)$. We need only show that $\tau$ is an endomorphism of $G$. Now let $a, b \in G$. Then $\tau(a + \mathcal{O}_{\sigma, b}) = \tau(a) + \tau(\mathcal{O}_{\sigma, b})$. Therefore if $\mathcal{O}_{\sigma, b}$ is a singleton set $\tau(a + b) = \tau(a) + \tau(b)$. Since $G$ is abelian the preceding equation holds if $\mathcal{O}_{\sigma, a}$ is a singleton set.

Suppose $\tau(a + b) \neq \tau(a) + \tau(b)$. Then neither $\mathcal{O}_{\sigma, a}$ nor $\mathcal{O}_{\sigma, b}$ are singleton sets and therefore each has two elements since $\sigma^2 = \mathrm{Id}_G$. From the equation $\tau(a + \mathcal{O}_{\sigma, b}) = \tau(a) + \tau(\mathcal{O}_{\sigma, b})$ we deduce the first two of four equations, namely: $\tau(a + b) = \tau(a) + \tau(\sigma(b))$ and $\tau(a + \sigma(b)) = \tau(a) + \tau(b)$. Since $G$ is commutative, with the roles of $a$ and $b$ reversed, these equations can be written $\tau(a + b) = \tau(\sigma(a)) + \tau(b)$ and $\tau(\sigma(a) + b) = \tau(a) + \tau(b)$. Using the first and third we deduce
$$
\tau(\sigma(a + b)) = \sigma(\tau(a + b)) = \sigma(\tau(a) + \tau(\sigma(b))) = \tau(\sigma(a)) + \tau(b) = \tau(a + b).
$$
Since $\tau$ is one-one $\sigma(a + b) = a + b$; therefore $a + b$ is a fixed point of $\sigma$. Using the second and fourth equations we deduce that $a + \sigma(b)$ is a fixed point of $\sigma$. Therefore $2a + b + \sigma(b)$ is a fixed point. Consequently $2a$ is a fixed point of $\sigma$ since $b + \sigma(b)$ is.

By assumption $(2m + 1)a = 0$ for some $m \geq 0$. Since $2a$ is a fixed point of $\sigma$ so is $-a = m(2a)$ and therefore $\sigma(a) = a$. But this means $\mathcal{O}_{\sigma, a} = \{ a \}$, a contradiction. We have shown that $\tau(a + b) = \tau(a) + \tau(b)$ after all.
\qed
\medskip

We next show in many cases the question of whether $\mathrm{Aut}_\sigma (G) = \mathrm{Sym}_\sigma(G)$ reduces to the case when $\sigma$ has one fixed point. Suppose that $G$ is a finite abelian group and $\sigma \in \mathrm{Aut}_{\mathrm{Group}}(G)$. Let $G_\sigma$ be the subgroup of fixed points of $\sigma$ and set $\overline{G} = G/G_\sigma$. Let $\overline{\sigma} \in \mathrm{Aut}_{\mathrm{Group}}(\overline{G})$ be defined by $\overline{\sigma}(a + G_\sigma) = \sigma(a) + G_\sigma$ for all $a \in G$.

Let $\tau \in \mathrm{Sym}_\sigma(G)$. Since $\tau$ permutes the $(\sigma)$-orbits of $G$ it permutes the fixed points of $\sigma$. Hence $\tau(G_\sigma) = G_\sigma$. We have noted $\tau (a + b) = \tau(a) + \tau(b)$ whenever $a$ or $b$ generates a singleton $(\sigma)$-orbit, that is whenever $a$ or $b$ is a fixed point of $\sigma$. Therefore $\overline{\tau} : \overline{G} \longrightarrow \overline{G}$ given by $\overline{\tau}(a + G_\sigma) = \tau(a) + G_\sigma$ for all $a \in G$ is well-defined. It is easy to see that $\overline{\tau} \in \mathrm{Sym}_{\overline{\sigma}}(\overline{G})$.
\begin{Theorem}\label{TheoremReduction}
Suppose $G$ is a finite abelian group and $\sigma \in \mathrm{Aut}_{\mathrm{Group}}(G)$. Let $\overline{G} = G/G_\sigma$ and suppose that the order of $\sigma$ and $|G|$ are relatively prime.
\begin{enumerate}
\item[{\rm (a)}]
$\overline{\sigma}$ has one fixed point in $\overline{G}$.
\item[{\rm (b)}] Let $\tau \in \mathrm{Sym}_\sigma(G)$. If $\overline{\tau} \in \mathrm{Aut}_{\overline{\sigma}}(\overline{G})$ then $\tau \in \mathrm{Aut}_\sigma(G)$.
\end{enumerate}
\end{Theorem}

\pf
We first show part (a). Suppose $a \in G$ and $\overline{\sigma}(a + G_\sigma) = a + G_\sigma$. Then $\sigma(a) - a \in G_\sigma$. Set $f = \sigma(a) - a$. Then $\sigma (a) = a + f$. Consequently $a$ and $a + f$ generate $\mathcal{O}_{\sigma, a}$. Let $|\mathcal{O}_{\sigma, a}| = r$. Expressing the sum of the elements of $\mathcal{O}_{\sigma, a}$ in two ways gives $\sum_{\ell = 0}^{r-1}\sigma^\ell(a) = \sum_{\ell = 0}^{r-1}\sigma^\ell (a + f) = \sum_{\ell = 0}^{r-1}\sigma^\ell(a) + rf$ which means $rf = 0$. Since $r$ divides the order of $\sigma$ it follows that $r$ and $|G|$ are relatively prime. Therefore $f = 0$ which means $\sigma (a) = a$. Thus $a \in G_\sigma$ which means $a + G_\sigma = G_\sigma$. Part (a) is established.

Assume the hypothesis of part (b). Let $a, b \in G$. Writing $\overline{x} = x + G_\sigma$ for $x \in G$ we have $\overline{\tau}(\overline{a}  + \overline{b}) = \overline{\tau}(\overline{a}) +  \overline{\tau}(\overline{a})$, or $\tau (a + b) - \tau(a) - \tau(b) \in G_\sigma$. Let $f \in G_\sigma$ be this difference. Thus $\tau (a + b) - \tau(a) = \tau(b) + f$. Now $\tau(a + b) = \tau(a) + \tau(\sigma^i(b))$ for some $0 \leq i < r$, where $r = |\mathcal{O}_{\sigma, b}|$. The two equations preceding the last imply $\tau(\sigma^i(b)) = \tau(b) + f$. Since $\tau(\sigma^i(b)) = \sigma^i(\tau(b))$ both $\tau(b)$ and $\tau(b) + f$ generate the same $(\sigma)$-orbit. In the proof of part (a) we established $f = 0$; hence $\tau(a + b) = \tau(a) + \tau(b)$. We have shown that $\tau \in \mathrm{Aut}_\sigma(G)$.
\qed

\end{document}